# Symmetric weak multicategories

Volodymyr Lyubashenko[*]

December 10, 2025

**Abstract**

A multicategory is what remains of a monoidal category when monoidal product is not available. A weak multicategory means that hom-sets are in fact categories, and in place of usual equations, there are natural isomorphisms, which have to satisfy their own equations. A symmetric weak multicategory implies a weak multicategory with a weak (up to a cocycle) action of symmetric groups. [1]

## 0.1 Conventions

We use symmetric strict cartesian category $\mathcal{S}et$ (monoidal with cartesian product as the monoidal product) equivalent to the ordinary category of sets $\widetilde{\mathcal{S}et}$ (which is cartesian, but not strict). Otherwise, we would have to mention the obvious bijection $(X \times Y) \times Z \to X \times (Y \times Z)$ explicitly. To an equality (or symmetry which we often denote by $=$) in $\mathcal{S}et$ there corresponds a canonical isomorphism in $\widetilde{\mathcal{S}et}$. Based on this $\mathcal{S}et$, we consider also the symmetric strict cartesian category $\mathcal{C}at$ of (essentially small) categories. We use also the symmetric strict monoidal category $\mathcal{S}_{\mathsf{sk}}$ skeletal for the category of finite totally ordered sets and arbitrary maps. Thus, $\operatorname{Ob}\mathcal{S}_{\mathsf{sk}} \cong \mathbb{N} = \{0, 1, 2, \dots\}$. The monoidal product of $I, J \in \operatorname{Ob}\mathcal{S}_{\mathsf{sk}}$ is their lexicographic product, which is the cartesian product of sets $I \times J$, equipped with the ordering $(i_1, j_1) <'' (i_2, j_2)$ iff $(j_1 < j_2$ or $(j_1 = j_2$ and $i_1 < i_2))$. We denote this monoidal product by $I \times J$ although it is not a cartesian one. The strict monoidal subcategory $\mathcal{O}_{\mathsf{sk}} \subset \mathcal{S}_{\mathsf{sk}}$ has the same objects, but only order-preserving maps.

## 0.2 Introduction

Unbiased form of symmetric multicategories was given in [BLM08, Definition 3.7], [Lyu23, Section 1.3] following Leinster [Lei03, Definition A.2.1]. When one passes from ordinary multicategories (in which morphisms form a set) to weak multicategories (whose morphisms form a category and structure maps are functors) one replaces equalities with natural isomorphisms. These isomorphisms, in their turn, have to satisfy some equations. In the search of these equations there is a guiding principle: each loop must be contractible. That is, any automorphism (of a functor constructed from compositions and insertion of a unit 1-morphism and other structural components) built from structure data must be equal to the identity automorphism. We implement this principle in the definitions of a symmetric weak multicategory (Definition 1.1.1), a symmetric weak multifunctor (Definition 1.1.6), a multinatural transformation (Definition 3.1.1) and a modification (Definition 3.1.3).

The usual definition of symmetric multicategories (or operads) uses explicit action of the symmetric group permuting the arguments. In Leinster's definition of symmetric multicategories explicit action of symmetric groups is not used. Instead he constructs this action in [Lei03, Lemma A.2.2] departing from the structure data of a symmetric multicategory. We generalize this action to symmetric weak multicategory so that a permutation acts by an equivalence of categories (Corollary 2.1.5), the identity permutation acts by an equivalence isomorphic to the identity functor (Remark 2.1.1), the composition of two permutations acts by an equivalence isomorphic to the composition of the actions of two permutations-factors in the opposite order (Corollary 2.1.3). Moreover, this isomorphism satisfies the non-abelian cocycle identity (Proposition 2.1.7). In this sense the action of a symmetric group by permutation of arguments on hom-categories of a symmetric weak multicategory is weak. This action of symmetries $(12)\colon I \times L \to L \times I \in \mathcal{S}et$ intertwines the compositions indexed by $\psi \times I$ and $I \times \psi$ (Proposition 2.2.1).

[*] *Email address:* lub@imath.kiev.ua

The general theory of lax monoids in non-symmetric case is provided by Day and Street [DS03]. It seems that replacing their category $\boldsymbol{\Delta}$ (algebraists' simplicial category coinciding with our $\mathcal{O}_{\mathsf{sk}}$) with the category $\mathcal{S}_{\mathsf{sk}}$ would yield the symmetric case. However we pursue simpler goal of symmetric weak multicategories.

Weak Cat-operads were also studied by Došen and Petrić [DP15]. They use a different system of operations adapted to non-symmetric case. One may guess that their notions are related to our non-symmetric weak multicategories, obtained by replacing the category $\mathcal{S}_{\mathsf{sk}}$ with $\mathcal{O}_{\mathsf{sk}}$ everywhere, but to establish the precise relation would require too much efforts.

**Acknowledgement.** I thank Paul Taylor for writing a package 'diagrams.sty', which I use a lot. I am grateful for clarifying discussions to Prof. Dr. Anna Beliakova and to Prof. Dr. Giovanni Felder. Excellent environment for work in the University of Zurich and of the Swiss Federal institute of Technology in Zurich were due to the generosity of the National Centre of Competence in Research SwissMAP of the Swiss National Science Foundation (grant number 205607), to whom I express my deep gratitude. The author thanks Prof. Dr. Christoph Schweigert and Faculty of Mathematics, Computer Science and Natural Sciences of the University of Hamburg for interest in the subject and excellent conditions of work. I thank Prof. Tom Leinster at the University of Edinburgh for fruitful discussions. The author wishes to thank the Isaac Newton Institute and London Mathematical Society for the financial support and the School of Mathematics at the University of Edinburgh for their hospitality. I am really grateful to the Armed Forces of Ukraine who gave me the possibility to work quietly on the subject.

# 1 Symmetric weak multicategories and multifunctors

## 1.1 Symmetric weak multicategories

**1.1.1 Definition.** Symmetric multicategory $\mathsf{C}$ weakly enriched in $\mathcal{C}at$ is

— a set of objects $\operatorname{Ob}\mathsf{C}$;

— a category $\mathsf{C}\big((X_i)_{i\in I};Y\big)$ for each tuple $\big((X_i)_{i\in I},Y\big)$ of objects of $\mathsf{C}$;

— a composition functor $\mu_\phi\colon \big[\prod_{j\in J}\mathsf{C}\big((X_i)_{i\in\phi^{-1}j};Y_j\big)\big] \times \mathsf{C}\big((Y_j)_{j\in J};Z\big) \to \mathsf{C}\big((X_i)_{i\in I};Z\big)$ for each map $\phi\colon I \to J$ from $\mathcal{S}_{\mathsf{sk}}$ and objects $X_i, Y_j, Z \in \operatorname{Ob}\mathsf{C}$, $i \in I$, $j \in J$;

— an identity 1-morphism $1_X \in \operatorname{Ob}\mathsf{C}(X;X)$ for each object $X \in \operatorname{Ob}\mathsf{C}$;

— a natural isomorphism (associator)

$$\big[\textstyle\prod_{j\in J}\mathsf{C}\big((X_i)_{i\in\phi^{-1}j};Y_j\big)\big] \times \big[\prod_{k\in K}\mathsf{C}\big((Y_j)_{j\in\psi^{-1}k};Z_k\big)\big] \times \mathsf{C}\big((Z_k)_{k\in K};W\big)$$

$$\wr\|$$

$$\big[\textstyle\prod_{k\in K}\big(\big[\prod_{j\in\psi^{-1}k}\mathsf{C}\big((X_i)_{i\in\phi_k^{-1}j};Y_j\big)\big] \times \mathsf{C}\big((Y_j)_{j\in\psi^{-1}k};Z_k\big)\big)\big] \times \mathsf{C}\big((Z_k)_{k\in K};W\big)$$

$$\Big\downarrow{\scriptstyle(\prod_{k\in K}\mu_{\phi_k})\times 1} \qquad\qquad\qquad\qquad \Big\downarrow{\scriptstyle 1\times\mu_\psi}$$

$$\Big[\prod_{k\in K}\mathsf{C}\big((X_i)_{i\in(\phi\psi)^{-1}k};Z_k\big)\Big] \times \mathsf{C}\big((Z_k)_{k\in K};W\big) \overset{\nu_{\phi,\psi}}{\Longrightarrow} \Big[\prod_{j\in J}\mathsf{C}\big((X_i)_{i\in\phi^{-1}j};Y_j\big)\Big] \times \mathsf{C}\big((Y_j)_{j\in J};W\big)$$

$$\searrow{\scriptstyle\mu_{\phi\cdot\psi}} \qquad \swarrow{\scriptstyle\mu_\phi}$$

$$\mathsf{C}\big((X_i)_{i\in I};W\big) \tag{1.1.1}$$

for each pair of composable maps $I \xrightarrow{\phi} J \xrightarrow{\psi} K \in \mathcal{S}_{\mathsf{sk}}$ and objects $X_i, Y_j, Z_k, W \in \operatorname{Ob}\mathsf{C}$, $i \in I$, $j \in J$, $k \in K$ (here $\phi_k = \phi|_{(\phi\psi)^{-1}(k)} \colon (\phi\psi)^{-1}(k) \to \psi^{-1}(k)$, $k \in K$);

— a natural isomorphism

$$\zeta_I \colon \big[\mathsf{C}((X_i)_{i\in I}; Z) \xrightarrow{1\times \mathrm{i}_Z} \mathsf{C}((X_i)_{i\in I}; Z) \times \mathsf{C}(Z;Z) \xrightarrow{\mu_{\triangledown\colon I\to\mathbf{1}}} \mathsf{C}((X_i)_{i\in I}; Z)\big] \to \mathrm{Id};$$

— a natural isomorphism

$$\theta_I \colon \big[\mathsf{C}((X_i)_{i\in I}; Z) \xrightarrow{(\prod_{i\in I}\mathrm{i}_{X_i})\times 1} \big(\prod_{i\in I}\mathsf{C}(X_i;X_i)\big) \times \mathsf{C}((X_i)_{i\in I}; Z) \xrightarrow{\mu_{\mathrm{id}_I}} \mathsf{C}((X_i)_{i\in I}; Z)\big] \to \mathrm{Id}$$

such that

— for each triple of composable maps $I \xrightarrow{\phi} J \xrightarrow{\psi} K \xrightarrow{\xi} L \in \mathcal{S}_{\mathsf{sk}}$ and objects $X_i$, $Y_j$, $Z_k$, $W_l$, $U \in \operatorname{Ob}\mathsf{C}$, $i \in I$, $j \in J$, $k \in K$, $l \in L$ we have a commutative cube shown on the following page. Here $\phi_l = \phi| \colon (\phi\psi\xi)^{-1}l \to (\psi\xi)^{-1}l$. Notice that

$$\big(\prod_{k\in K}\mu_{\phi_k}\big)\times 1_{\prod_{l\in L}\mathsf{C}((Z_k)_{k\in\xi^{-1}l};W_l)} \times 1 = \big(\prod_{l\in L}\prod_{k\in\xi^{-1}l}\mu_{\phi_k}\big)\times 1_{\prod_{l\in L}\mathsf{C}((Z_k)_{k\in\xi^{-1}l};W_l)} \times 1$$
$$= \prod_{l\in L}\big(\prod_{k\in\xi^{-1}l}\mu_{\phi_k}\times 1_{\mathsf{C}((Z_k)_{k\in\xi^{-1}l};W_l)}\big)\times 1.$$

A system of notations besides pasting is proposed by Bartlett [Bar14]. We apply it to the Cartesian 2-category $\mathcal{C}at$. The difference with his approach is that we read diagrams downwards. The source of each natural transformation is indicated by an oval box. The equation at Figure 1 takes the form

$$\prod_{j\in J}\mathsf{C}\big((X_i)_{i\in\phi^{-1}j};Y_j\big)\ \prod_{k\in K}\mathsf{C}\big((Y_j)_{j\in\psi^{-1}k};Z_k\big)\ \prod_{l\in L}\mathsf{C}\big((Z_k)_{k\in\xi^{-1}l};W_l\big)\ \mathsf{C}\big((W_l)_{l\in L};U\big) \qquad (1.1.2)$$

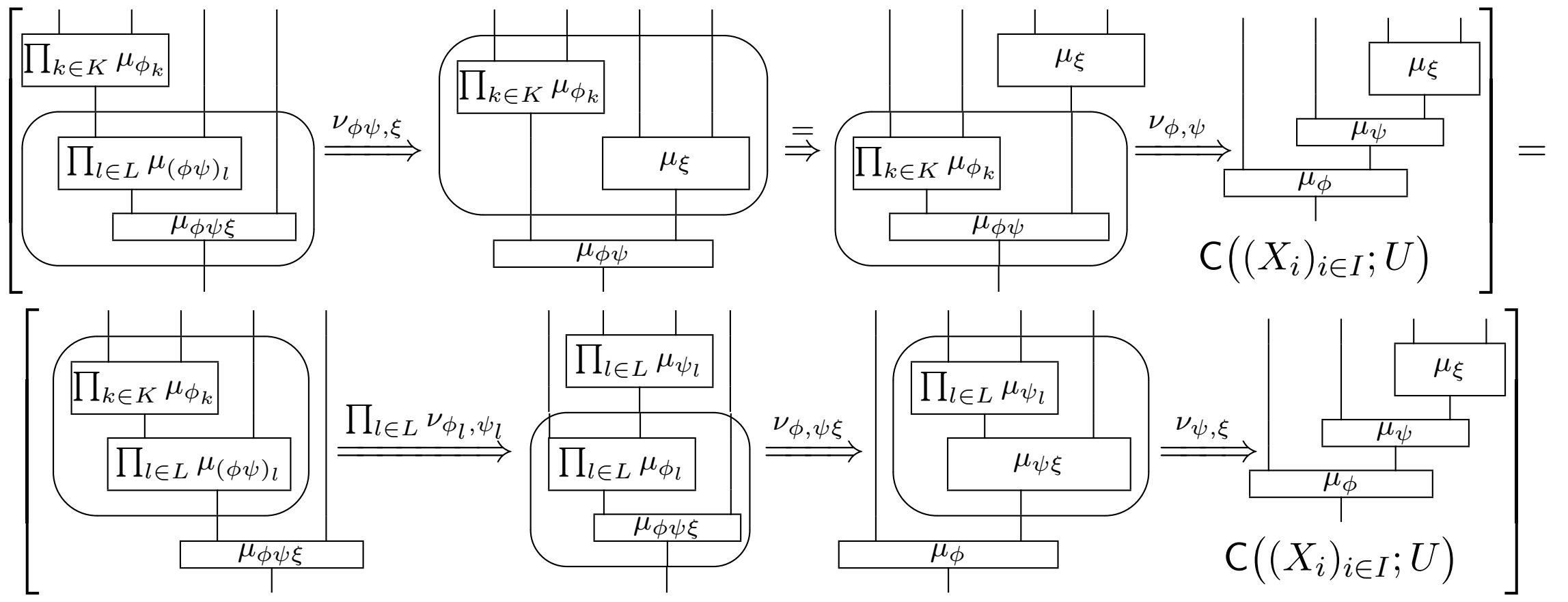

(here $\phi_l = \phi|_{(\phi\psi\xi)^{-1}(l)} \colon (\phi\psi\xi)^{-1}(l) \to (\psi\xi)^{-1}(l)$, $l \in L$).

— for maps $I \xrightarrow{\phi} J \xrightarrow{\mathrm{id}} J \in \mathcal{S}_{\mathsf{sk}}$ and objects $X_i, Y_j, Z \in \operatorname{Ob}\mathsf{C}$, $i \in I$, $j \in J$, the pasting

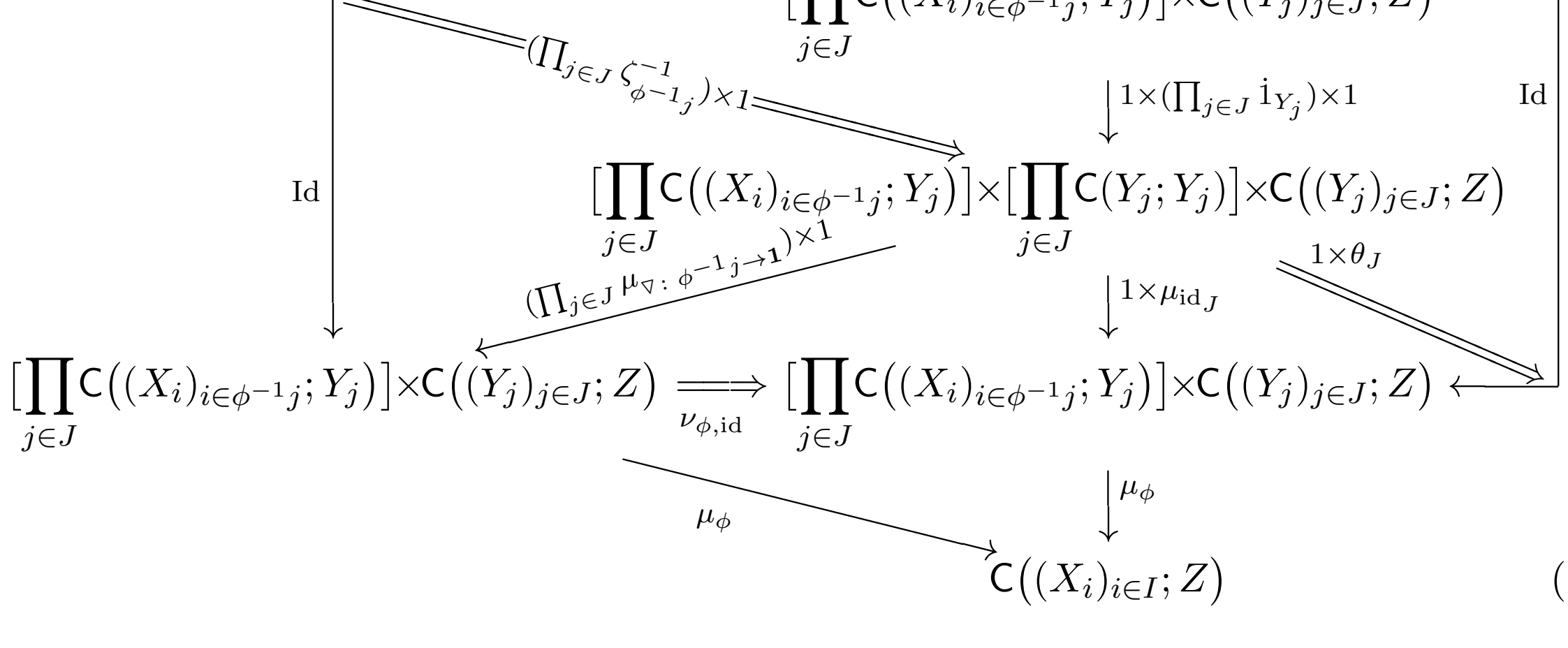

(1.1.3)

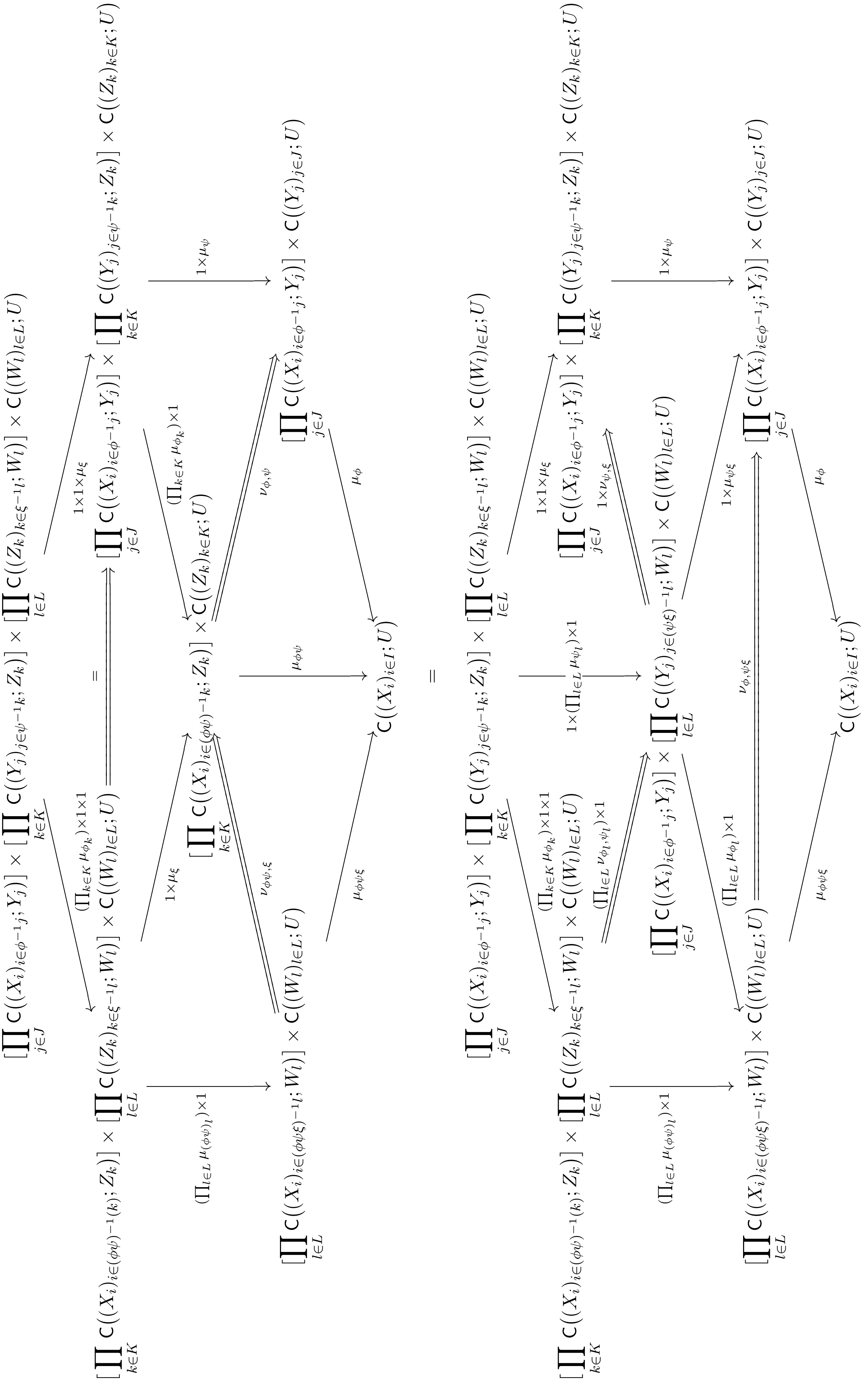

Figure 1: A pentagon equation for $\nu$

is equal to id: $\mu_\phi \to \mu_\phi \colon \bigl[\prod_{j\in J} \mathsf{C}\bigl((X_i)_{i\in\phi^{-1}j}; Y_j\bigr)\bigr] \times \mathsf{C}\bigl((Y_j)_{j\in J}; Z\bigr) \to \mathsf{C}\bigl((X_i)_{i\in I}; Z\bigr)$;

Plain weak multicategories are defined as above with $\mathcal{S}_{\mathsf{sk}}$ replaced with $\mathcal{O}_{\mathsf{sk}}$.

In Bartlett's notation [Bar14] the relation between $\zeta$, $\nu$ and $\theta$ takes the form:

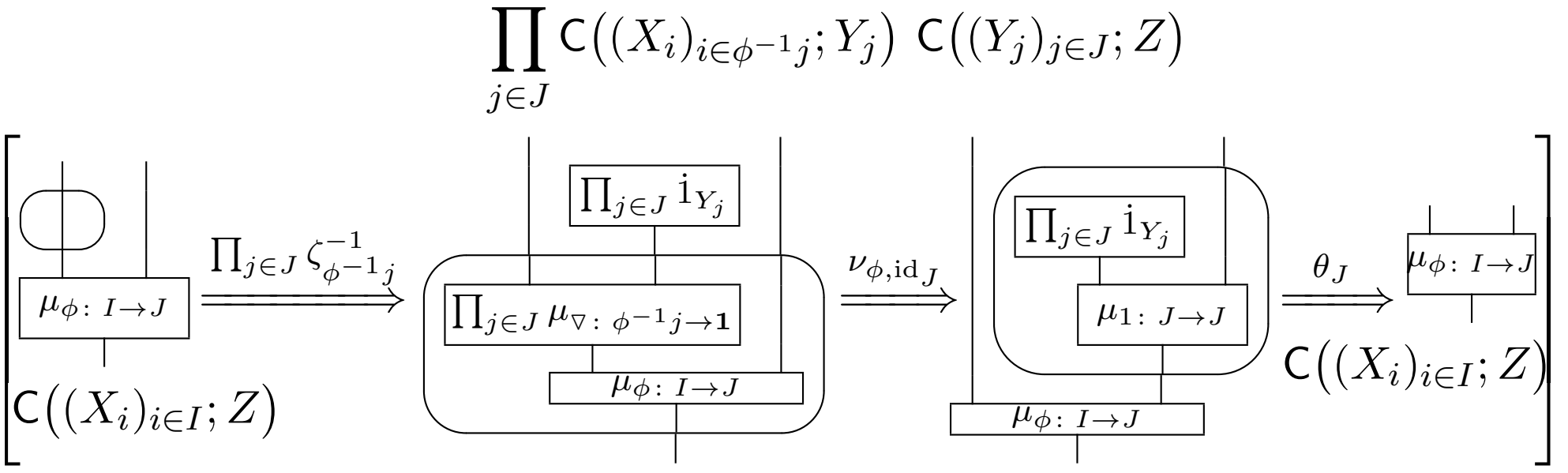

(1.1.4)

Equivalently,

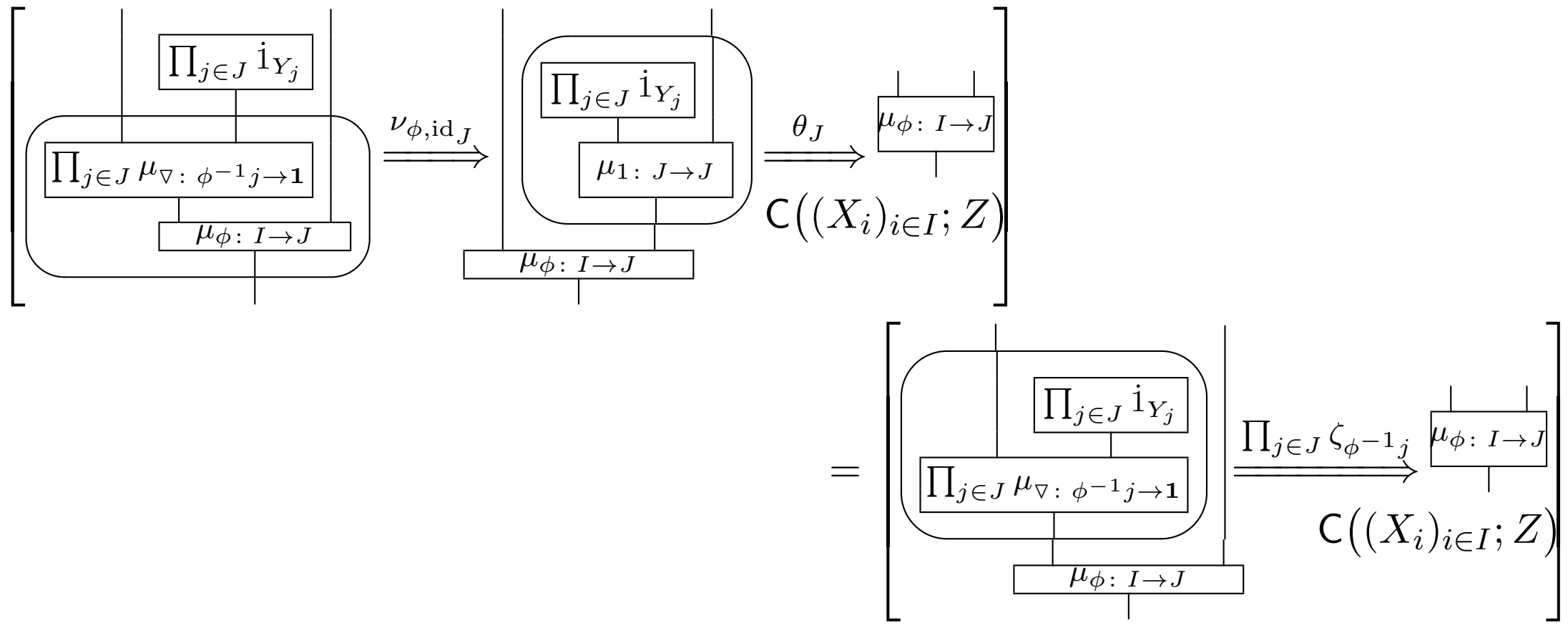

. (1.1.5)

**1.1.2 Example.** The following data give a symmetric weak multicategory $\mathsf{C}$:

— $\operatorname{Ob}\mathsf{C}$ = class of small categories;

— $\mathsf{C}\bigl((\mathcal{C}_i)_{i\in I}; \mathcal{D}\bigr) = \underline{\mathcal{C}at}(\prod_{i\in I}\mathcal{C}_i, \mathcal{D})$ = category of functors $\prod_{i\in I}\mathcal{C}_i \to \mathcal{D}$ and of natural transformations of such for each tuple $\bigl((\mathcal{C}_i)_{i\in I}, \mathcal{D}\bigr)$ of small categories;

For a map $\phi \colon I \to J \in \mathcal{S}_{\mathsf{sk}}$ define a bijection $\sigma(\phi) \colon I \to \coprod_{j\in J}\phi^{-1}j \colon i \mapsto (i, \phi i) \in \coprod_J I = I\times J$ as the graph of $\phi$ (see [Lyu23, (2.1.1) and above]).

— the composition functor

$$\mu_\phi \colon \Bigl[\prod_{j\in J}\underline{\mathcal{C}at}\Bigl(\prod_{i\in\phi^{-1}j}\mathcal{C}_i, \mathcal{D}_j\Bigr)\Bigr] \times \underline{\mathcal{C}at}\Bigl(\prod_{j\in J}\mathcal{D}_j, \mathcal{E}\Bigr) \to \underline{\mathcal{C}at}\Bigl(\prod_{i\in I}\mathcal{C}_i, \mathcal{E}\Bigr),$$
$$\bigl((F_j)_{j\in J}, G\bigr) \mapsto c_{\sigma(\phi)} \cdot \Bigl(\prod_{j\in J} F_j\Bigr) \cdot G,$$

for each map $\phi \colon I \to J$ from $\mathcal{S}_{\mathsf{sk}}$ and small categories $\mathcal{C}_i$, $\mathcal{D}_j$, $\mathcal{E}$ with $i \in I$, $j \in J$, where $c_{\sigma(\phi)} \colon \prod_{i\in I}\mathcal{C}_i \to \prod_{j\in J}\prod_{i\in\phi^{-1}j}\mathcal{C}_i$ is a permutation of factors, corresponding to the bijection $\sigma(\phi)$;

— the identity 1-morphism $1_{\mathcal{C}} \in \operatorname{Ob}\underline{\mathcal{C}at}(\mathcal{C}, \mathcal{C})$ is chosen as the identity functor $\mathrm{Id}_{\mathcal{C}} \colon \mathcal{C} \to \mathcal{C}$;

With this composition diagram (1.1.1) commutes, so we may choose

— $\nu_{\phi,\psi} = \mathrm{id}$;

This makes equation (1.1.2) to hold trivially.

— the right unitor is taken as $\zeta_I = \mathrm{id}$ (note that $\sigma(\triangledown \colon I \to \mathbf{1}) = \mathrm{id}_I$).

— the left unitor is taken as $\theta_I = \mathrm{id}$ (note that $\sigma(\mathrm{id}_I) = \mathrm{id}_I$).

With these choices equation (1.1.3) holds obviously.

We may call *strict* a symmetric weak multicategory for which $\nu_{\phi,\psi} = \mathrm{id}$, $\zeta_I = \mathrm{id}$ and $\theta_I = \mathrm{id}$. Thus $\mathsf{C}$ in this example is strict.

**1.1.3 Proposition.** *For maps $K \xrightarrow{\mathrm{id}} K \xrightarrow{\xi} L \in \mathcal{S}_{\mathsf{sk}}$ and objects $Y_k, W_l, U \in \mathrm{Ob}\,\mathsf{C}$, $k \in K$, $l \in L$, the pasting*

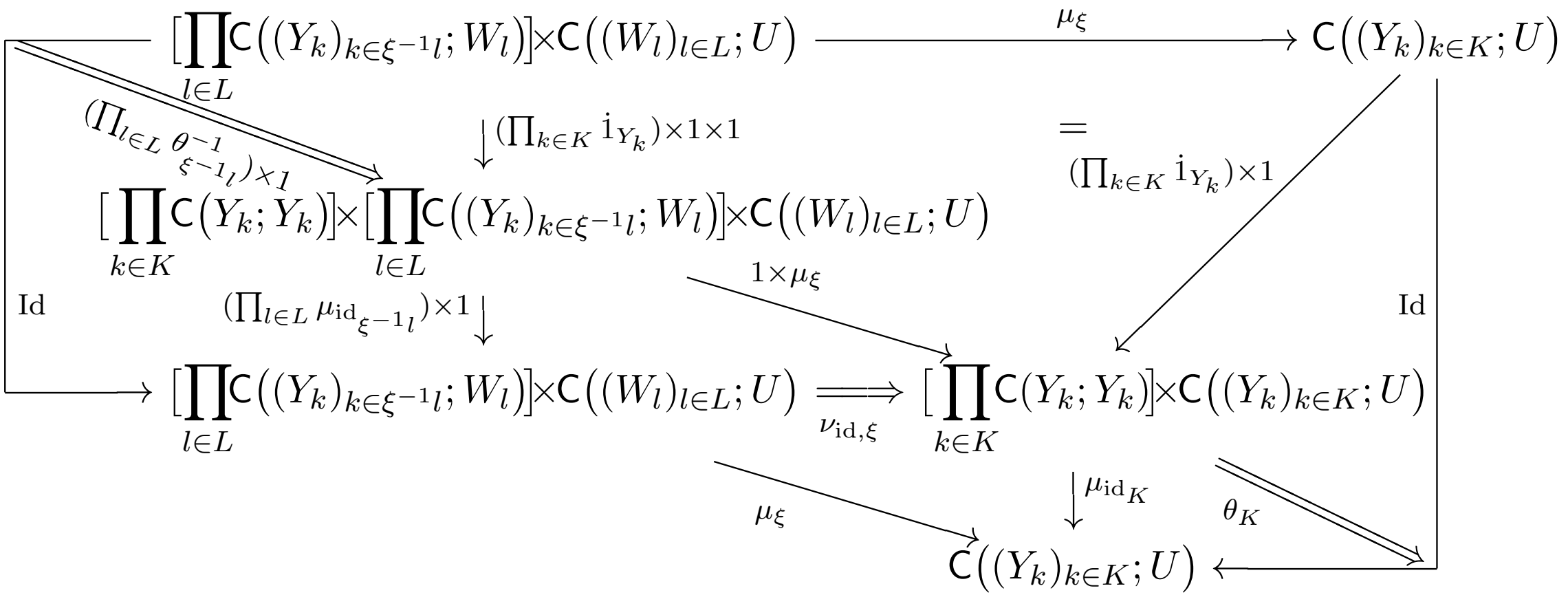

(1.1.6)

*is equal to*

$$\mathrm{id}\colon \mu_\xi \to \mu_\xi\colon [\prod_{l\in L}\mathsf{C}\big((Y_k)_{k\in\xi^{-1}l};W_l\big)] \times \mathsf{C}\big((W_l)_{l\in L};U\big) \to \mathsf{C}\big((Y_k)_{k\in K};U\big). \tag{1.1.7}$$

*The equality* (1.1.6)=(1.1.7) *can also be written in the form*

$$\prod_{l\in L}\mathsf{C}\big((Y_k)_{k\in\xi^{-1}l};W_l\big)\ \mathsf{C}\big((W_l)_{l\in L};U\big)$$

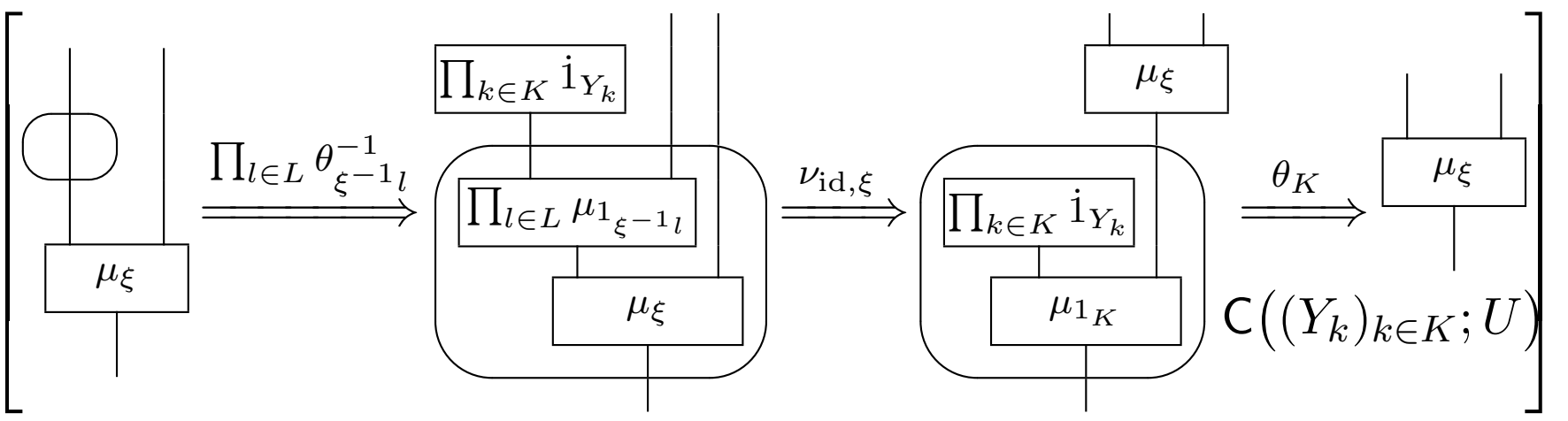

(1.1.8)

*Proof.* We follow the proof of Theorem 7 of Kelly's article [Kel64]. Let $I \xrightarrow{\phi} J \xrightarrow{\mathrm{id}} K \xrightarrow{\xi} L$. Assume that $J = K$, $\psi = \mathrm{id}$, $Y_k = Z_k$ for $k \in K$. By axiom (1.1.2) we have the equation

$$\prod_{j\in J}\mathsf{C}\big((X_i)_{i\in\phi^{-1}j};Y_j\big)\ \prod_{l\in L}\mathsf{C}\big((Z_k)_{k\in\xi^{-1}l};W_l\big)\ \mathsf{C}\big((W_l)_{l\in L};U\big) \tag{1.1.9}$$

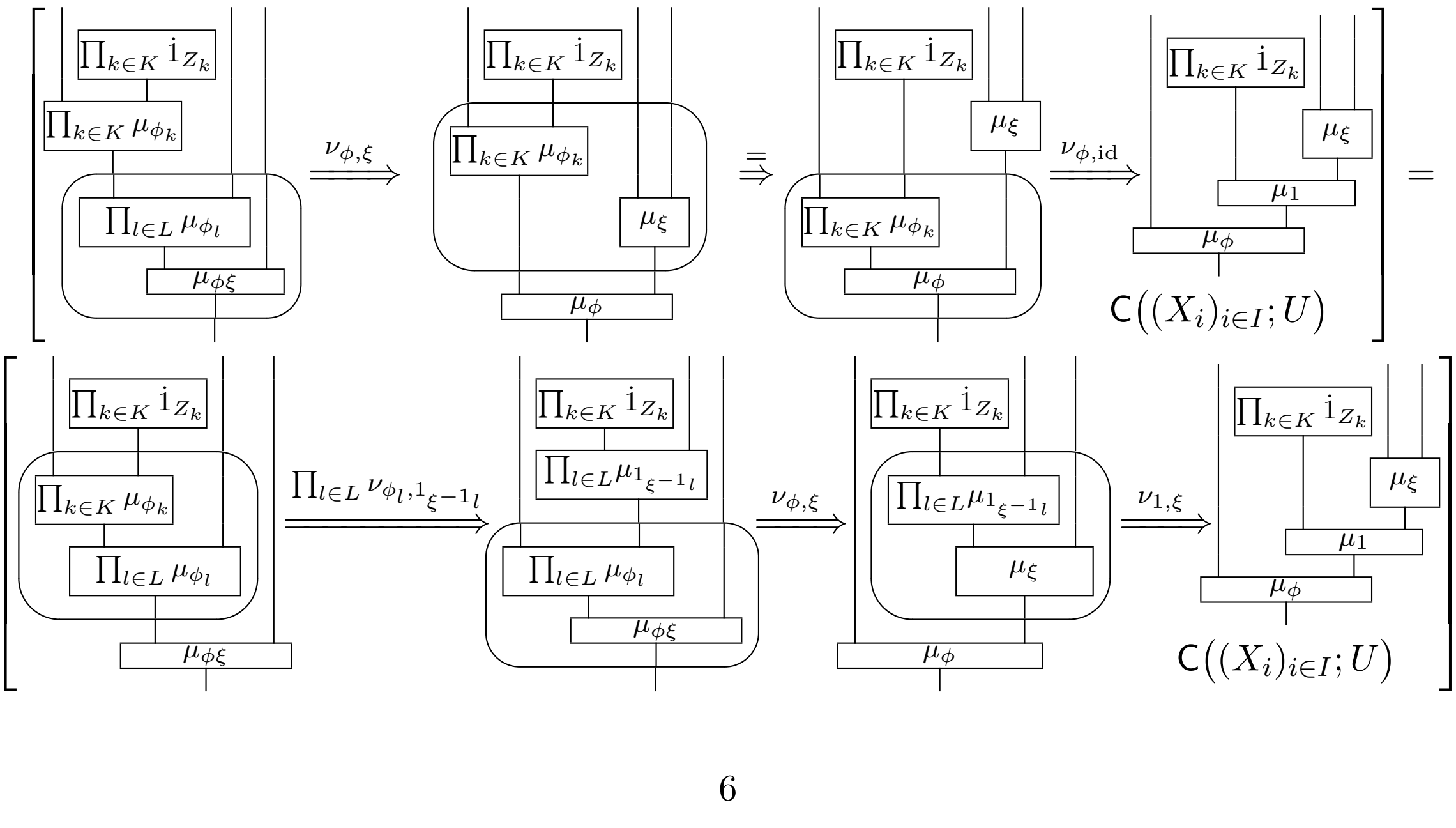

(here $\phi_l = \phi|_{(\phi\xi)^{-1}(l)} \colon (\phi\xi)^{-1}(l) \to \xi^{-1}(l)$, $l \in L$).

Axiom (1.1.3) written for $I \xrightarrow{\phi} J \xrightarrow{\mathrm{id}} J \in \mathcal{S}_{\mathsf{sk}}$ can be presented as the equation

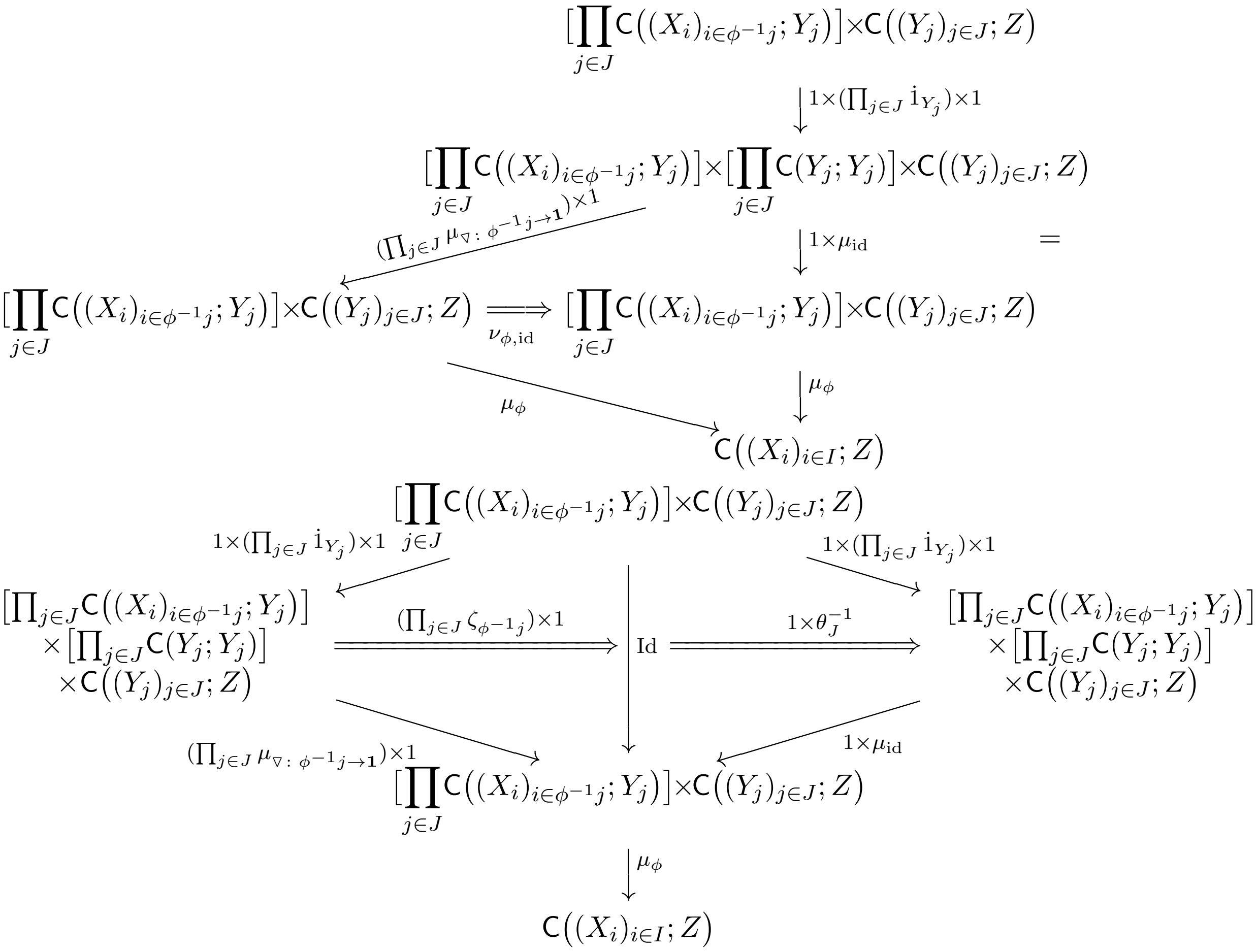

(1.1.10)

Using it, we transform (1.1.9) to the following identity

$$\prod_{j\in J} \mathsf{C}\big((X_i)_{i\in\phi^{-1}j}; Y_j\big) \ \prod_{l\in L} \mathsf{C}\big((Z_k)_{k\in\xi^{-1}l}; W_l\big) \ \mathsf{C}\big((W_l)_{l\in L}; U\big)$$

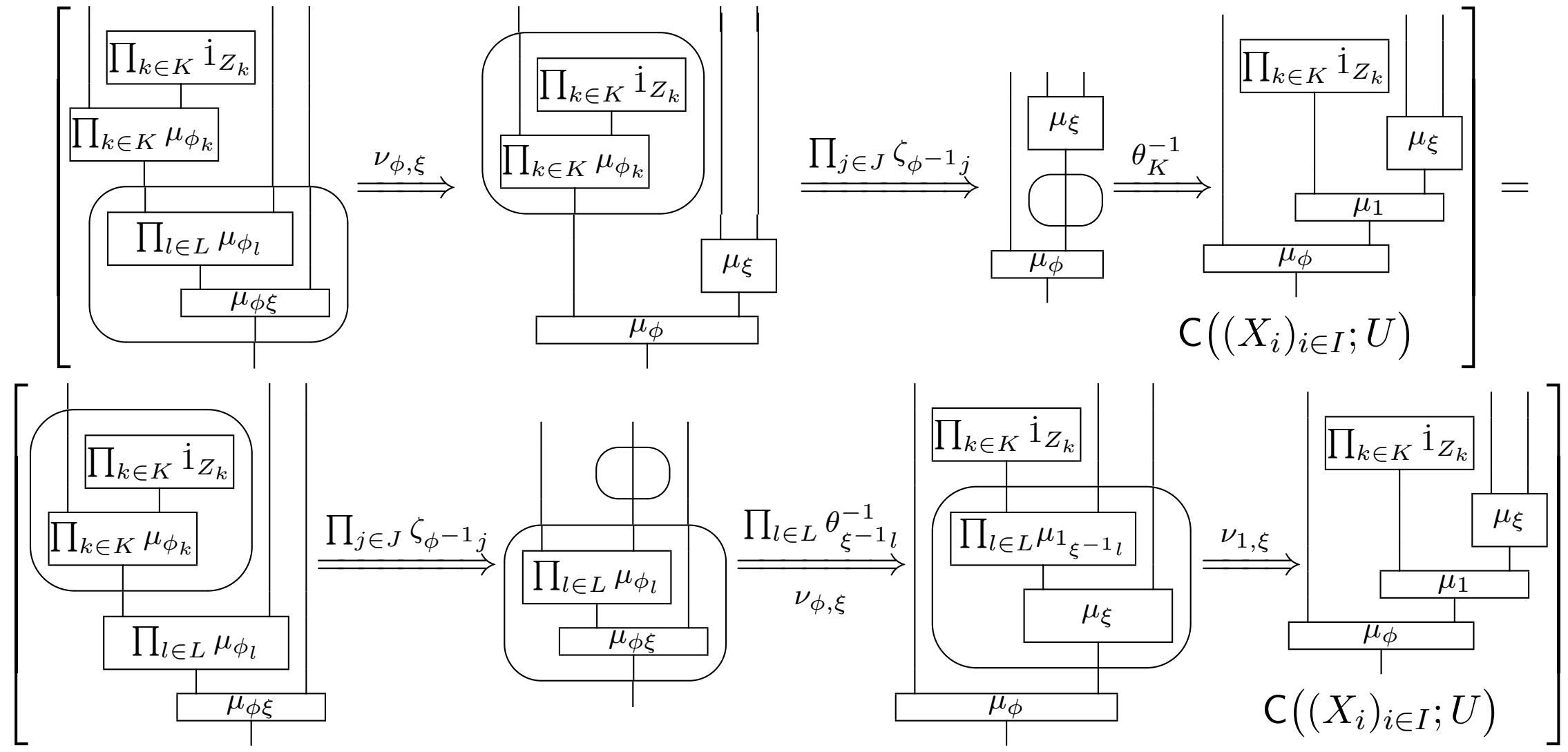

Since $\nu_{\phi,\xi}$ and $\zeta_{\phi^{-1}j}$ are invertible, we may write the equation as the following

$$\prod_{j\in J} \mathsf{C}\big((X_i)_{i\in\phi^{-1}j}; Y_j\big) \ \prod_{l\in L} \mathsf{C}\big((Z_k)_{k\in\xi^{-1}l}; W_l\big) \ \mathsf{C}\big((W_l)_{l\in L}; U\big)$$

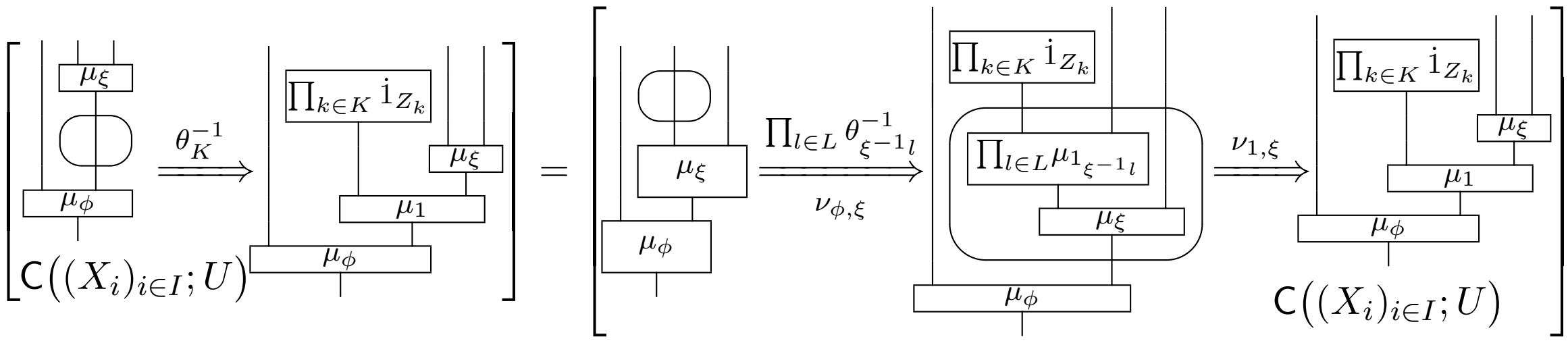

The transformation $\theta_K^{-1}$ can be moved to the right hand side. Thus, the pasting

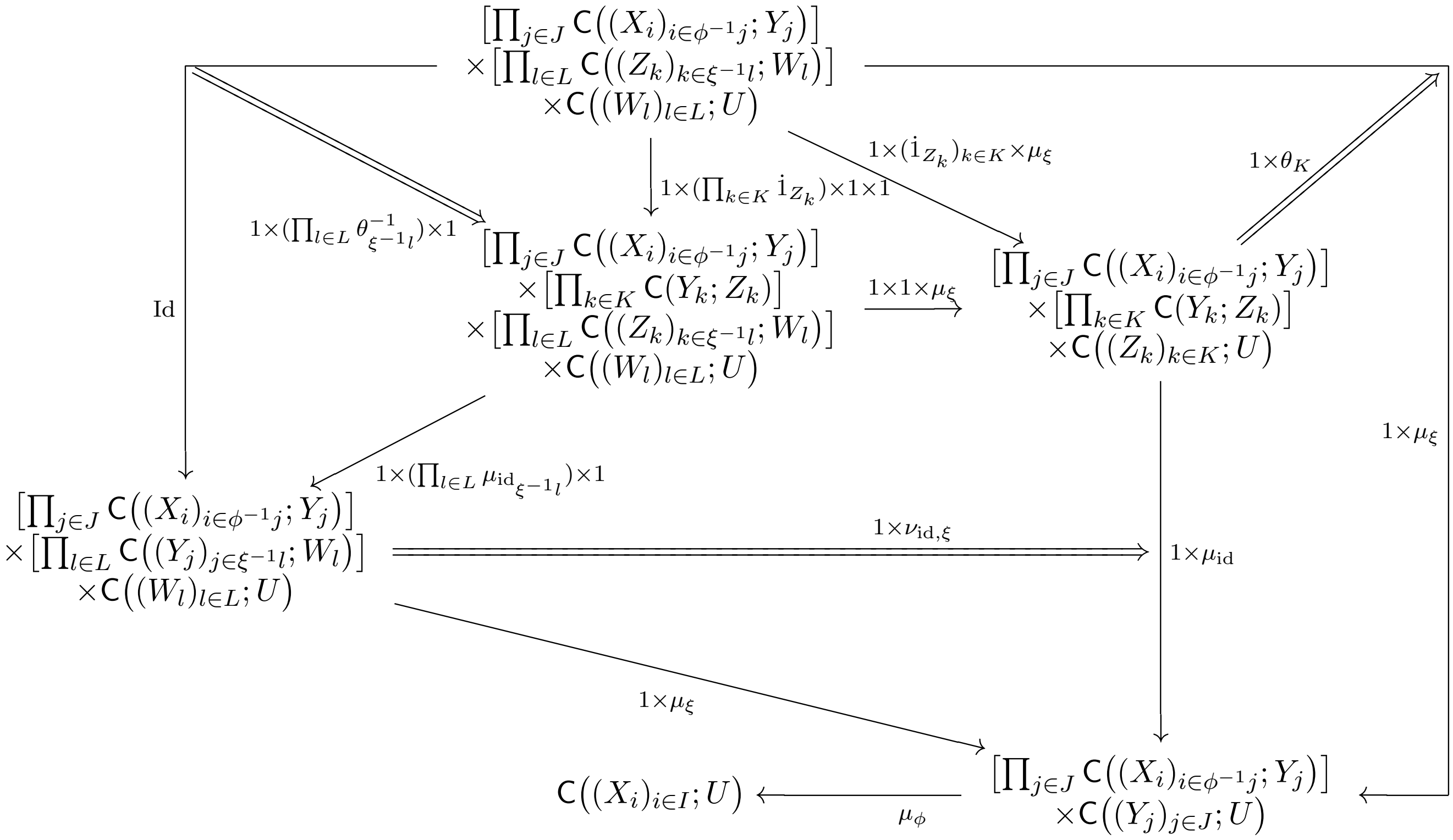

is equal to the identity automorphism of the composite functor

$$\big\{[\prod_{j\in J}\mathsf{C}\big((X_i)_{i\in\phi^{-1}j};Y_j\big)]\times[\prod_{l\in L}\mathsf{C}\big((Y_j)_{j\in\xi^{-1}l};W_l\big)]\times\mathsf{C}\big((W_l)_{l\in L};U\big)$$
$$\xrightarrow{1\times\mu_\xi}[\prod_{j\in J}\mathsf{C}\big((X_i)_{i\in\phi^{-1}j};Y_j\big)]\times\mathsf{C}\big((Y_j)_{j\in J};U\big)\xrightarrow{\mu_\phi}\mathsf{C}\big((X_i)_{i\in I};U\big)\big\}.$$

Let us consider the particular case $I\xrightarrow{\mathrm{id}}J\xrightarrow{\mathrm{id}}K\xrightarrow{\xi}L$. Assume that $I=J=K$, $\phi=\mathrm{id}=\psi$, $X_k=Y_k=Z_k$ for $k\in K$. Compose the above with

$$[\prod_{l\in L}\mathsf{C}\big((Y_j)_{j\in\xi^{-1}l};W_l\big)]\times\mathsf{C}\big((W_l)_{l\in L};U\big)\xrightarrow{(\mathrm{i}_{Y_j})_{j\in J}\times1\times1}$$
$$[\prod_{j\in J}\mathsf{C}\big((X_i)_{i\in\phi^{-1}j};Y_j\big)]\times[\prod_{l\in L}\mathsf{C}\big((Y_j)_{j\in\xi^{-1}l};W_l\big)]\times\mathsf{C}\big((W_l)_{l\in L};U\big).$$

Since $\theta_K$ is invertible the above statement implies that (1.1.6) equals (1.1.7). □

**1.1.4 Proposition.** *For maps* $I\xrightarrow{\phi}J\xrightarrow{\triangledown}\mathbf{1}\in\mathcal{S}_{\mathsf{sk}}$ *and objects* $X_i,Y_j,U\in\operatorname{Ob}\mathsf{C}$, $i\in I$, $j\in J$, *the*

*pasting*

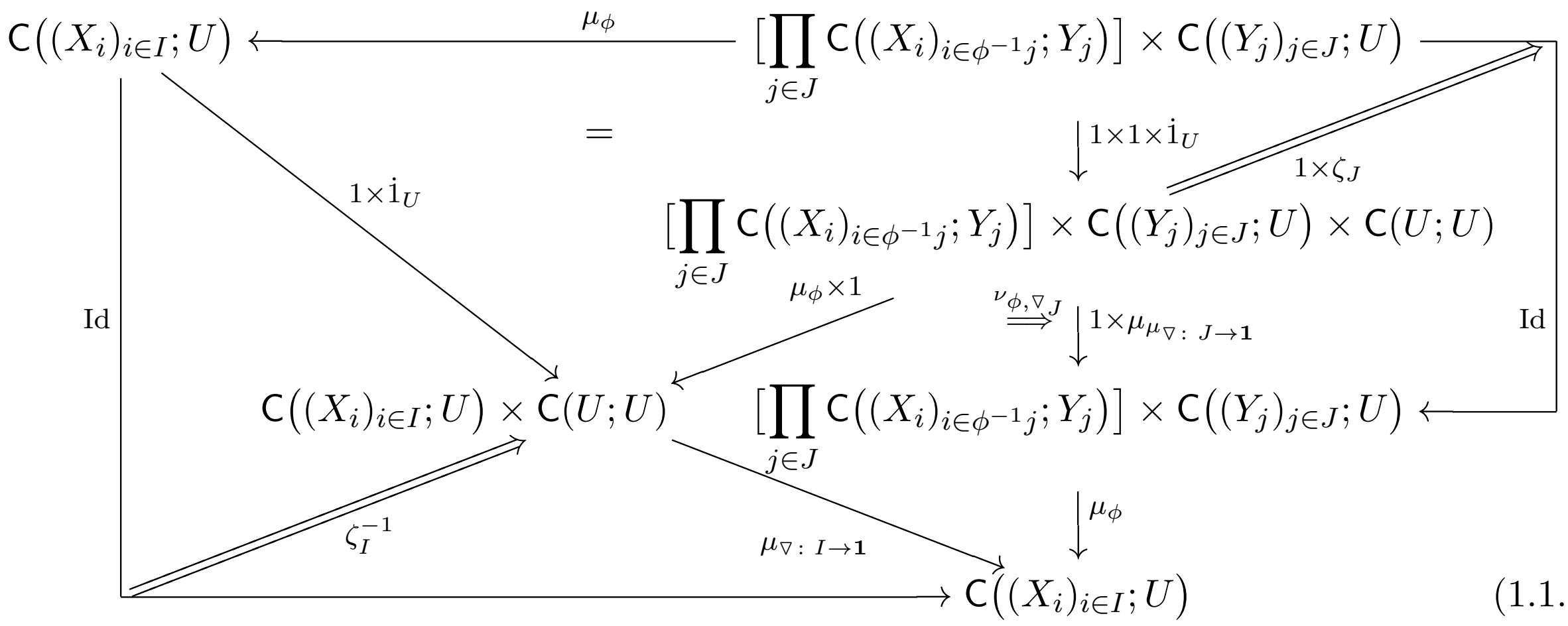

(1.1.11)

*is equal to* $\mathrm{id}\colon \mu_\phi \to \mu_\phi\colon [\prod_{j\in J}\mathsf{C}\big((X_i)_{i\in\phi^{-1}j};Y_j\big)] \times \mathsf{C}\big((Y_j)_{j\in J};U\big) \to \mathsf{C}\big((X_i)_{i\in I};U\big)$. *This equality can also be written in the form*

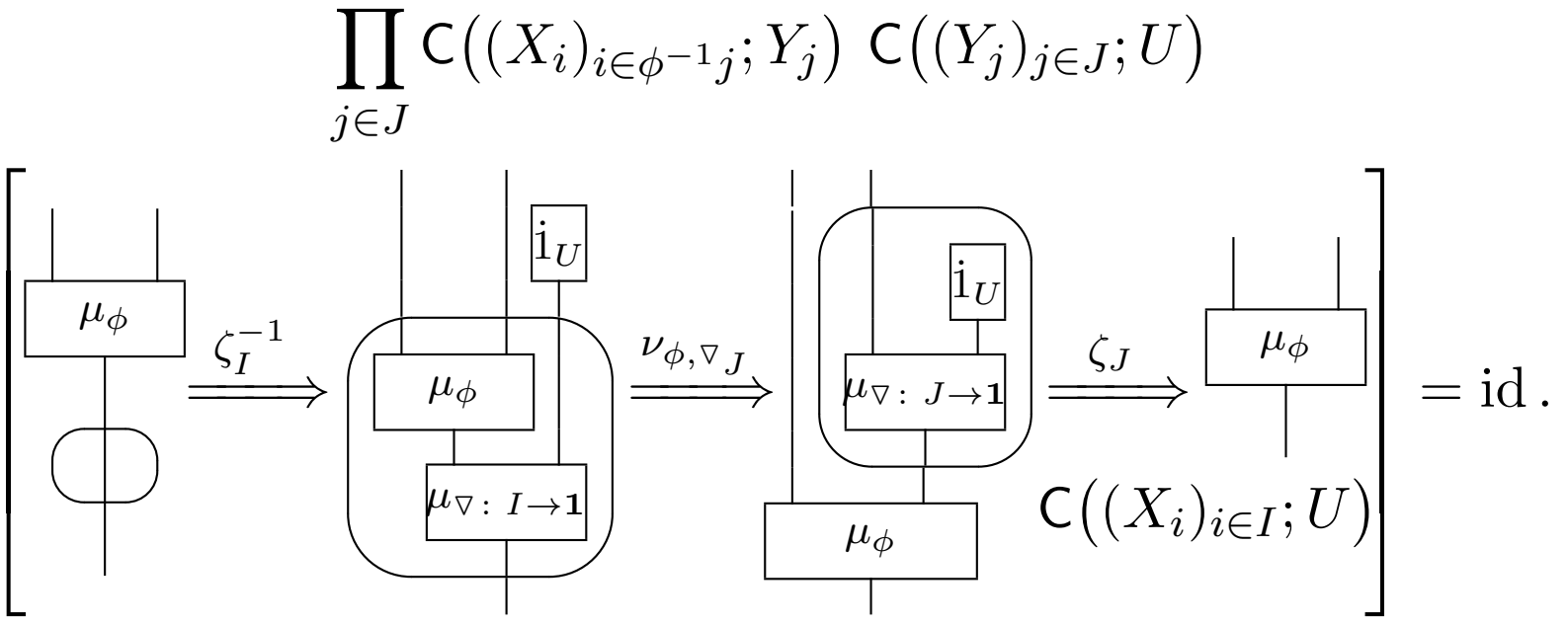

(1.1.12)

*Proof.* Consider the case $I \xrightarrow{\phi} J \xrightarrow{\psi} K \xrightarrow{\mathrm{id}} K \in \mathcal{S}_{\mathsf{sk}}$ of Figure 1. We are given objects $X_i$, $Y_j$, $Z_k$, $U \in \operatorname{Ob}\mathsf{C}$, $i \in I$, $j \in J$, $k \in K$ and assume that $W_k = Z_k$. Compose the diagram with the map $1 \times 1 \times (\mathrm{i}_{Z_k})_{k\in K} \times 1\colon [\prod_{j\in J}\mathsf{C}\big((X_i)_{i\in\phi^{-1}j};Y_j\big)] \times [\prod_{k\in K}\mathsf{C}\big((Y_j)_{j\in\psi^{-1}k};Z_k\big)] \times \mathsf{C}\big((Z_k)_{k\in K};U\big) \to [\prod_{j\in J}\mathsf{C}\big((X_i)_{i\in\phi^{-1}j};Y_j\big)] \times [\prod_{k\in K}\mathsf{C}\big((Y_j)_{j\in\psi^{-1}k};Z_k\big)] \times [\prod_{k\in K}\mathsf{C}(Z_k;Z_k)] \times \mathsf{C}\big((Z_k)_{k\in K};U\big)$. Thus the following equality holds identically

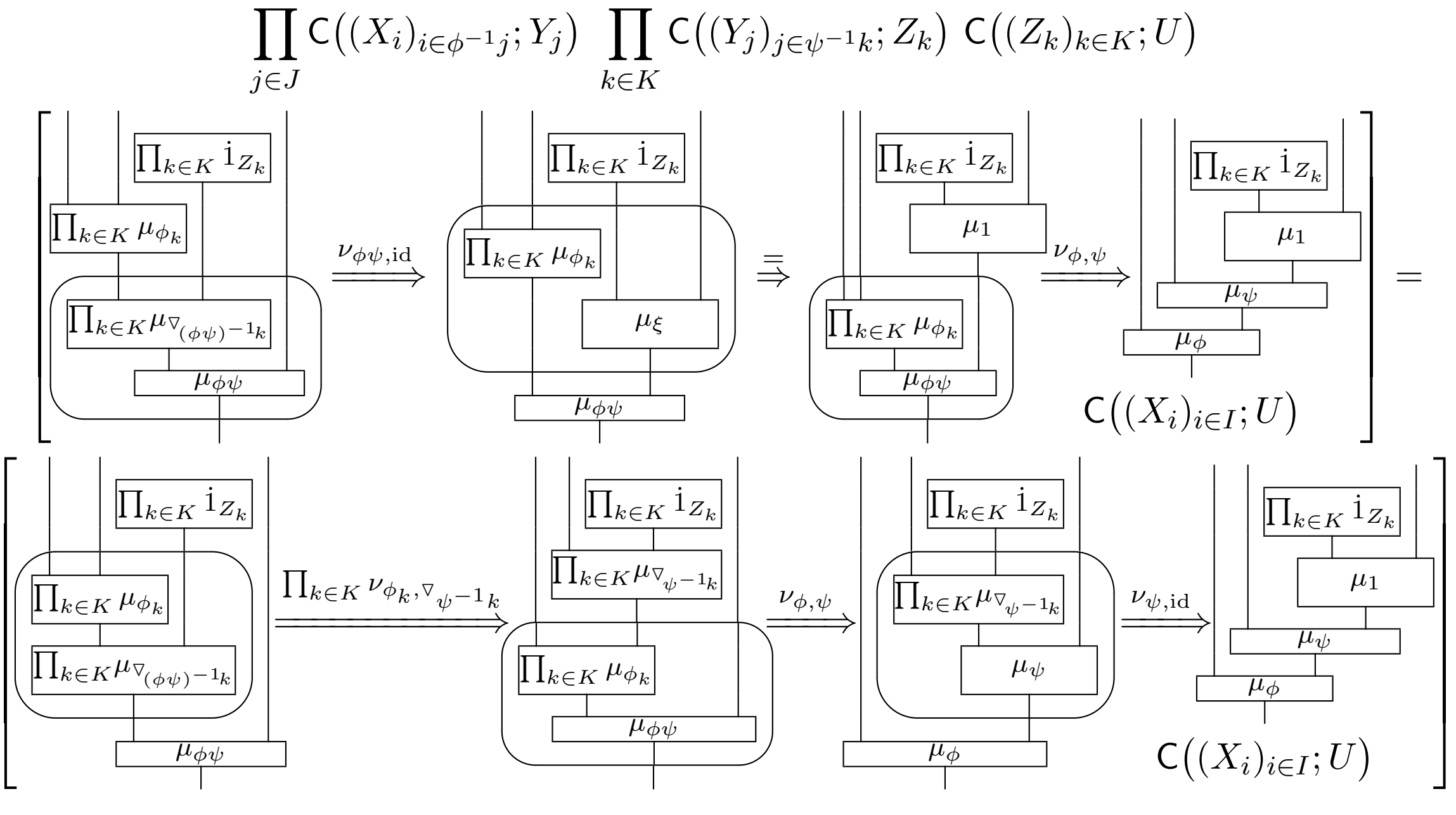

(here $\phi_k = \phi|_{(\phi\psi)^{-1}k}\colon (\phi\psi)^{-1}k \to \psi^{-1}k$, $k \in K$, and $\nabla_S = \nabla\colon S \to \mathbf{1}$).

Let us rewrite this equation using axiom (1.1.10).

$$\prod_{j\in J}\mathsf{C}\big((X_i)_{i\in\phi^{-1}j};Y_j\big)\ \prod_{k\in K}\mathsf{C}\big((Y_j)_{j\in\psi^{-1}k};Z_k\big)\ \mathsf{C}\big((Z_k)_{k\in K};U\big)$$

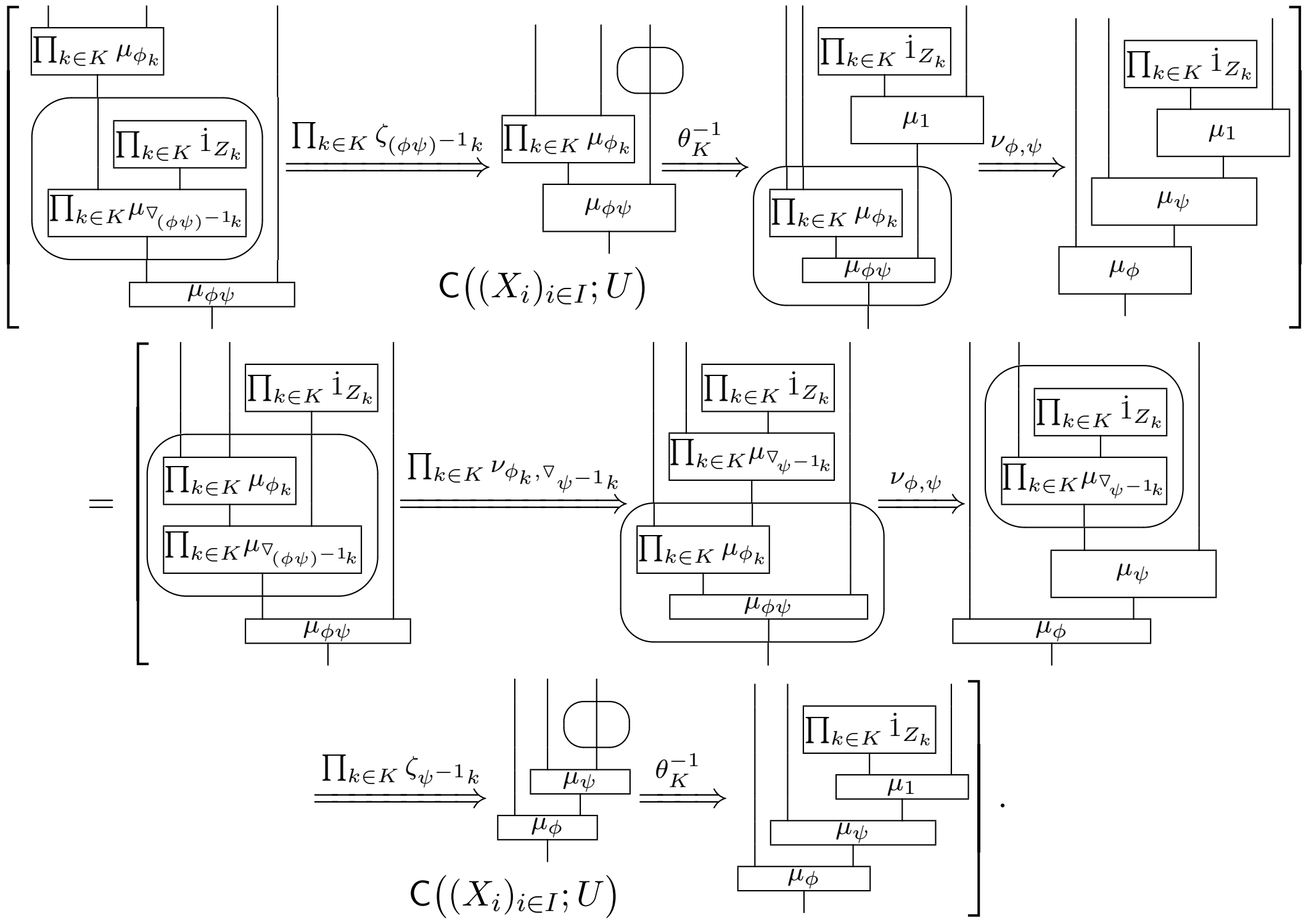

Since $\nu_{\phi,\psi}$ and $\theta_k$ are invertible, we reduce the equation to

$$\prod_{j\in J}\mathsf{C}\big((X_i)_{i\in\phi^{-1}j};Y_j\big)\ \prod_{k\in K}\mathsf{C}\big((Y_j)_{j\in\psi^{-1}k};Z_k\big)\ \mathsf{C}\big((Z_k)_{k\in K};U\big)$$

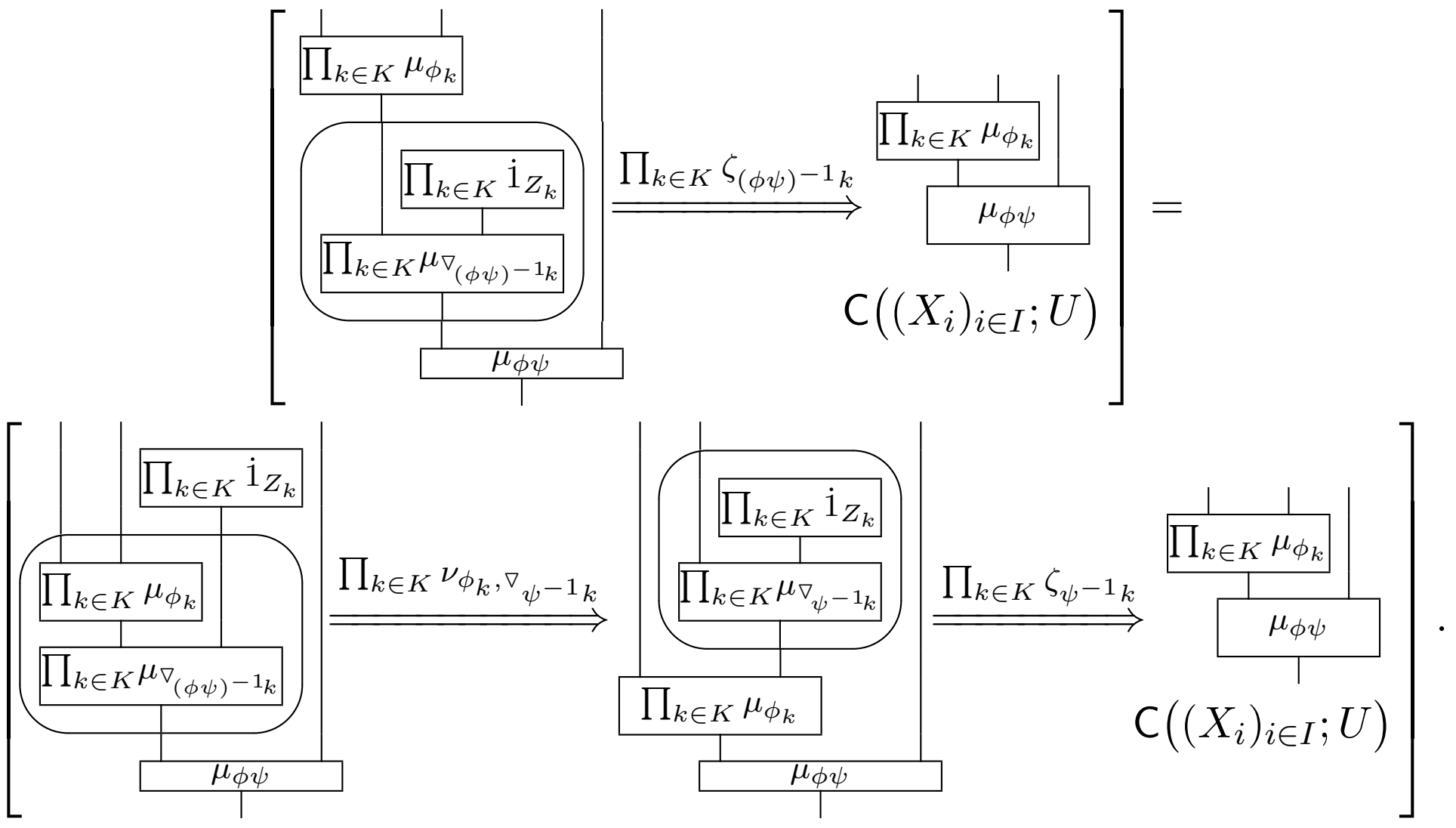

So the composition

$$\prod_{j\in J}\mathsf{C}\big((X_i)_{i\in\phi^{-1}j};Y_j\big)\ \prod_{k\in K}\mathsf{C}\big((Y_j)_{j\in\psi^{-1}k};Z_k\big)\ \mathsf{C}\big((Z_k)_{k\in K};U\big)$$

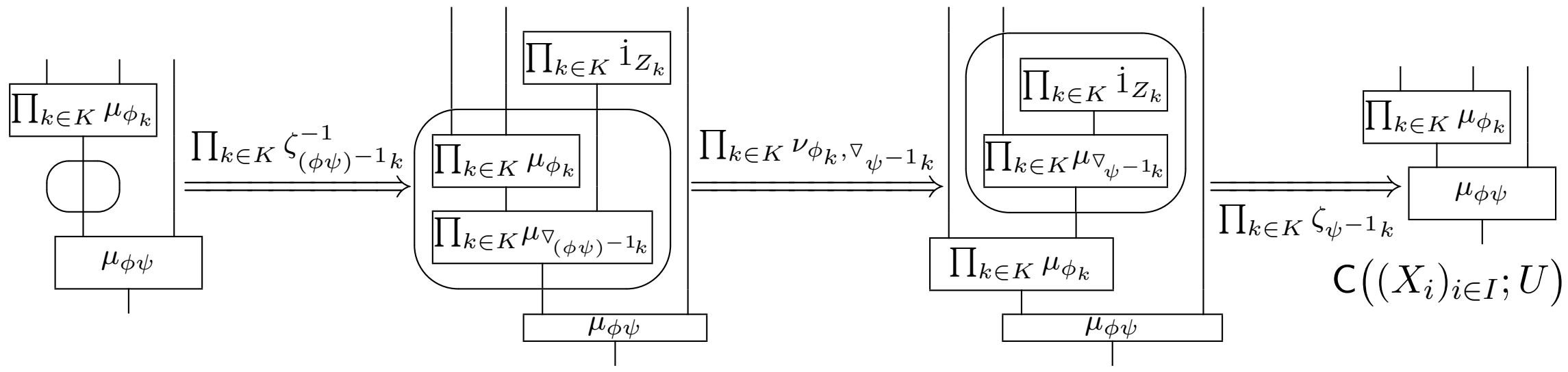

equals the identity transformation.

Suppose now that $K = L = \mathbf{1}$, so we deal with the maps $I \xrightarrow{\phi} J \xrightarrow{\triangledown} \mathbf{1} \xrightarrow{\mathrm{id}} \mathbf{1} \in \mathcal{S}_{\mathsf{sk}}$. Assume that $Z_1 = U$. Precompose the above with the map $1\times 1\times \dot{1}_U : \bigl[\prod_{j\in J} \mathsf{C}\bigl((X_i)_{i\in\phi^{-1}j}; Y_j\bigr)\bigr] \times \mathsf{C}\bigl((Y_j)_{j\in J}; U\bigr) \to \bigl[\prod_{j\in J} \mathsf{C}\bigl((X_i)_{i\in\phi^{-1}j}; Y_j\bigr)\bigr] \times \mathsf{C}\bigl((Y_j)_{j\in J}; U\bigr) \times \mathsf{C}(U;U)$. We deduce that

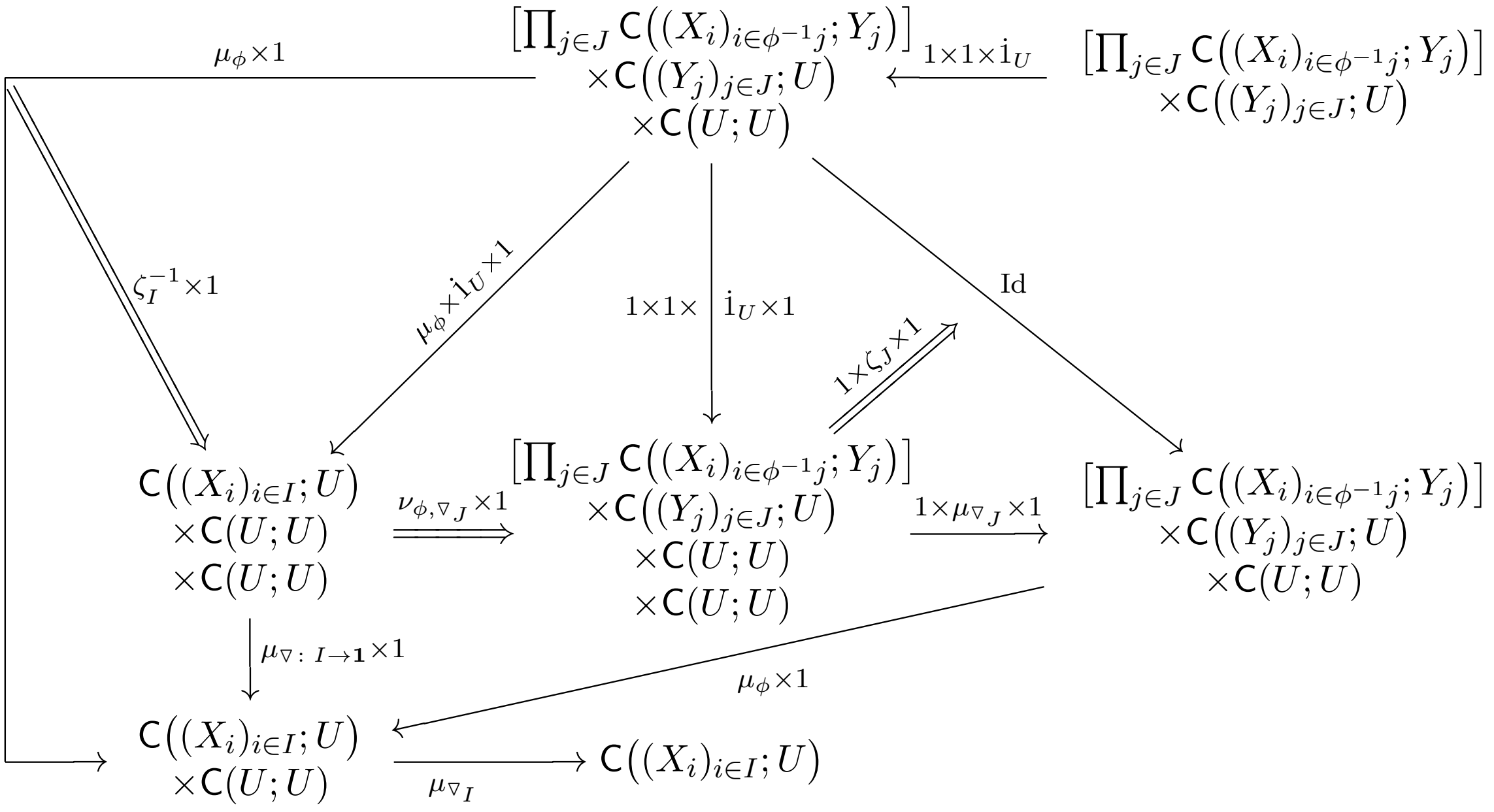

is equal to the identity transformation of the functor

$$\Bigl[\prod_{j\in J} \mathsf{C}\bigl((X_i)_{i\in\phi^{-1}j}; Y_j\bigr)\Bigr] \times \mathsf{C}\bigl((Y_j)_{j\in J}; U\bigr) \xrightarrow{1\times 1\times \dot{1}_U} \Bigl[\prod_{j\in J} \mathsf{C}\bigl((X_i)_{i\in\phi^{-1}j}; Y_j\bigr)\Bigr] \times \mathsf{C}\bigl((Y_j)_{j\in J}; U\bigr) \times \mathsf{C}(U;U) \xrightarrow{\mu_\phi\times 1} \mathsf{C}\bigl((X_i)_{i\in I}; U\bigr) \times \mathsf{C}(U;U) \xrightarrow{\mu_\triangledown} \mathsf{C}\bigl((X_i)_{i\in I}; U\bigr).$$

Since $\zeta_I$ is invertible, we obtain that

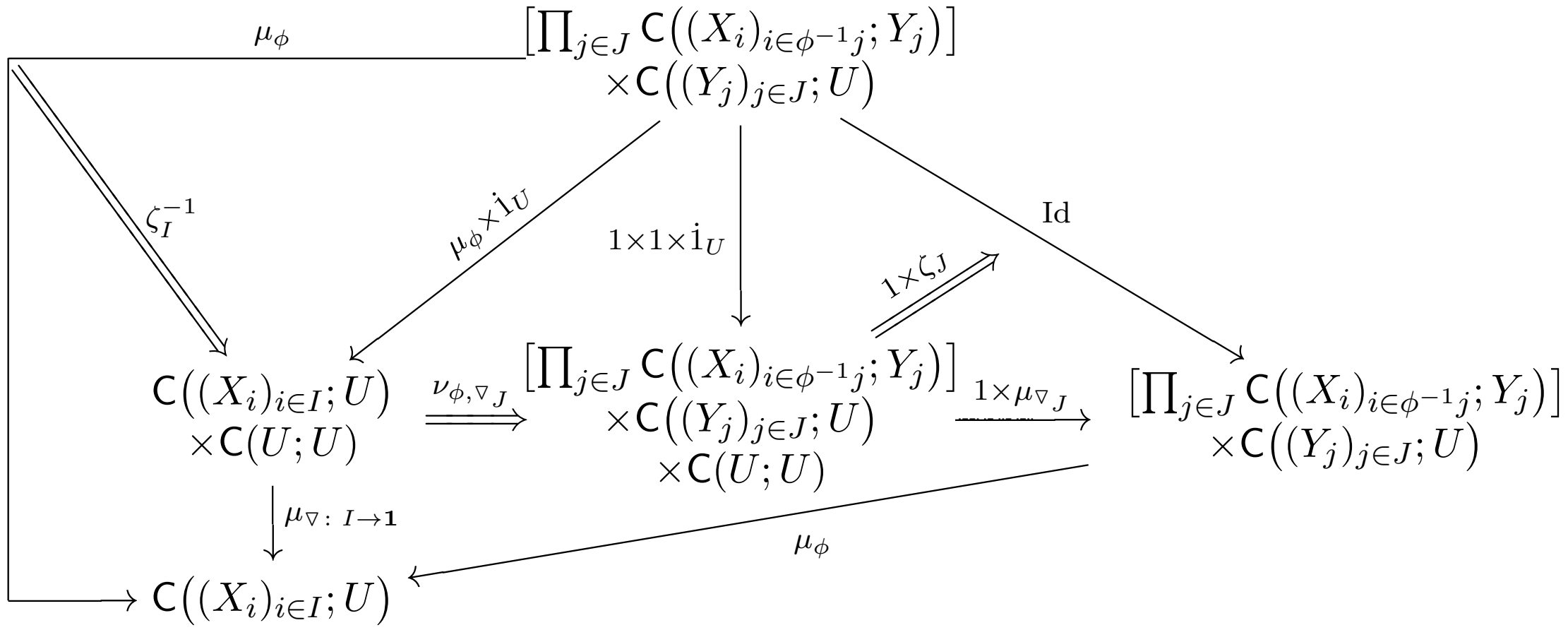

is equal to the identity transformation of the functor

$$\mu_\phi : \Bigl[\prod_{j\in J} \mathsf{C}\bigl((X_i)_{i\in\phi^{-1}j}; Y_j\bigr)\Bigr] \times \mathsf{C}\bigl((Y_j)_{j\in J}; U\bigr) \to \mathsf{C}\bigl((X_i)_{i\in I}; U\bigr). \tag{1.1.13}$$

That is, (1.1.11) equals the identity transformation of (1.1.13). □

**1.1.5 Proposition.** *Each weak multicategory* $\mathsf{C}$ *has an underlying bicategory* $\mathsf{C}_1$ *with*

— $\operatorname{Ob} \mathsf{C}_1 = \operatorname{Ob} \mathsf{C}$;

— $\mathsf{C}_1(X,Y) = \mathsf{C}(X;Y)$;

— *the composition* $\kappa_{X,Y,U} = \mu_{\mathbf{1}\to\mathbf{1}} : \mathsf{C}(X;Y) \times \mathsf{C}(Y;Z) \to \mathsf{C}(X;Z)$;

— *the identity 1-morphism* $1_X \in \operatorname{Ob}\mathsf{C}(X;X)$ *for each object* $X \in \operatorname{Ob}\mathsf{C}$;

— *the associator natural isomorphism*

$$\begin{array}{ccccc} & & \mathsf{C}(X;Y)\times\mathsf{C}(Y;Z)\times\mathsf{C}(Z;W) & & \\ & \swarrow^{\mu_{\mathbf{1}\to\mathbf{1}}\times 1} & & \searrow^{1\times\mu_{\mathbf{1}\to\mathbf{1}}} & \\ \mathsf{C}(X;Z)\times\mathsf{C}(Z;W) & & \xRightarrow{\nu_{\mathrm{id}_\mathbf{1},\mathrm{id}_\mathbf{1}}} & & \mathsf{C}(X;Y)\times\mathsf{C}(Y;W) \\ & \searrow_{\mu_{\mathbf{1}\to\mathbf{1}}} & & \swarrow_{\mu_{\mathbf{1}\to\mathbf{1}}} & \\ & & \mathsf{C}(X;W) & & \end{array}$$

*for all objects* $X, Y, Z, W \in \operatorname{Ob}\mathsf{C}$;

— *the right unitor natural isomorphism*

$$\zeta_\mathbf{1} : \big[\mathsf{C}(X;Z) \xrightarrow{1\times \mathrm{i}_Z} \mathsf{C}(X;Z)\times\mathsf{C}(Z;Z) \xrightarrow{\mu_{\mathbf{1}\to\mathbf{1}}} \mathsf{C}(X;Z)\big] \to \mathrm{Id};$$

— *the left unitor natural isomorphism*

$$\theta_\mathbf{1} : \big[\mathsf{C}(X;Z) \xrightarrow{\mathrm{i}_X\times 1} \mathsf{C}(X;X)\times\mathsf{C}(X;Z) \xrightarrow{\mu_{\mathbf{1}\to\mathbf{1}}} \mathsf{C}(X;Z)\big] \to \mathrm{Id}\,.$$

*Proof.* The associator isomorphism satisfies the associativity coherence equation [Bén67, §1.1(A.C.)] as a particular case of Figure 1. The identity condition [Bén67, §1.1(I.C.)] is satisfied as a particular case of (1.1.3). See Remark 2.1.6. □

**1.1.6 Definition.** A symmetric weak multifunctor $F = (F, F^\phi, F_X) : \mathsf{C} \to \mathsf{D}$ consists of

— a mapping $\operatorname{Ob} F : \operatorname{Ob}\mathsf{C} \to \operatorname{Ob}\mathsf{D}$;

— functors $F = F_{(X_i)_{i\in I};Y} : \mathsf{C}\big((X_i)_{i\in I};Y\big) \to \mathsf{D}\big((FX_i)_{i\in I};FY\big)$;

— a natural isomorphism for each $\phi : I \to J \in \mathcal{S}_{\mathsf{sk}}$

$$\begin{array}{ccc} \big[\prod_{j\in J}\mathsf{C}\big((X_i)_{i\in\phi^{-1}j};Y_j\big)\big]\times\mathsf{C}\big((Y_j)_{j\in J};Z\big) & \xrightarrow{\mu^{\mathsf{C}}_\phi} & \mathsf{C}\big((X_i)_{i\in I};Z\big) \\ {\scriptstyle [\prod_{j\in J}F_{(X_i)_{i\in\phi^{-1}j};Y_j}]\times F_{(Y_j)_{j\in J};Z}}\big\downarrow & \overset{F^\phi_{(X_i),(Y_j),Z}}{\Swarrow} & \big\downarrow{\scriptstyle F_{(X_i)_{i\in I};Z}} \\ \big[\prod_{j\in J}\mathsf{D}\big((FX_i)_{i\in\phi^{-1}j};FY_j\big)\big]\times\mathsf{D}\big((FY_j)_{j\in J};FZ\big) & \xrightarrow{\mu^{\mathsf{D}}_\phi} & \mathsf{D}\big((FX_i)_{i\in I};FZ\big) \end{array} \tag{1.1.14}$$

— an isomorphism for each object $X \in \operatorname{Ob}\mathsf{C}$

$$\begin{array}{ccccc} & \nearrow^{\mathrm{i}^{\mathsf{C}}_X} & \mathsf{C}(X;X) & \searrow^{F_{X;X}} & \\ \mathbb{1} & & \Downarrow F_X & & \mathsf{D}(FX;FX) \\ & & \xrightarrow[\mathrm{i}^{\mathsf{D}}_{FX}]{} & & \end{array} \tag{1.1.15}$$

such that the composition is preserved up to isomorphisms (equation on the next page holds with $\phi_k = \phi| : \phi^{-1}\psi^{-1}k \to \psi^{-1}k$, $k \in K$) and the units are preserved up to isomorphism. For the sake of

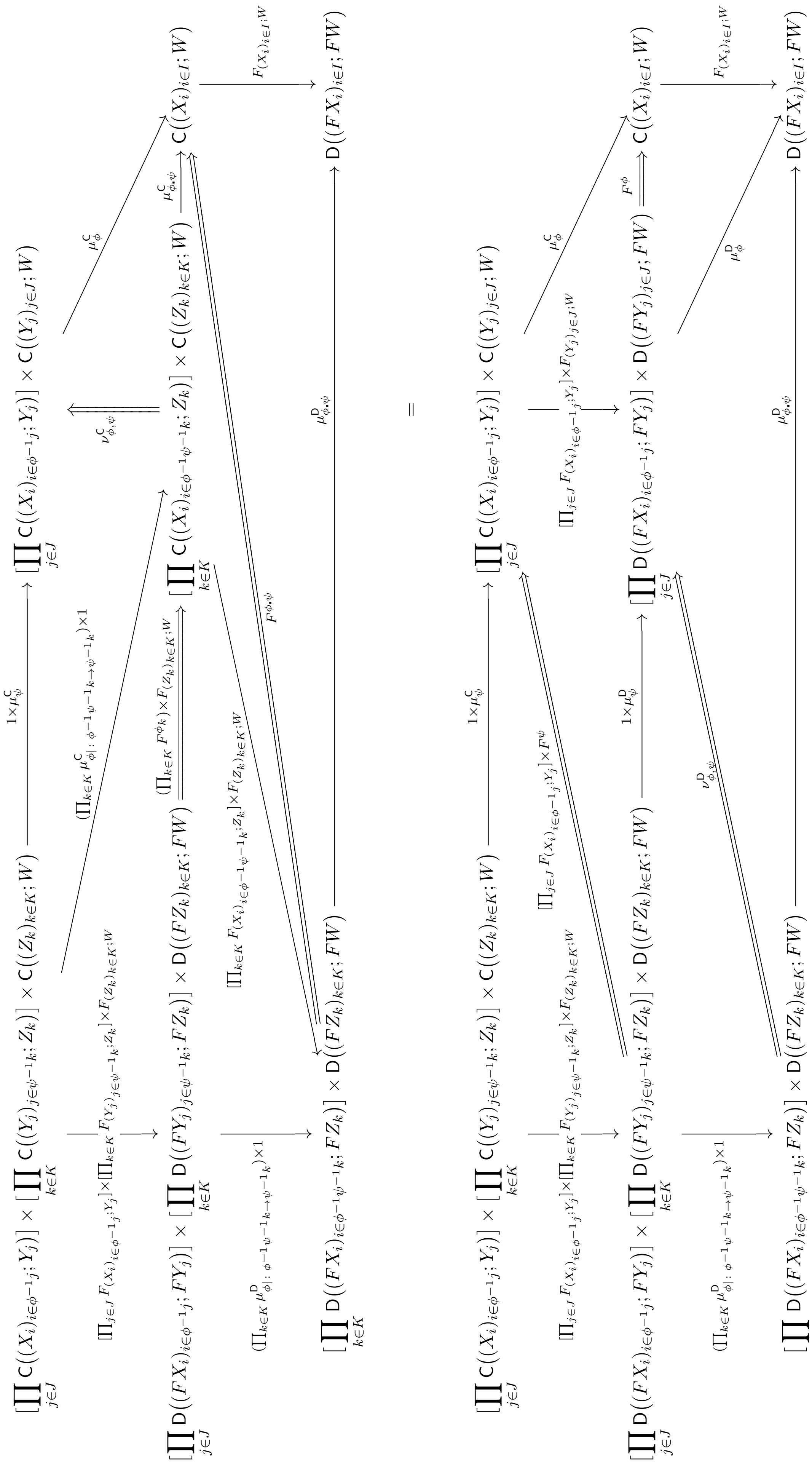

Figure 2: Coherence of weak multifunctor with the composition

uniqueness of such isomorphism we impose two relations

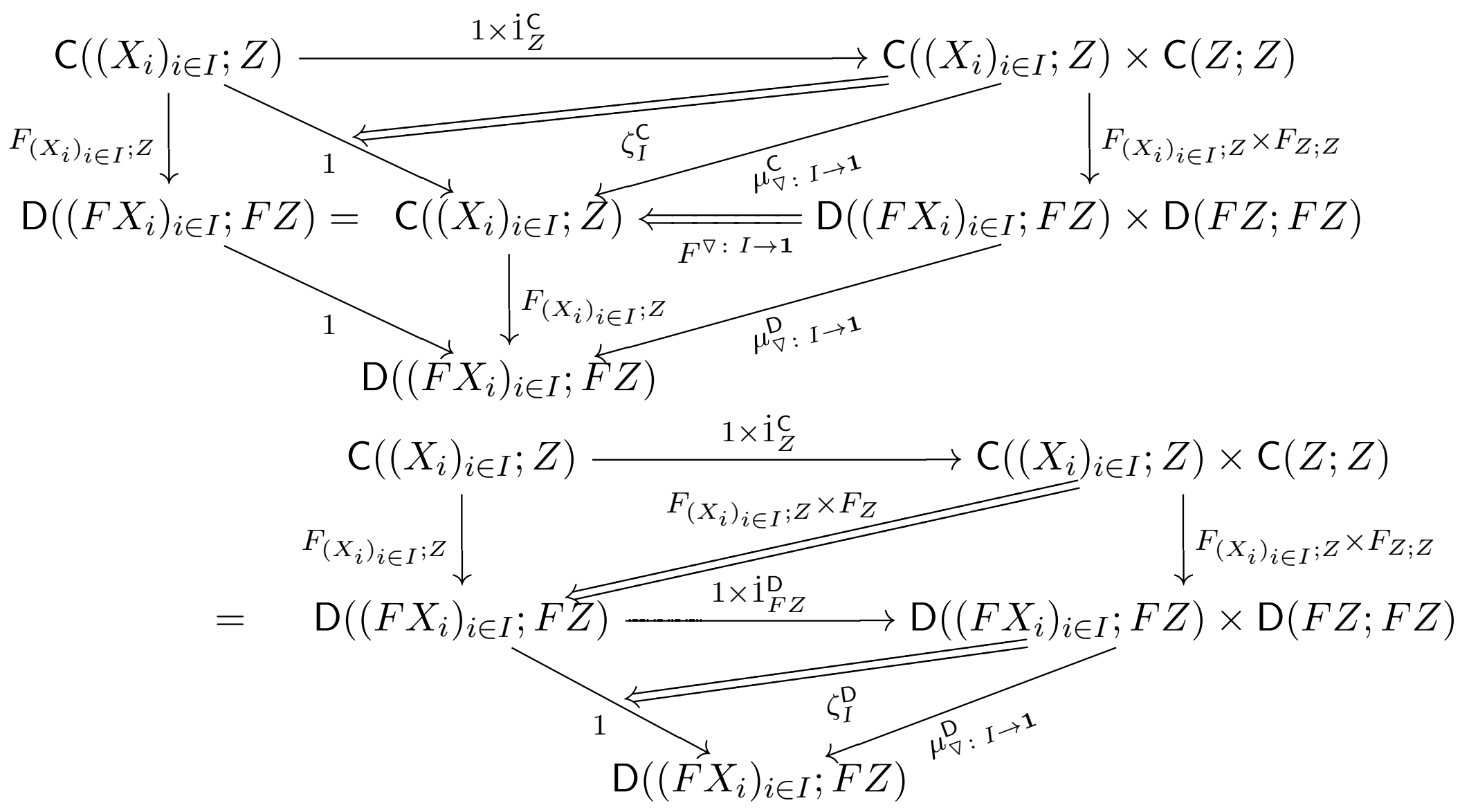

, (1.1.16)

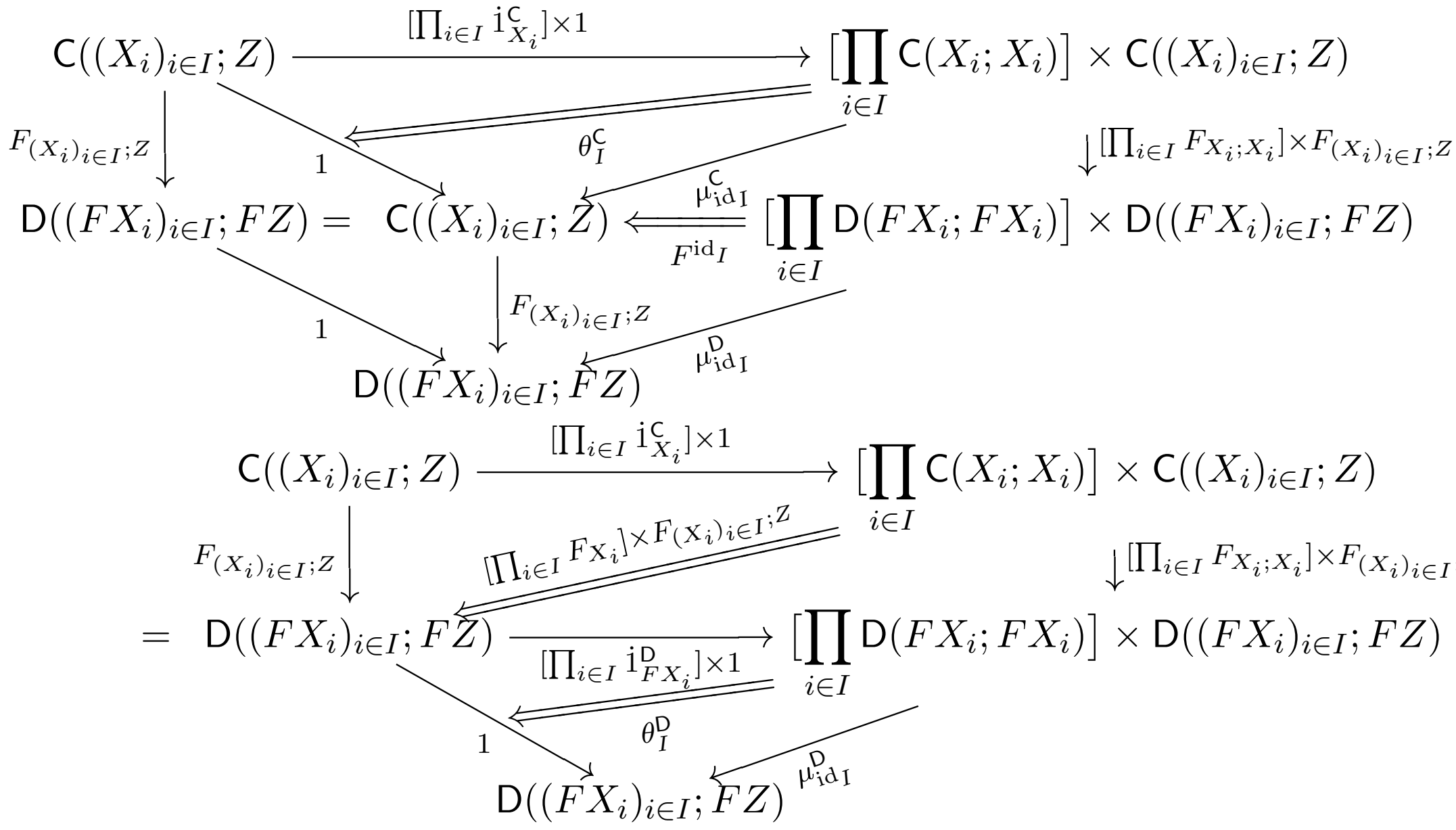

. (1.1.17)

In notations of Bartlett [Bar14] we write Figure 2 as

$$\prod_{j\in J} \mathsf{C}\big((X_i)_{i\in\phi^{-1}j};Y_j\big) \prod_{k\in K} \mathsf{C}\big((Y_j)_{j\in\psi^{-1}k};Z_k\big)\ \mathsf{C}\big((Z_k)_{k\in K};W\big) \tag{1.1.18}$$

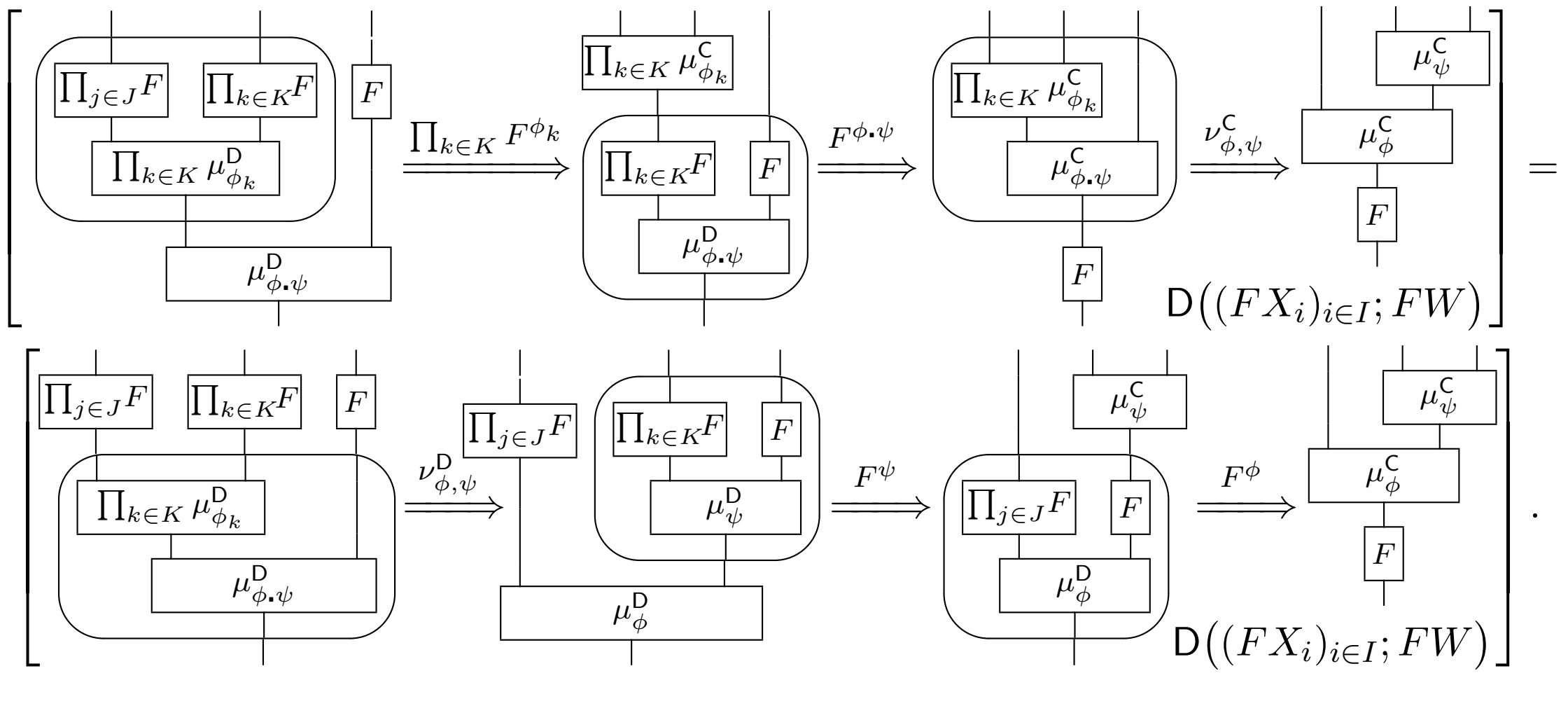

Equation (1.1.16) can be written as

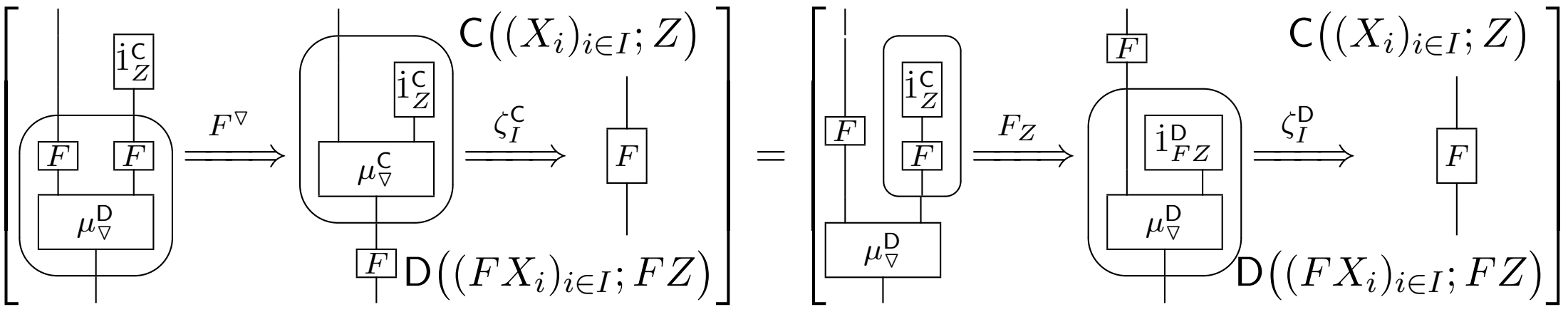

(1.1.19)

Equation (1.1.17) can be written as

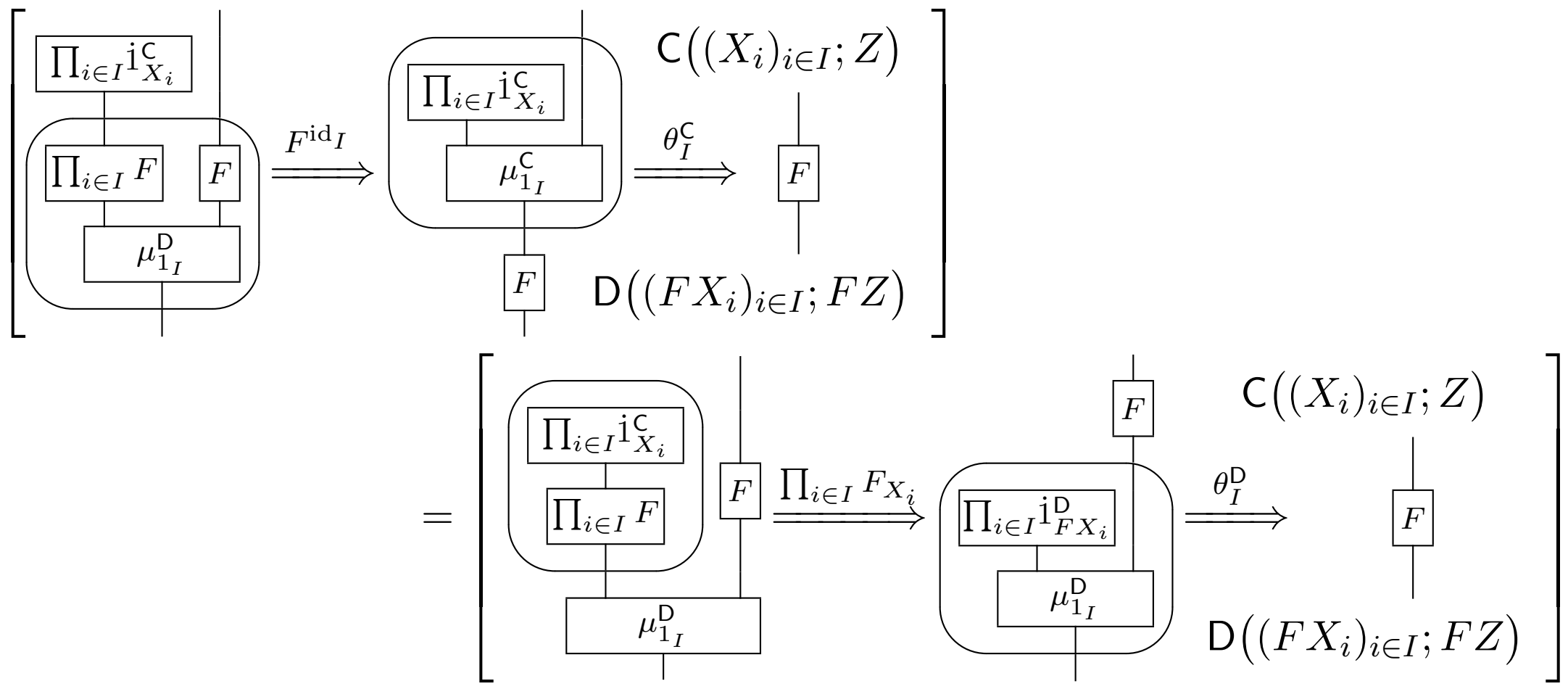

(1.1.20)

Note that weak multifunctors resemble (and are a generalisation of) homomorphisms of bicategories [Bén67, Remark 4.2] or [Lei98, § 1.1]. Having two symmetric weak multifunctors $F\colon \mathsf{C}\to\mathsf{D}$, $G\colon \mathsf{D}\to\mathsf{E}$, we define their composition $K$ with the components

— $\operatorname{Ob} K = \big(\operatorname{Ob}\mathsf{C} \xrightarrow{\operatorname{Ob}F} \operatorname{Ob}\mathsf{D} \xrightarrow{\operatorname{Ob}G} \operatorname{Ob}\mathsf{E}\big)$;

$$K_{(X_i)_{i\in I};Y} = \big[\mathsf{C}((X_i)_{i\in I};Y) \xrightarrow{F_{(X_i)_{i\in I};Y}} \mathsf{D}((FX_i)_{i\in I};FY) \xrightarrow{G_{(FX_i)_{i\in I};FY}} \mathsf{E}((GFX_i)_{i\in I};GFY)\big]; \quad (1.1.21)$$

— $$K^\phi = \begin{array}{ccc} \big[\prod_{j\in J}\mathsf{C}\big((X_i)_{i\in\phi^{-1}j};Y_j\big)\big]\times\mathsf{C}\big((Y_j)_{j\in J};Z\big) & \xrightarrow{\mu_\phi^{\mathsf{C}}} & \mathsf{C}\big((X_i)_{i\in I};Z\big) \\ \big\downarrow{\scriptstyle [\prod_{j\in J}F_{(X_i)_{i\in\phi^{-1}j};Y_j}]\times F_{(Y_j)_{j\in J};Z}} & \overset{F^\phi}{\Longrightarrow} & \big\downarrow{\scriptstyle F_{(X_i)_{i\in I};Z}} \\ \big[\prod_{j\in J}\mathsf{D}\big((FX_i)_{i\in\phi^{-1}j};FY_j\big)\big]\times\mathsf{D}\big((FY_j)_{j\in J};FZ\big) & \xrightarrow{\mu_\phi^{\mathsf{D}}} & \mathsf{D}\big((FX_i)_{i\in I};FZ\big) \\ \big\downarrow{\scriptstyle [\prod_{j\in J}G_{(FX_i)_{i\in\phi^{-1}j};FY_j}]\times G_{(FY_j)_{j\in J};FZ}} & \overset{G^\phi}{\Longrightarrow} & \big\downarrow{\scriptstyle G_{(FX_i)_{i\in I};FZ}} \\ \big[\prod_{j\in J}\mathsf{E}\big((GFX_i)_{i\in\phi^{-1}j};GFY_j\big)\big]\times\mathsf{E}\big((GFY_j)_{j\in J};GFZ\big) & \xrightarrow{\mu_\phi^{\mathsf{E}}} & \mathsf{E}\big((GFX_i)_{i\in I};GFZ\big) \end{array};$$

— $$K_X = \begin{array}{ccc} \mathbb{1} & \xrightarrow{\mathrm{i}_X^{\mathsf{C}}} & \mathsf{C}(X;X) \\ \Big\| & \overset{F_X}{\Longleftarrow} & \big\downarrow{\scriptstyle F_{X;X}} \\ \mathbb{1} & \xrightarrow{\mathrm{i}_{FX}^{\mathsf{D}}} & \mathsf{D}(FX;FX) \\ \Big\| & \overset{G_{FX}}{\Longleftarrow} & \big\downarrow{\scriptstyle G_{FX;FX}} \\ \mathbb{1} & \xrightarrow{\mathrm{i}_{GFX}^{\mathsf{E}}} & \mathsf{E}(GFX;GFX) \end{array}.$$

Stacking one commutative cube on top of another, we prove that equation at Figure 2 holds for $\chi$. Similarly, stacking one commutative prism (1.1.16) on top of another, we get commutative prism (1.1.16) for $\kappa$. Likewise, stacking one commutative prism (1.1.17) on top of another, we get commutative prism (1.1.17) for $\kappa$.

**1.1.7 Remark.** The identity multifunctor $(\mathrm{Id},\mathrm{id},\mathrm{id})\colon \mathsf{C} \to \mathsf{C}$ consists of identity function $\mathrm{id} : \mathrm{Ob}\,\mathsf{C} \to \mathrm{Ob}\,\mathsf{C}$, identity functors $\mathrm{Id}\colon \mathsf{C}((X_i)_{i\in I};Y) \to \mathsf{C}((X_i)_{i\in I};Y)$ and identity transformations $\mathrm{Id}^\phi = \mathrm{id}$ and $\mathrm{Id}_X = \mathrm{id}$. Summing up, there is a category of symmetric weak multicategories whose morphisms are symmetric weak multifunctors.

This category is closed under arbitrary small products. In fact, let $(\mathsf{C}_k)_{k\in K}$ be a family of symmetric weak multicategories. Define a symmetric weak multicategory $\mathsf{C}$ with $\mathrm{Ob}\,\mathsf{C} = \prod_{k\in K}\mathrm{Ob}\,\mathsf{C}_k$. An object is typically denoted $X_i = (X_i^k)^{k\in K} \in \mathrm{Ob}\,\mathsf{C}$. Define $\mathsf{C}\bigl((X_i)_{i\in I};Y\bigr) = \mathsf{C}\bigl((X_i^k)_{i\in I}^{k\in K};(Y^k)^{k\in K}\bigr) \overset{\text{def}}{=} \prod_{k\in K}\mathsf{C}_k\bigl((X_i^k)_{i\in I};Y^k\bigr)$. The composition and the units are componentwise. The projection together with identity transformations $(\mathrm{pr}_k,\mathrm{id},\mathrm{id})\colon \mathsf{C} \to \mathsf{C}_k$ is a symmetric weak multifunctor. These projections turn $\mathsf{C}$ into the product $\prod_{k\in K}\mathsf{C}_k$ in the category of symmetric weak multicategories.

# 2 Action of symmetric groups on a symmetric weak multicategory

## 2.1 Symmetric group action functors

Let $\beta\colon J \to K \in \mathcal{S}_{\mathsf{sk}}$ be a bijection. Let $(Y_j)_{j\in J}$, $(Z_k)_{k\in K}$, $W$ be (families of) objects of a symmetric weak multicategory $\mathsf{V}$ such that $Z_k = Y_{\beta^{-1}k}$. Similarly to [Lei03, Lemma A.2.2] and to [Lyu23, § A.1] define a functor

$$r_\beta = \bigl\{\mathsf{V}\bigl((Z_k)_{k\in K};W\bigr) \xrightarrow{(\mathrm{i}_{Z_k})_{k\in K}\times 1} \bigl[\prod_{k\in K}\mathsf{V}(Y_{\beta^{-1}k};Z_k)\bigr] \times \mathsf{V}\bigl((Z_k)_{k\in K};W\bigr) \xrightarrow{\mu_\beta} \mathsf{V}\bigl((Y_j)_{j\in J};W\bigr)\bigr\}.$$

**2.1.1 Remark.** Let $\beta = 1\colon K \to K$. Then there is an isomorphism $\theta_K\colon r_{1_K} \to \mathrm{Id}$.

**2.1.2 Proposition.** *Let, furthermore, $\gamma = (I \xrightarrow{\alpha} J \xrightarrow{\beta} K) \in \mathcal{S}_{\mathsf{sk}}$, where $\beta$ is a bijection, and $(X_i)_{i\in I}$ be a family of objects of $\mathsf{V}$. Then there is an isomorphism*

$$\psi_{\alpha,\beta}\colon \mu_\gamma \to \bigl\{\bigl[\prod_{k\in K}\mathsf{V}\bigl((X_i)_{i\in\gamma^{-1}k};Y_{\beta^{-1}k}\bigr)\bigr] \times \mathsf{V}\bigl((Y_{\beta^{-1}k})_{k\in K};W\bigr) \xrightarrow{\prod_{\beta^{-1}}\times r_\beta} \bigl[\prod_{j\in J}\mathsf{V}\bigl((X_i)_{i\in\alpha^{-1}(j)};Y_j\bigr)\bigr] \times \mathsf{V}\bigl((Y_j)_{j\in J};W\bigr) \xrightarrow{\mu_\alpha} \mathsf{V}\bigl((X_i)_{i\in I};W\bigr)\bigr\}. \qquad (2.1.1)$$

*Proof.* Applying the associativity isomorphism from (1.1.1) for maps $I \xrightarrow{\alpha} J \xrightarrow{\beta} K$ we get the sought isomorphism on the following page. In Bartlett's notation

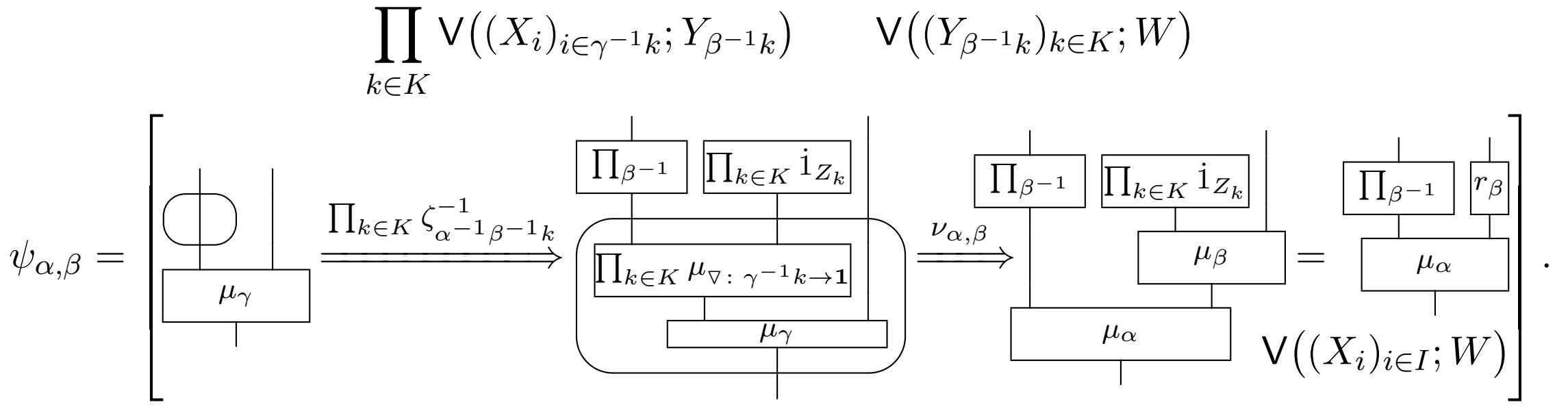

□

**2.1.3 Corollary.** *Assume that both $\alpha$ and $\beta$ are bijections from $\mathcal{S}_{\mathsf{sk}}$, $\gamma = (I \xrightarrow{\alpha} J \xrightarrow{\beta} K)$. Then there is an isomorphism*

$$\phi_{\alpha,\beta}\colon r_\gamma \to \bigl[\mathsf{V}\bigl((Y_{\beta^{-1}k})_{k\in K};W\bigr) \xrightarrow{r_\beta} \mathsf{V}\bigl((Y_j)_{j\in J};W\bigr) \xrightarrow{r_\alpha} \mathsf{V}\bigl((Y_{\alpha i})_{i\in I};W\bigr)\bigr].$$

*Proof.* Consider $X_i = Y_{\alpha i}$, hence, $Y_j = X_{\alpha^{-1}j}$ as well as $Z_k = Y_{\beta^{-1}k}$. Rewrite (2.1.1) as

$$\psi_{\alpha,\beta}\colon \mu_\gamma \to \bigl\{\bigl[\prod_{k\in K}\mathsf{V}\bigl(X_{\gamma^{-1}k};Y_{\beta^{-1}k}\bigr)\bigr] \times \mathsf{V}\bigl((Y_{\beta^{-1}k})_{k\in K};W\bigr) \xrightarrow{\prod_{\beta^{-1}}\times r_\beta} \bigl[\prod_{j\in J}\mathsf{V}\bigl(X_{\alpha^{-1}j};Y_j\bigr)\bigr] \times \mathsf{V}\bigl((Y_j)_{j\in J};W\bigr) \xrightarrow{\mu_\alpha} \mathsf{V}\bigl((X_i)_{i\in I};W\bigr)\bigr\}. \qquad (2.1.2)$$

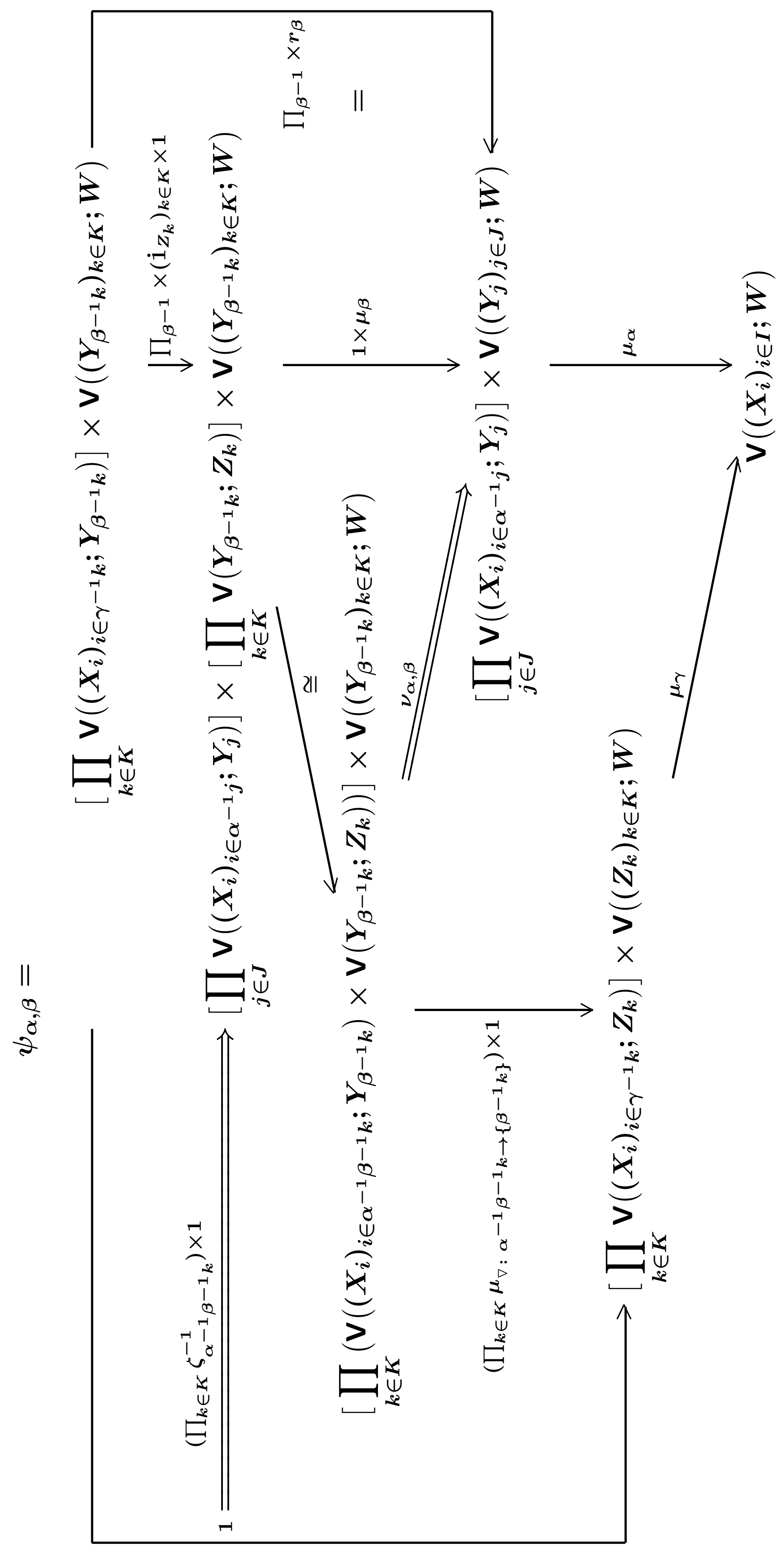

Figure 3: Action of symmetric groups on a symmetric weak multicategory

Substitute $(1_{Y_{\beta^{-1}k}})_{k\in K}$ into the first factor. We get from the source of (2.1.2)

$$\big\{\mathsf{V}\big((X_{\gamma^{-1}k})_{k\in K};W\big) \xrightarrow{(\mathrm{i}_{X_{\gamma^{-1}k}})_{k\in K}\times 1} \big[\prod_{k\in K}\mathsf{V}(X_{\gamma^{-1}k};X_{\gamma^{-1}k})\big]\times\mathsf{V}\big((X_{\gamma^{-1}k})_{k\in K};W\big) \xrightarrow{\mu_\gamma} \mathsf{V}\big((X_i)_{i\in I};W\big)\big\} = r_\gamma.$$

From the target of (2.1.2) we get the functor

$$\begin{aligned}
&\big\{\mathsf{V}\big((Y_{\beta^{-1}k})_{k\in K};W\big) \xrightarrow{(\mathrm{i}_{Y_{\beta^{-1}k}})_{k\in K}\times 1} \big[\prod_{k\in K}\mathsf{V}(X_{\gamma^{-1}k};Y_{\beta^{-1}k})\big]\times\mathsf{V}\big((Y_{\beta^{-1}k})_{k\in K};W\big)\\
&\qquad \xrightarrow{\prod_{\beta^{-1}}\times r_\beta} \big[\prod_{j\in J}\mathsf{V}(X_{\alpha^{-1}j};Y_j)\big]\times\mathsf{V}\big((Y_j)_{j\in J};W\big) \xrightarrow{\mu_\alpha} \mathsf{V}\big((X_i)_{i\in I};W\big)\big\}\\
&= \big\{\mathsf{V}\big((Y_{\beta^{-1}k})_{k\in K};W\big) \xrightarrow{r_\beta} \mathsf{V}\big((Y_j)_{j\in J};W\big)\\
&\qquad \xrightarrow{(\mathrm{i}_{Y_j})_{j\in J}\times 1} \big[\prod_{j\in J}\mathsf{V}(X_{\alpha^{-1}j};Y_j)\big]\times\mathsf{V}\big((Y_j)_{j\in J};W\big) \xrightarrow{\mu_\alpha} \mathsf{V}\big((X_i)_{i\in I};W\big)\big\}\\
&\qquad = \big\{\mathsf{V}\big((Y_{\beta^{-1}k})_{k\in K};W\big) \xrightarrow{r_\beta} \mathsf{V}\big((Y_j)_{j\in J};W\big) \xrightarrow{r_\alpha} \mathsf{V}\big((X_i)_{i\in I};W\big)\big\}.
\end{aligned}$$

Hence, the isomorphism $\phi_{\alpha,\beta}\colon r_{\alpha\centerdot\beta}\to r_\beta\centerdot r_\alpha$. It is presented on the next page. We have used

**2.1.4 Lemma.** *We have* $\zeta_{\mathbf{1}}=\theta_{\mathbf{1}}\colon (1_X,1_X)\mu_{\mathbf{1}\to\mathbf{1}}\to 1_X\in\mathsf{C}(X;X)$.

*Proof.* For any object $X$ of $\mathsf{C}$ the category $\mathsf{C}(X;X)$ is a monoidal category with the unit object $1_X$ (exercise). By a result of Kelly [Kel64] we have the equality

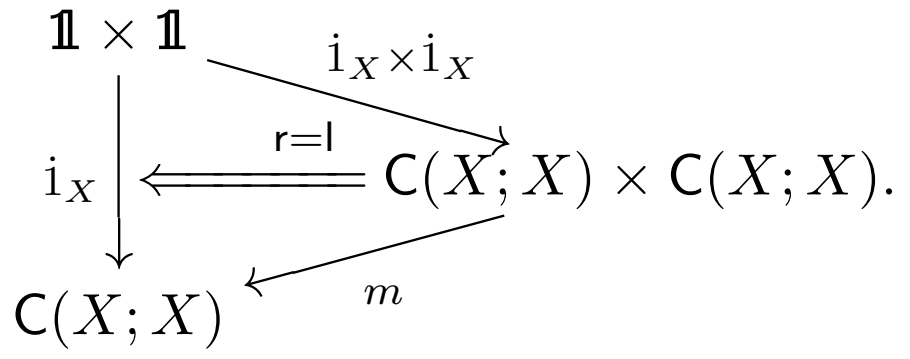

That is, $\zeta_{\mathbf{1}}=\theta_{\mathbf{1}}\colon (1_X,1_X)\mu_{\mathbf{1}\to\mathbf{1}}\to 1_X$. □

□

In Bartlett's notation

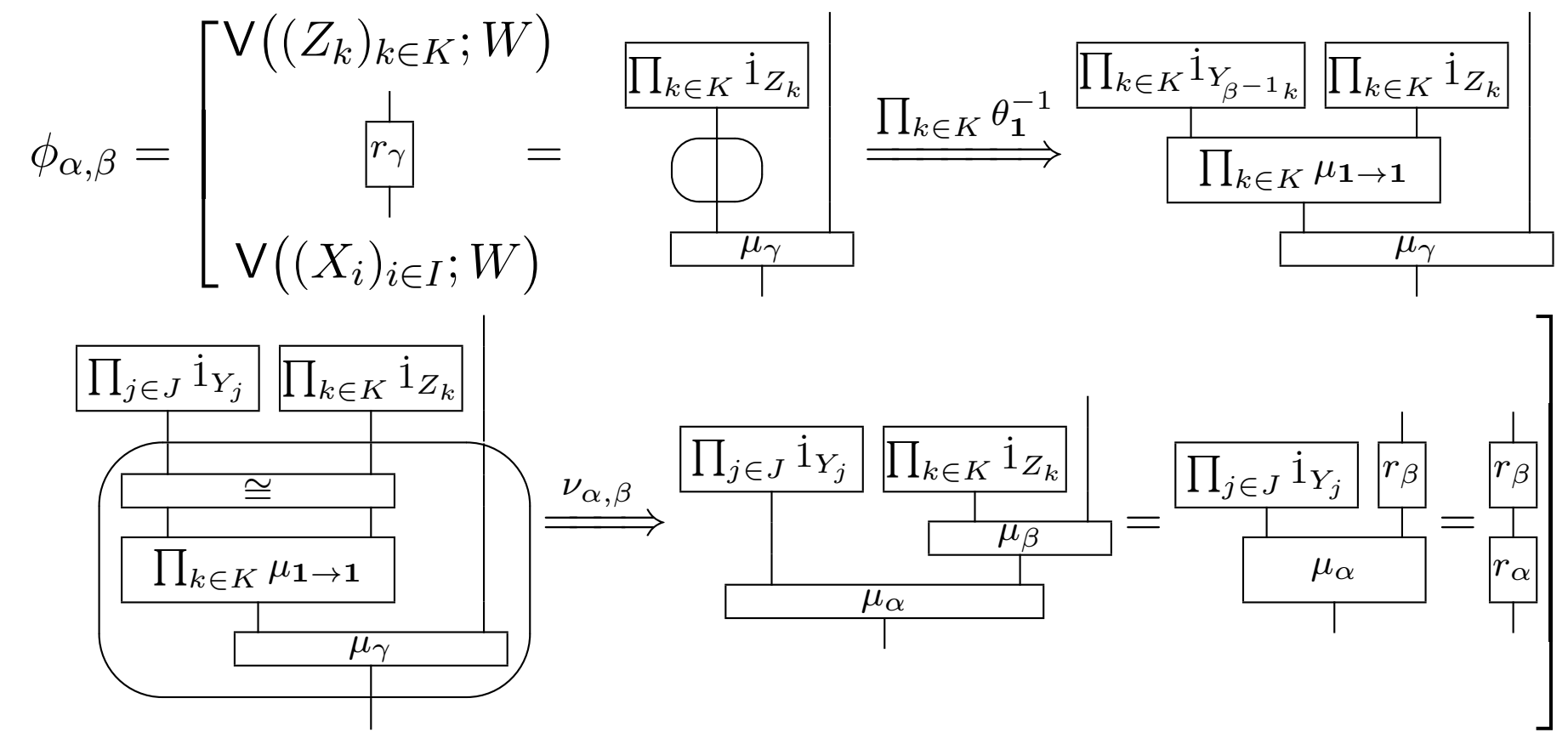

(2.1.3)

**2.1.5 Corollary.** *For all bijections* $\beta\colon J\to K\in\mathcal{S}_{\mathsf{sk}}$ *and* $Z_k=Y_{\beta^{-1}k}$ *the functor*

$$r_\beta\colon \mathsf{V}\big((Z_k)_{k\in K};W\big)\to\mathsf{V}\big((Y_j)_{j\in J};W\big)$$

*is an equivalence.*

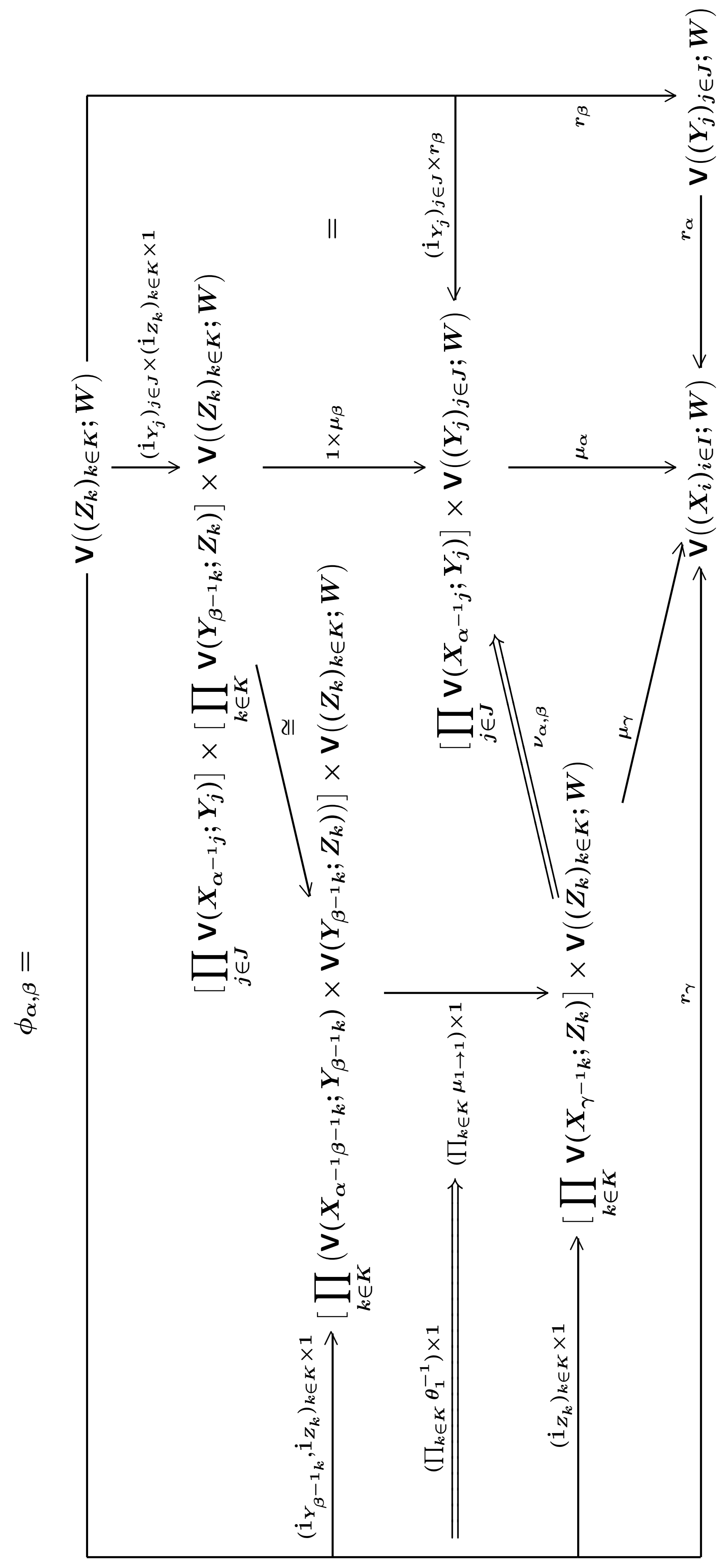

Figure 4: Cocycle for the action of symmetric groups on a symmetric weak multicategory

**2.1.6 Remark.** Each (symmetric or plain) weak multicategory $\mathsf{C}$ has an underlying bicategory (see Proposition 1.1.5). In fact, consider only index set $I = J = K = L = \mathbf{1}$ and the only map $\phi = \psi = \xi \colon \mathbf{1} \to \mathbf{1}$. So we consider the collection of objects $\operatorname{Ob}\mathsf{C}$, categories $\mathsf{C}(X;Y)$, composition functors $\mu = \mu_{\mathbf{1}\to\mathbf{1}} \colon \mathsf{C}(X;Y) \times \mathsf{C}(Y;Z) \to \mathsf{C}(X;Z)$, identity 1-morphisms $1_X \in \operatorname{Ob}\mathsf{C}(X;X)$, associator natural isomorphisms

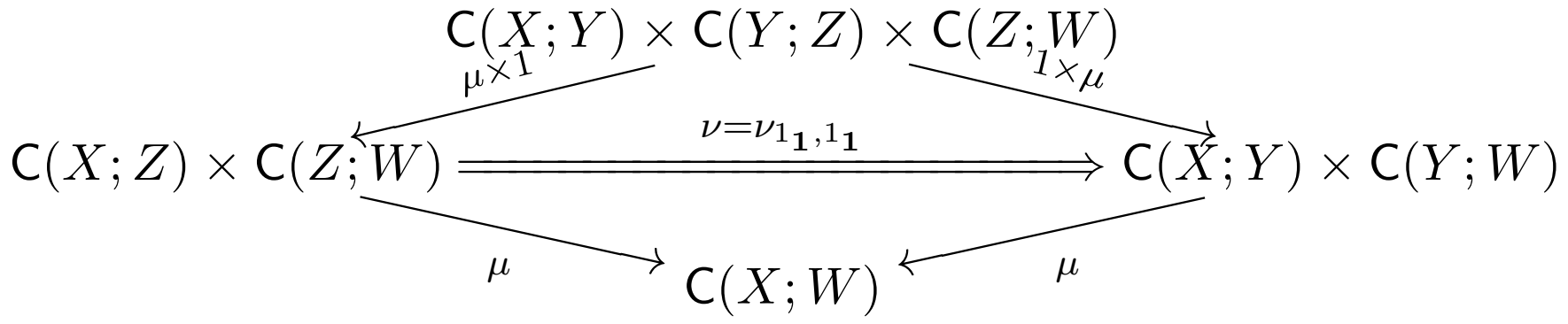

inverse right unitor natural isomorphisms

$$\zeta = \zeta_{\mathbf{1}} \colon \big[\mathsf{C}(X;Z) \xrightarrow{1\times \mathrm{i}_Z} \mathsf{C}(X;Z)\times\mathsf{C}(Z;Z) \xrightarrow{\mu} \mathsf{C}(X;Z)\big] \to \mathrm{Id},$$

inverse left unitor natural isomorphisms

$$\theta = \theta_{\mathbf{1}} \colon \big[\mathsf{C}(X;Z) \xrightarrow{\mathrm{i}_X\times 1} \mathsf{C}(X;X)\times\mathsf{C}(X;Z) \xrightarrow{\mu} \mathsf{C}(X;Z)\big] \to \mathrm{Id}$$

such that (1.1.2) and (1.1.4) hold (for $I = J = K = L = \mathbf{1}$).

**2.1.7 Proposition.** *Let $I \xrightarrow{\alpha} J \xrightarrow{\beta} K \xrightarrow{\gamma} L \in \mathcal{S}_{\mathsf{sk}}$ be bijections. Let $(X_i)_{i\in I}$, $(Y_j)_{j\in J}$, $(Z_k)_{k\in K}$, $(W_l)_{l\in L}$, $U$ be (families of) objects of $\mathsf{V}$ such that $X_i = Y_{\alpha i}$, $Y_j = Z_{\beta j}$, $Z_k = W_{\gamma k}$ for all $i \in I$, $j \in J$, $k \in K$. Then the non-abelian cocycle identity holds for $\phi$:*

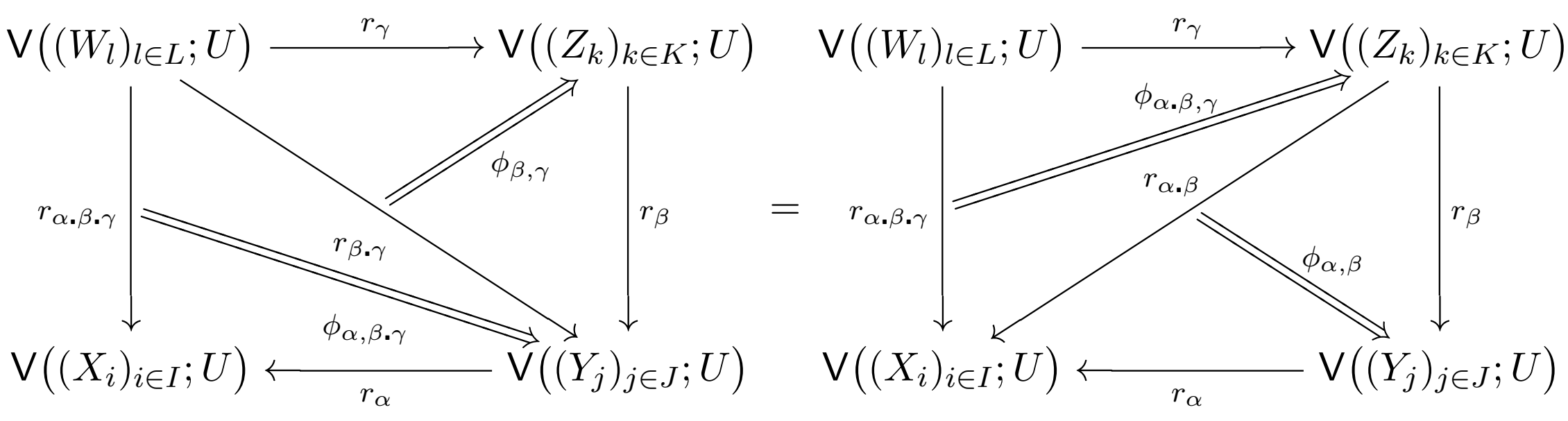

.

*Proof.* The left hand side is the composition

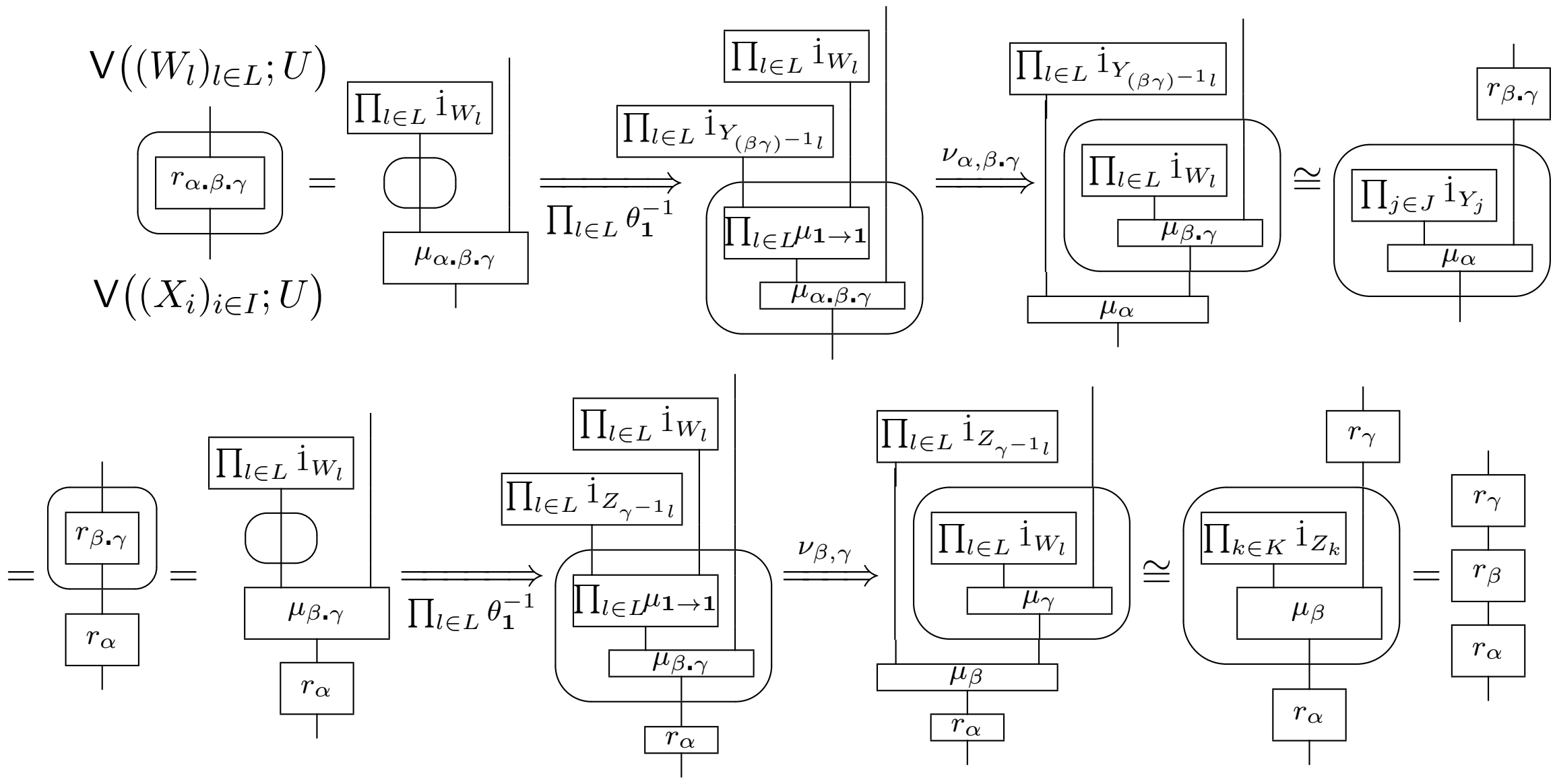

(2.1.4)

The right hand side is the composition

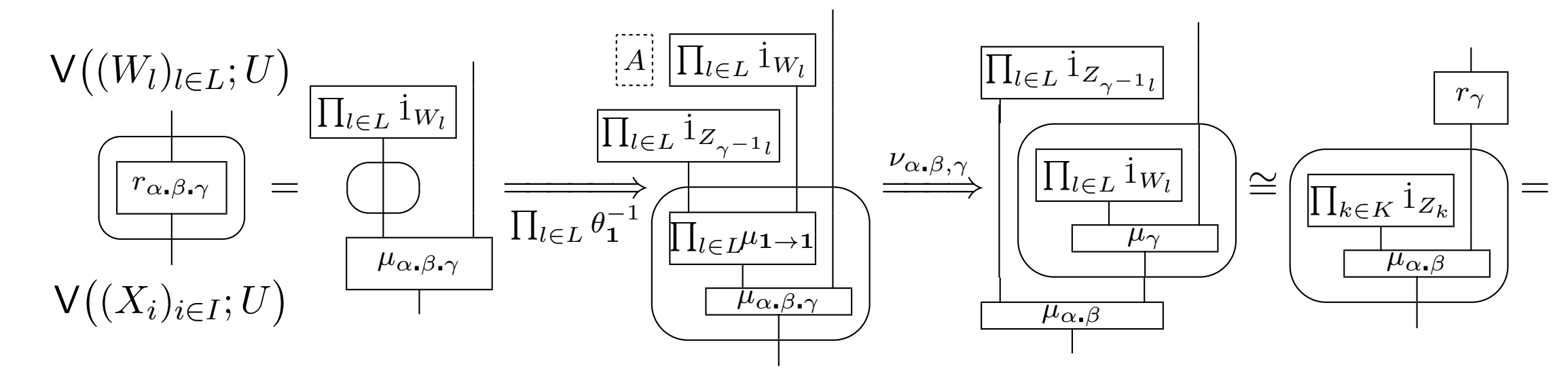

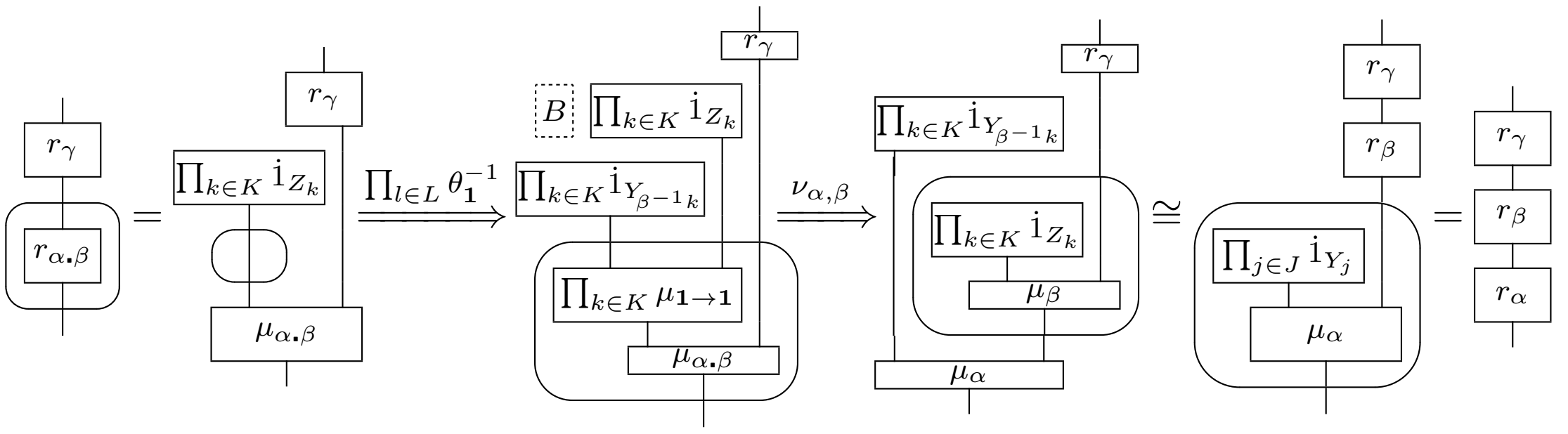

(2.1.5)

The left hand side can be transformed to

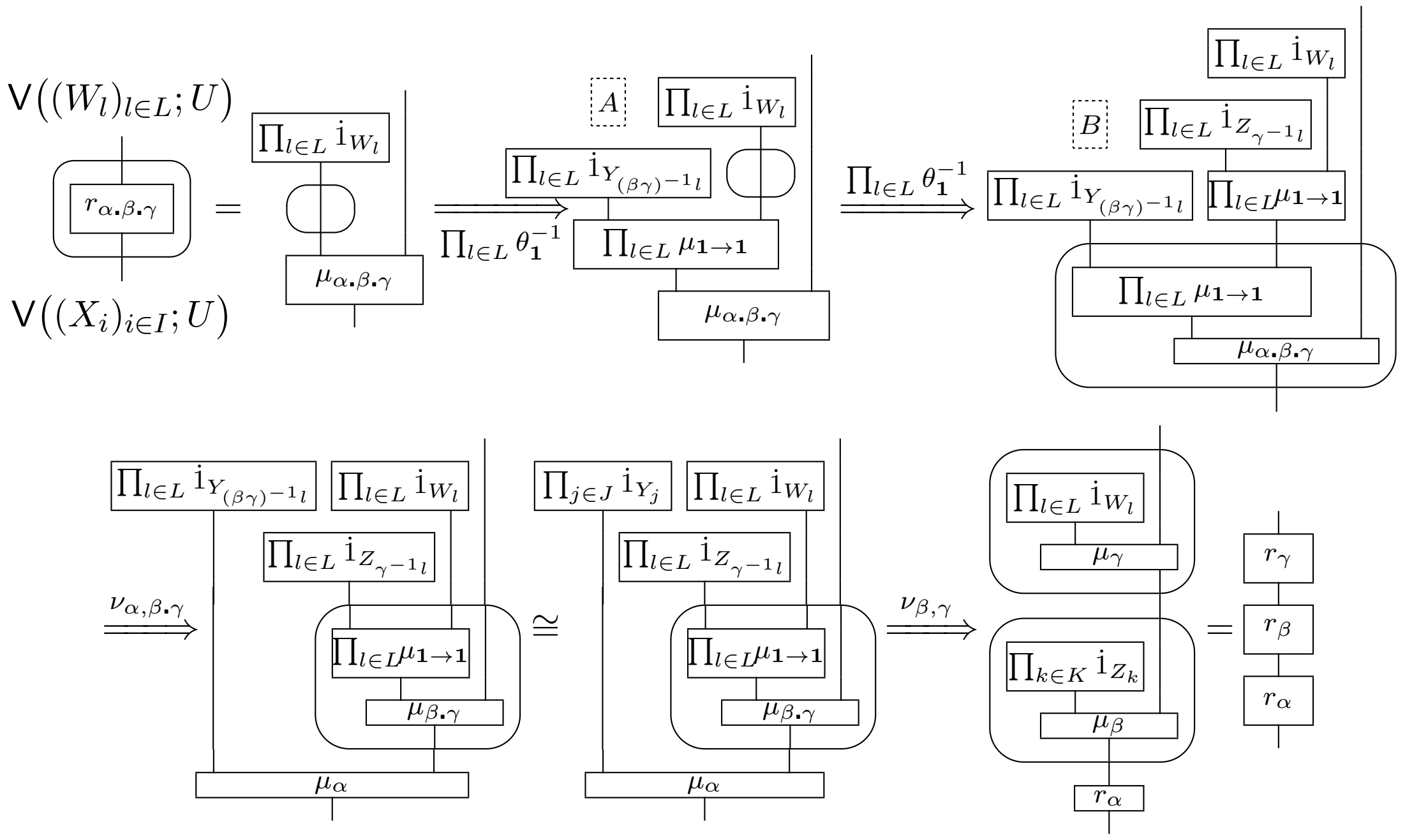

We may replace the above morphism $A \Rightarrow B$ with the following composition

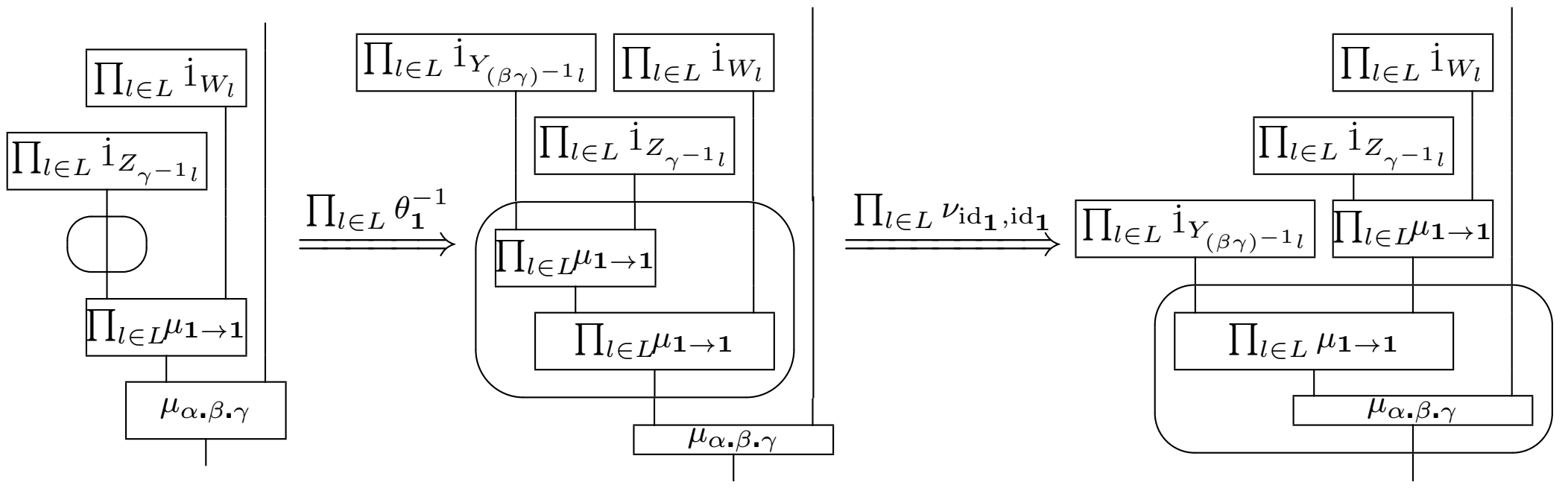

(2.1.6)

by the coherence (see Gurski [Gur13, Theorem 2.13]) for the bicategory $\mathsf{V}_1$ constructed in Proposition 1.1.5, for all involved arrows and 2-cells are structure elements of this bicategory.

Passage from (2.1.6) to

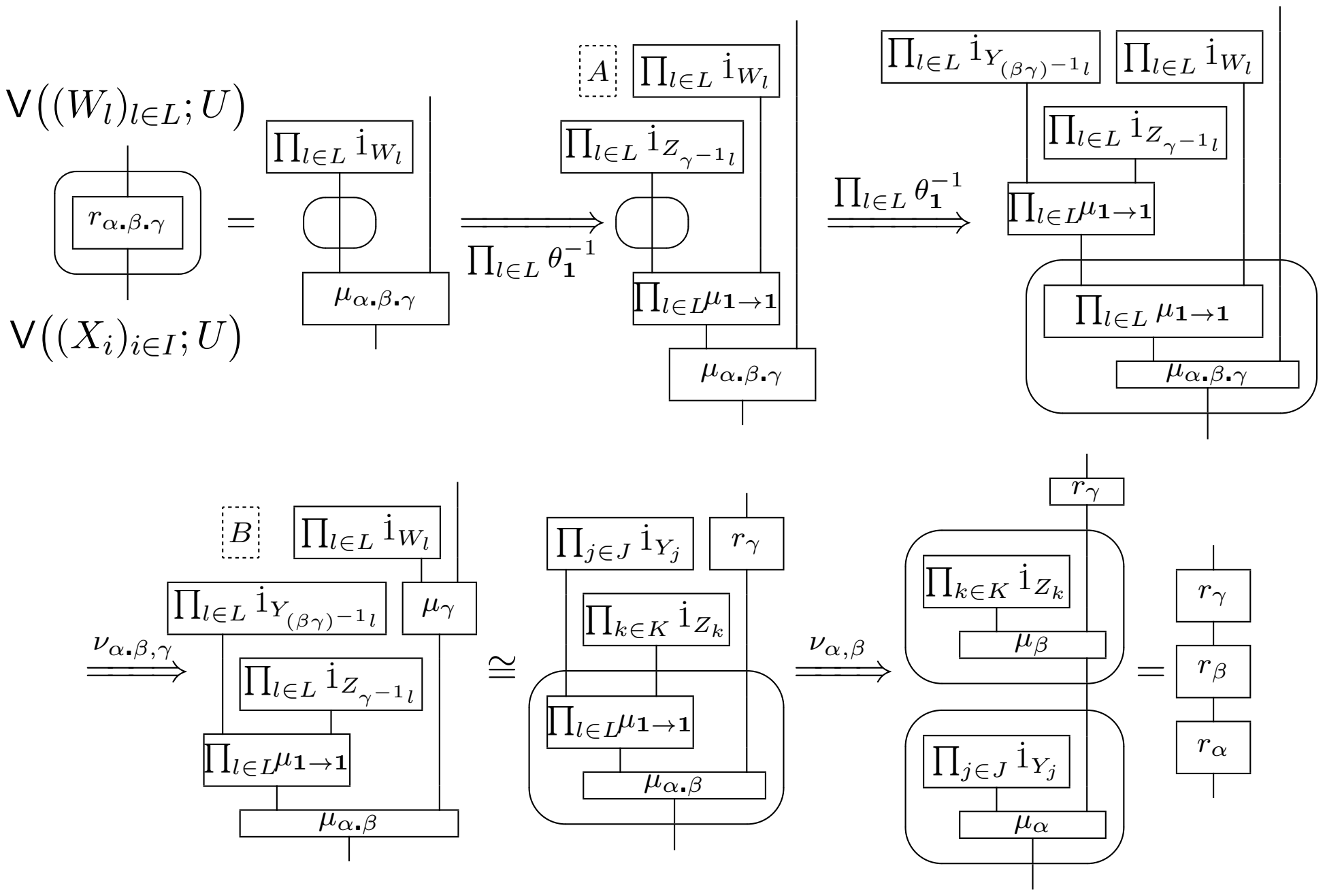

(2.1.7)

is ensured by (1.1.2).

Transformations (2.1.5) and (2.1.7) are equal because they have equal compositions $A \Rightarrow B$. Indeed, transformations $\nu_{\alpha.\beta,\gamma}$ and $\prod_{l\in L}\theta_{\mathbf{1}}^{-1}$ commute having non-intersecting supports. □

We have a natural isomorphism $\theta_I\colon r_{\mathrm{id}_I} \to \mathrm{id}$. Thus, we have a weak action of a symmetric group on the collection of homomorphism categories of a symmetric weak multicategory $\mathsf{V}$.

**2.1.8 Proposition.** *Let the square in $\mathcal{S}_{\mathsf{sk}}$, where vertical arrows are bijections,*

$$\begin{array}{ccc} I & \xrightarrow{\ \alpha\ } & J \\ {\scriptstyle\pi}\Big\downarrow{\scriptstyle\cong} & & {\scriptstyle\cong}\Big\downarrow{\scriptstyle\beta} \\ L & \xrightarrow{\ \gamma\ } & K \end{array}$$

*commute. Then there is the equivariance property*

$$\begin{aligned}
&\Big\{\Big[\prod_{k\in K}\mathsf{V}\big((X_{\pi^{-1}l})_{l\in\gamma^{-1}k};Y_{\beta^{-1}k}\big)\Big]\times\mathsf{V}\big((Y_{\beta^{-1}k})_{k\in K};W\big)\\
&\qquad\xrightarrow{\prod_{\beta^{-1}}\times 1}\Big[\prod_{j\in J}\mathsf{V}\big((X_{\pi^{-1}l})_{l\in\pi\alpha^{-1}j};Y_j\big)\Big]\times\mathsf{V}\big((Y_{\beta^{-1}k})_{k\in K};W\big)\\
&\qquad\xrightarrow{\prod_{j\in J}r_{\varpi_j}\times r_\beta}\Big[\prod_{j\in J}\mathsf{V}\big((X_i)_{i\in\alpha^{-1}j};Y_j\big)\Big]\times\mathsf{V}\big((Y_j)_{j\in J};W\big)\xrightarrow{\mu_\alpha}\mathsf{V}\big((X_i)_{i\in I};W\big)\Big\}\\
&\xrightarrow[\cong]{\nu_{\alpha,\beta}^{-1}\cdot\prod_{k\in K}\zeta_{\alpha^{-1}\beta^{-1}k}}\Big\{\Big[\prod_{k\in K}\mathsf{V}\big((X_{\pi^{-1}l})_{l\in\gamma^{-1}k};Y_{\beta^{-1}k}\big)\Big]\times\mathsf{V}\big((Y_{\beta^{-1}k})_{k\in K};W\big)\\
&\qquad\xrightarrow{\prod_{k\in K}r_{\pi_k}\times 1}\Big[\prod_{k\in K}\mathsf{V}\big((X_i)_{i\in\pi^{-1}\gamma^{-1}k};Y_{\beta^{-1}k}\big)\Big]\times\mathsf{V}\big((Y_{\beta^{-1}k})_{k\in K};W\big)\xrightarrow{\mu_{\pi.\gamma}}\mathsf{V}\big((X_i)_{i\in I};W\big)\Big\}\\
&\xrightarrow[\cong]{\nu_{\pi,\gamma}}\Big\{\Big[\prod_{k\in K}\mathsf{V}\big((X_{\pi^{-1}l})_{l\in\gamma^{-1}k};Y_{\beta^{-1}k}\big)\Big]\times\mathsf{V}\big((Y_{\beta^{-1}k})_{k\in K};W\big)\xrightarrow{\mu_\gamma}\mathsf{V}\big((X_{\pi^{-1}l})_{l\in L};W\big)\\
&\qquad\xrightarrow{r_\pi}\mathsf{V}\big((X_i)_{i\in I};W\big)\Big\}. \qquad (2.1.8)
\end{aligned}$$

*Here $\varpi_j = \pi|\colon \alpha^{-1}j \to \pi\alpha^{-1}j = \gamma^{-1}\beta j$ and $\pi_k = \varpi_{\beta^{-1}k} = \pi|\colon \pi^{-1}\gamma^{-1}k \to \gamma^{-1}k$ are bijections.*

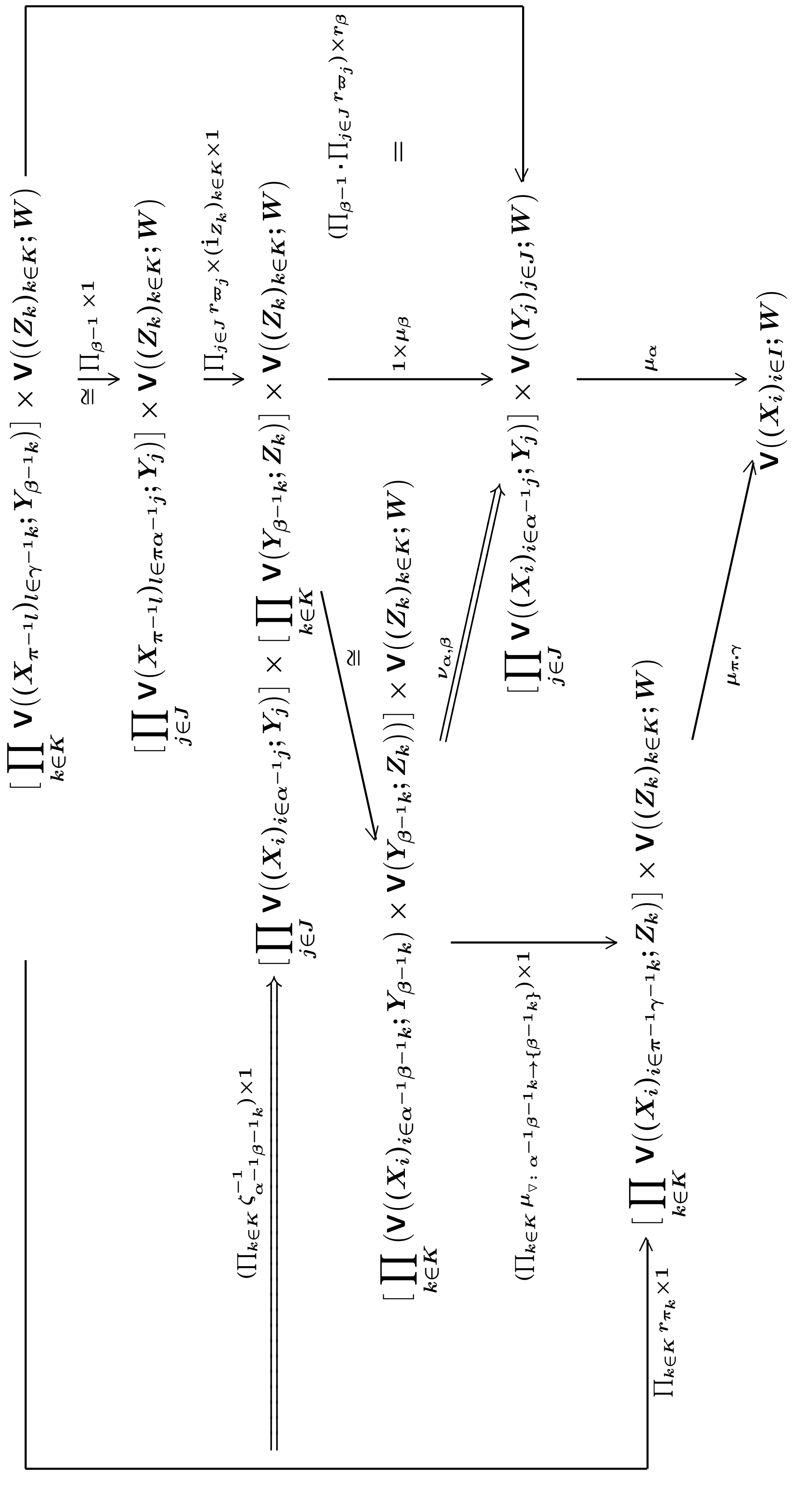

Figure 5: Equivariance of action of symmetric groups on a symmetric multicategory

*Proof.* Denote $Z_k = Y_{\beta^{-1}k}$. Using the associativity isomorphism from (1.1.1) for maps $I \xrightarrow{\alpha} J \xrightarrow{\beta} K$ we get the first isomorphism from (2.1.8) on the preceding page. In Bartlett's notation this isomorphism is the composition

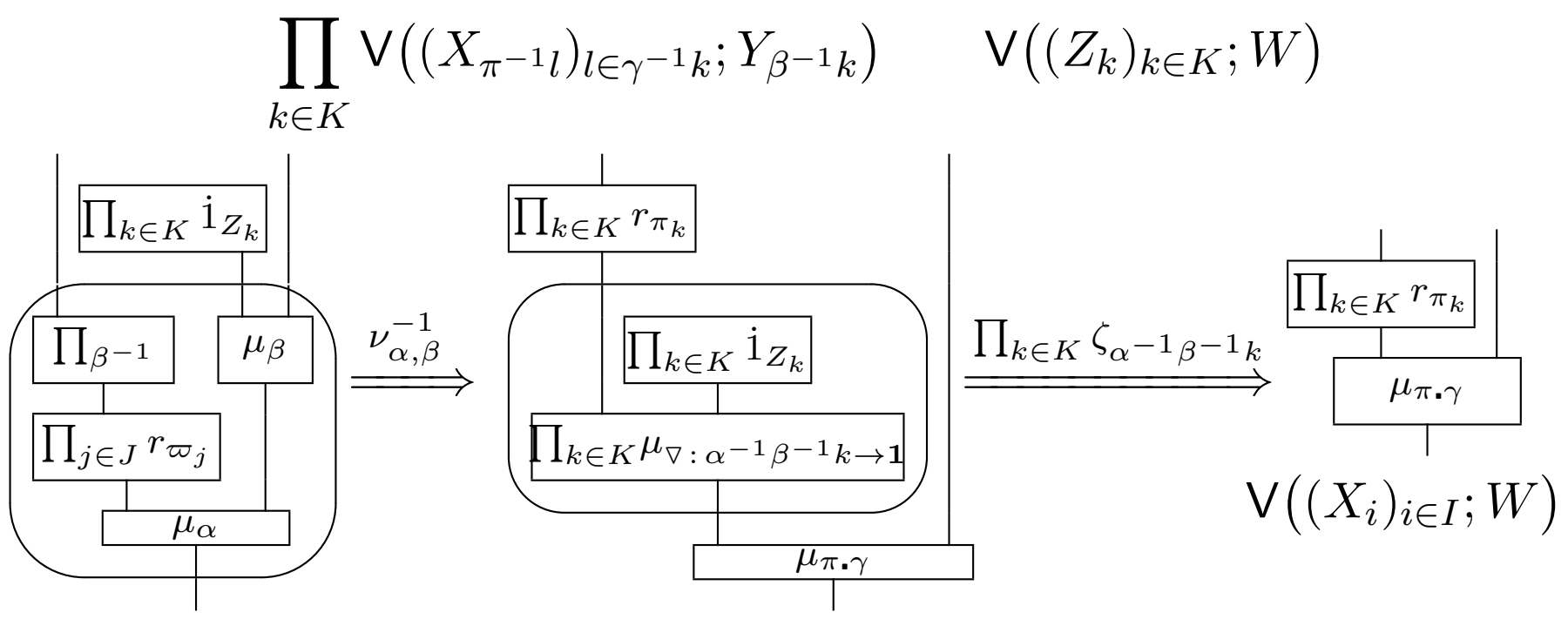

.

In order to obtain the second isomorphism from (2.1.8) we substitute into the former expression the definition of $r$:

$$\Big[\prod_{k\in K} \mathsf{V}\big((X_{\pi^{-1}l})_{l\in\gamma^{-1}k}; Z_k\big)\Big] \times \mathsf{V}\big((Z_k)_{k\in K}; W\big) \xrightarrow{\prod_{k\in K}[(\mathrm{i}_{X_{\pi^{-1}l}})_{l\in\gamma^{-1}k}\times 1]\times 1}$$

$$\prod_{k\in K}\Big[\prod_{l\in\gamma^{-1}k} \mathsf{V}(X_{\pi^{-1}l}; X_{\pi^{-1}l}) \times \mathsf{V}\big((X_{\pi^{-1}l})_{l\in\gamma^{-1}k}; Z_k\big)\Big] \times \mathsf{V}\big((Z_k)_{k\in K}; W\big)$$

$$\xrightarrow{\prod_{k\in K}\mu_{\pi_k}\times 1} \Big[\prod_{k\in K} \mathsf{V}\big((X_i)_{i\in\pi^{-1}\gamma^{-1}k}; Z_k\big)\Big] \times \mathsf{V}\big((Z_k)_{k\in K}; W\big) \xrightarrow{\mu_{\pi\centerdot\gamma}} \mathsf{V}\big((X_i)_{i\in I}; W\big).$$

Transforming this with the help of the associativity isomorphism from (1.1.1) for maps $I \xrightarrow{\pi} L \xrightarrow{\gamma} K$ we get

$$\Big\{\Big[\prod_{k\in K} \mathsf{V}\big((X_{\pi^{-1}l})_{l\in\gamma^{-1}k}; Z_k\big)\Big] \times \mathsf{V}\big((Z_k)_{k\in K}; W\big) \xrightarrow{(\mathrm{i}_{X_{\pi^{-1}l}})_{l\in L}\times 1\times 1}$$

$$\Big[\prod_{l\in L}\mathsf{V}(X_{\pi^{-1}l}; X_{\pi^{-1}l})\Big] \times \Big[\prod_{k\in K} \mathsf{V}\big((X_{\pi^{-1}l})_{l\in\gamma^{-1}k}; Z_k\big)\Big] \times \mathsf{V}\big((Z_k)_{k\in K}; W\big)$$

$$\xrightarrow{1\times\mu_\gamma} \Big[\prod_{l\in L}\mathsf{V}(X_{\pi^{-1}l}; X_{\pi^{-1}l})\Big] \times \mathsf{V}\big((X_{\pi^{-1}l})_{l\in L}; W\big) \xrightarrow{\mu_\pi} \mathsf{V}\big((X_i)_{i\in I}; W\big)\Big\}$$

$$= \Big\{\Big[\prod_{k\in K} \mathsf{V}\big((X_{\pi^{-1}l})_{l\in\gamma^{-1}k}; Z_k\big)\Big] \times \mathsf{V}\big((Z_k)_{k\in K}; W\big) \xrightarrow{\mu_\gamma} \mathsf{V}\big((X_{\pi^{-1}l})_{l\in L}; W\big)$$

$$\xrightarrow{(\mathrm{i}_{X_{\pi^{-1}l}})_{l\in L}\times 1} \Big[\prod_{l\in L}\mathsf{V}(X_{\pi^{-1}l}; X_{\pi^{-1}l})\Big] \times \mathsf{V}\big((X_{\pi^{-1}l})_{l\in L}; W\big) \xrightarrow{\mu_\pi} \mathsf{V}\big((X_i)_{i\in I}; W\big)\Big\}.$$

This is the last expression from (2.1.8) with expanded $r_\pi$. In graphical notations this isomorphism reads:

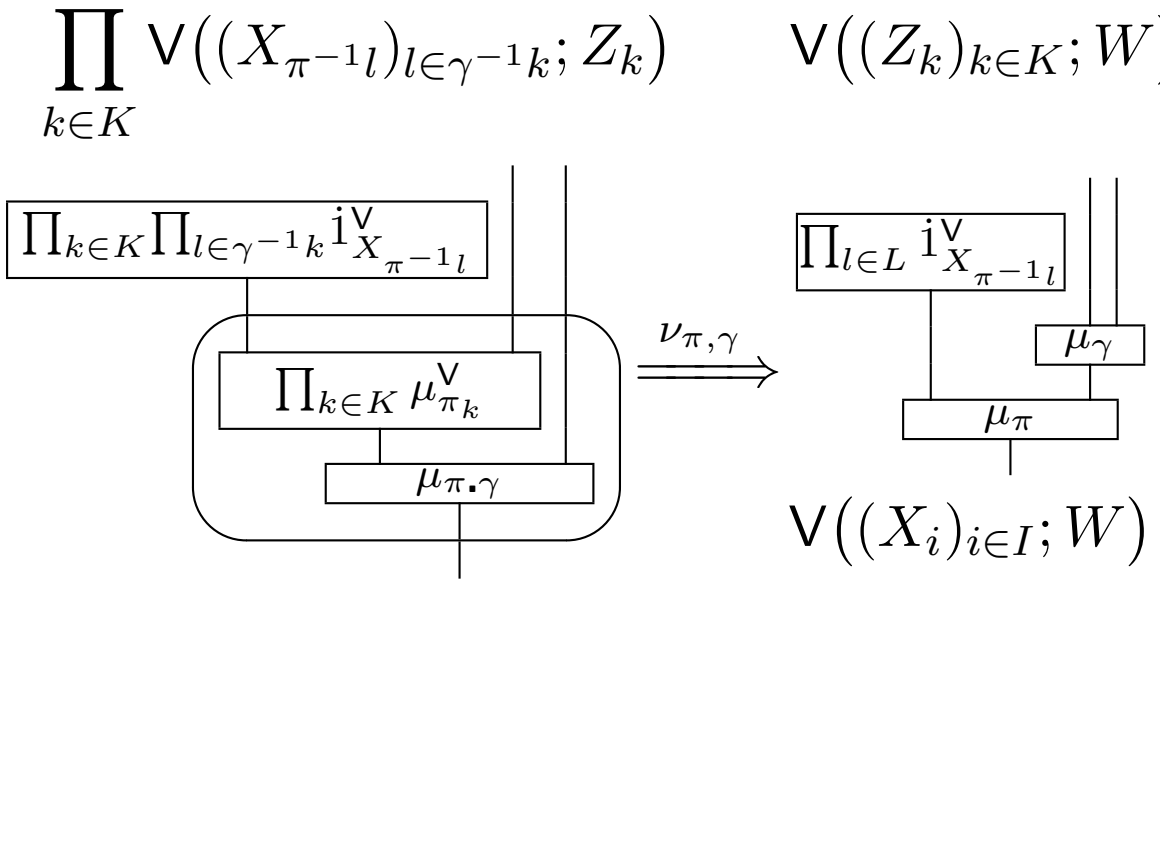

. (2.1.9)

□

## 2.2 Dependence on the order of totally ordered factors

For ordinary symmetric multicategories $\mathsf{V}$ we apply [Lyu23, (A.1.3)] to the square

$$\begin{array}{ccc} (I\times L,<'') & \xrightarrow{\mathrm{pr}_1} & I \\ {\scriptstyle\pi_{I,L}}\downarrow & & \| \\ (I\times L,<') & \xrightarrow{\mathrm{pr}_1} & I \end{array}$$

where $\pi_{I,L} = \mathrm{id}\colon (I\times L,<'') \to (I\times L,<') \in \mathcal{S}et$ is the identity map of the set $I\times L$, equipped with two different lexicographic orderings $(i_1,l_1) <'' (i_2,l_2)$ iff $l_1 < l_2$ or ($l_1 = l_2$ and $i_1 < i_2$). Another lexicographic ordering satisfies $(i_1,l_1) <' (i_2,l_2)$ iff $i_1 < i_2$ or ($i_1 = i_2$ and $l_1 < l_2$). The object $(I\times L,<') \in \mathrm{Ob}\,\mathcal{O}_{\mathsf{sk}}$ admits a different presentation $(L\times I,<'')$. For this presentation $\pi_{I,L}$ becomes the symmetry $\pi_{I,L} = (12)\colon (I\times L,<'') \to (L\times I,<'') \in \mathcal{S}et$. Then

$$\mu_{\mathrm{pr}_1\colon (I\times L,<'')\to I} = \big\{ \big[\prod_{i\in I} \mathsf{V}\big((X_{il})_{l\in L};Y_i\big)\big] \times \mathsf{V}\big((Y_i)_{i\in I};W\big) \xrightarrow{\mu_{\mathrm{pr}_1\colon (I\times L,<')\to I}} \mathsf{V}\big(((X_{il})_{l\in L})_{i\in I};W\big) \xrightarrow{r_{\pi_{I,L}}} \mathsf{V}\big(((X_{il})_{i\in I})_{l\in L};W\big) \big\}.$$

For symmetric weak multicategories $\mathsf{V}$ we have by Remark 2.1.1 and (2.1.8) the isomorphism

$$\mu_{\mathrm{pr}_1\colon (I\times L,<'')\to I} \xrightarrow{(\prod_{i\in I}\theta_L^{-1})\cdot\nu_{\pi_{I,L},\mathrm{pr}_1}} \big\{ \big[\prod_{i\in I} \mathsf{V}\big((X_{il})_{l\in L};Y_i\big)\big] \times \mathsf{V}\big((Y_i)_{i\in I};W\big) \xrightarrow{\mu_{\mathrm{pr}_1\colon (I\times L,<')\to I}} \mathsf{V}\big(((X_{il})_{l\in L})_{i\in I};W\big) \xrightarrow{r_{\pi_{I,L}}} \mathsf{V}\big(((X_{il})_{i\in I})_{l\in L};W\big) \big\}. \tag{2.2.1}$$

From here we get a natural isomorphism

$$\prod_{i\in I} \mathsf{V}\big((X_{il})_{l\in L};Y_i\big) \qquad \mathsf{V}\big((Y_i)_{i\in I};W\big) \tag{2.2.2}$$

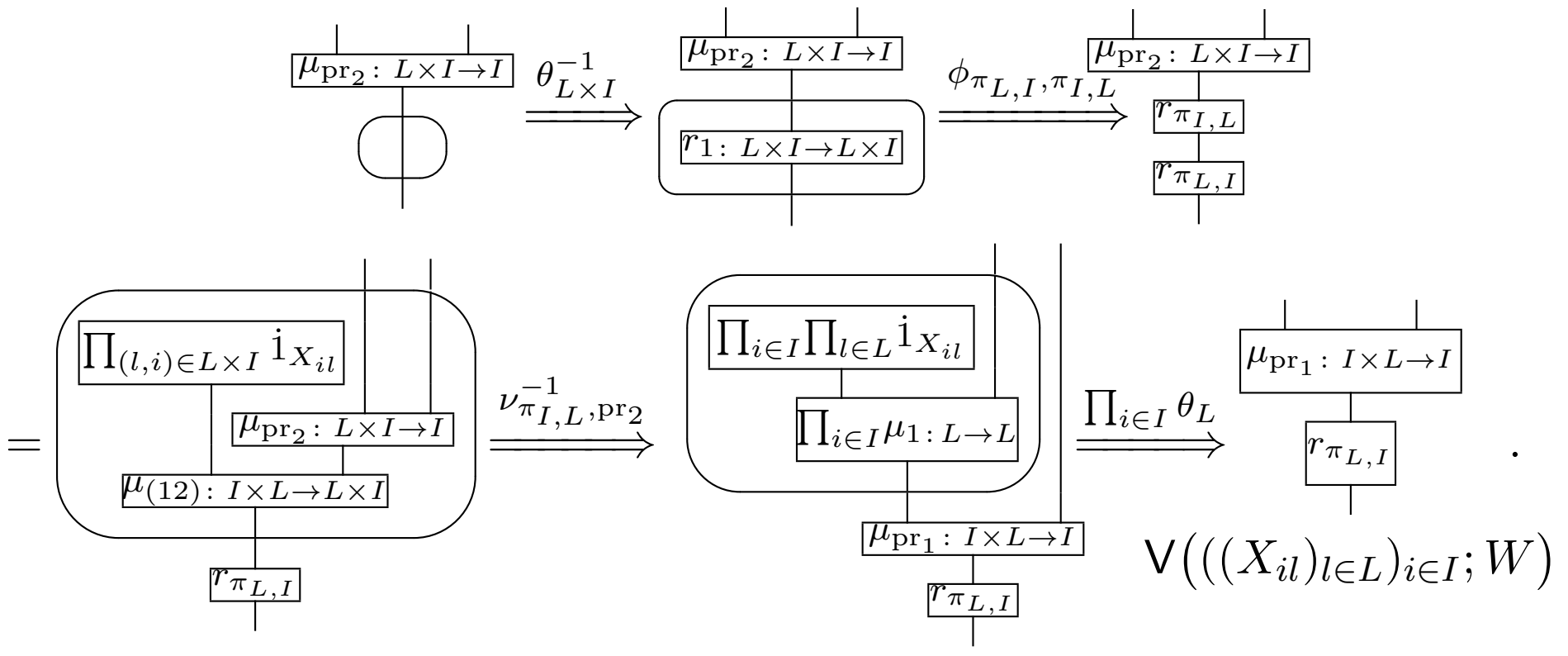

Here all products $S\times T$ of finite totally ordered sets are equipped with the lexicographic order $<''$.

Similarly we apply [Lyu23, (A.1.3)] for ordinary symmetric multicategories $\mathsf{V}$ for the square

$$\begin{array}{ccc} (I\times L,<'') & \xrightarrow{\mathrm{pr}_2} & L \\ {\scriptstyle\pi_{I,L}}\downarrow & & \| \\ (I\times L,<') & \xrightarrow{\mathrm{pr}_2} & L \end{array}$$

By [Lyu23, Proposition A.1.4] we conclude that

$$\mu_{\mathrm{pr}_2\colon (I\times L,<'')\to L} = \big\{ \big[\prod_{l\in L} \mathsf{V}\big((X_{il})_{i\in I};Y_l\big)\big] \times \mathsf{V}\big((Y_l)_{l\in L};W\big) \xrightarrow{\mu_{\mathrm{pr}_2\colon (I\times L,<')\to L}} \mathsf{V}\big(((X_{il})_{l\in L})_{i\in I};W\big) \xrightarrow{r_{\pi_{I,L}}} \mathsf{V}\big(((X_{il})_{i\in I})_{l\in L};W\big) \big\}.$$

For symmetric weak multicategories $\mathsf{V}$ we have by Remark 2.1.1 and (2.1.8) the isomorphism

$$\mu_{\mathrm{pr}_2\colon (I\times L,<'')\to L} \xrightarrow{(\prod_{l\in L}\theta_I^{-1})\cdot\nu_{\pi_{I,L},\mathrm{pr}_2}} \Big\{\big[\prod_{l\in L}\mathsf{V}\big((X_{il})_{i\in I};Y_l\big)\big]\times\mathsf{V}\big((Y_l)_{l\in L};W\big)$$

$$\xrightarrow{\mu_{\mathrm{pr}_2\colon (I\times L,<')\to L}} \mathsf{V}\big(((X_{il})_{l\in L})_{i\in I};W\big) \xrightarrow{r_{\pi_{I,L}}} \mathsf{V}\big(((X_{il})_{i\in I})_{l\in L};W\big)\Big\}. \quad (2.2.3)$$

From here we get a natural isomorphism

$$\prod_{l\in L}\mathsf{V}\big((X_{il})_{i\in I};Y_l\big) \qquad \mathsf{V}\big((Y_l)_{l\in L};W\big) \quad (2.2.4)$$

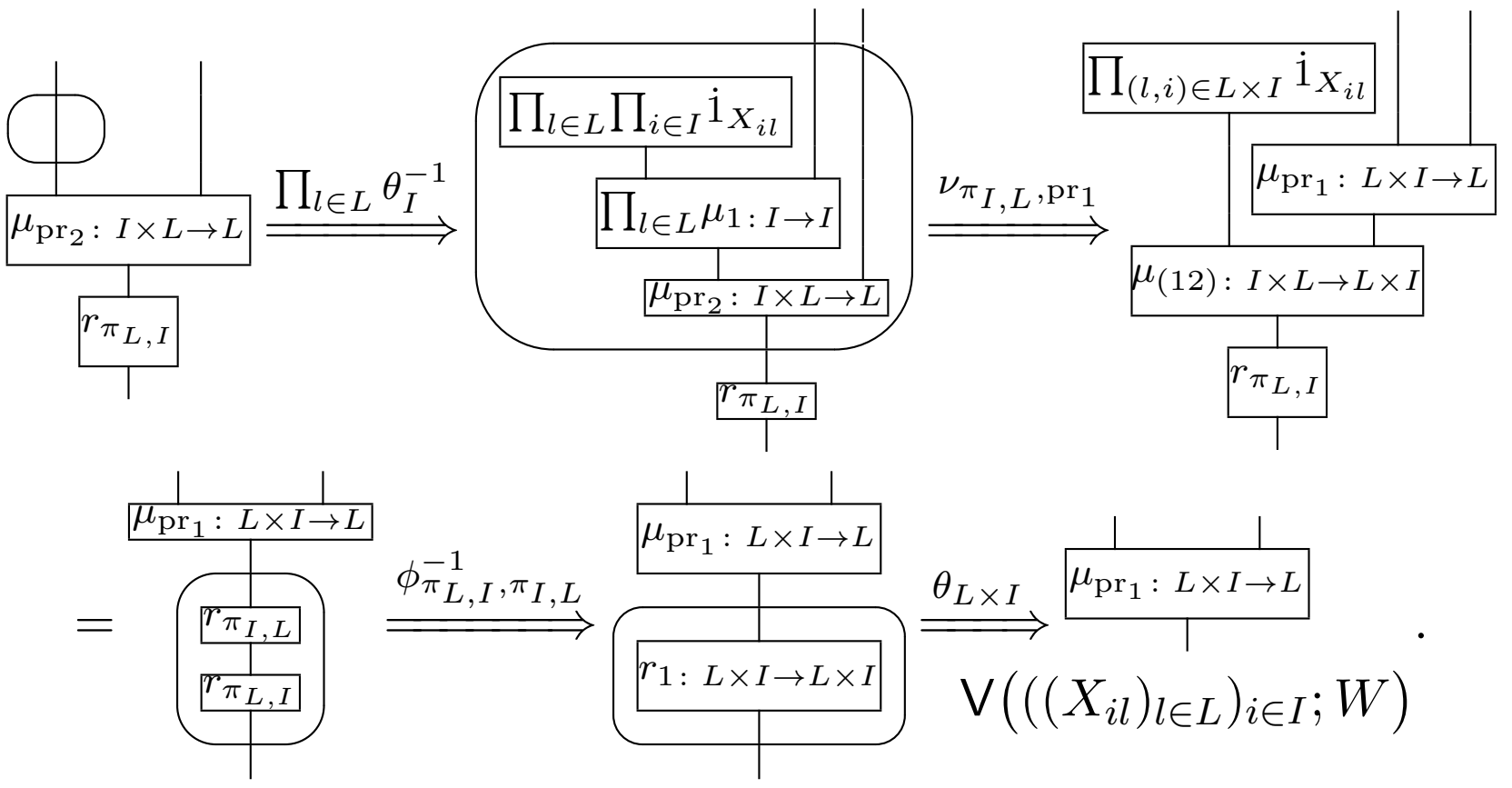

**2.2.1 Proposition.** *For any symmetric weak multicategory* $\mathsf{V}$ *we have an isomorphism of functors:*

$$\begin{array}{ccc}
\big[\prod_{(l,i)\in L\times I}\mathsf{V}\big((X_{ki})_{k\in\psi^{-1}l};Y_{li}\big)\big]\times\mathsf{V}\big(((Y_{li})_{l\in L})_{i\in I};Z\big) & \xrightarrow{\mu_{\psi\times I\colon K\times I\to L\times I}} & \mathsf{V}\big(((X_{ki})_{k\in K})_{i\in I};Z\big) \\
\Big\downarrow{\scriptstyle [\prod_{(l,i)\in L\times I} r_{1_{\psi^{-1}l}}]\times r_{\pi_{I,L}}} & \cong & \Big\downarrow{\scriptstyle r_{\pi_{I,K}}} \\
\big[\prod_{(i,l)\in I\times L}\mathsf{V}\big((X_{ki})_{k\in\psi^{-1}l};Y_{li}\big)\big]\times\mathsf{V}\big(((Y_{li})_{i\in I})_{l\in L};Z\big) & \xrightarrow{\mu_{I\times\psi\colon I\times K\to I\times L}} & \mathsf{V}\big(((X_{ki})_{i\in I})_{k\in K};Z\big)
\end{array}$$

*Proof.* This isomorphism is the composition of the following isomorphisms:

$$\prod_{(l,i)\in L\times I}\mathsf{V}\big((X_{ki})_{k\in\psi^{-1}l};Y_{li}\big) \qquad \mathsf{V}\big(((Y_{li})_{l\in L})_{i\in I};Z\big) \quad (2.2.5)$$

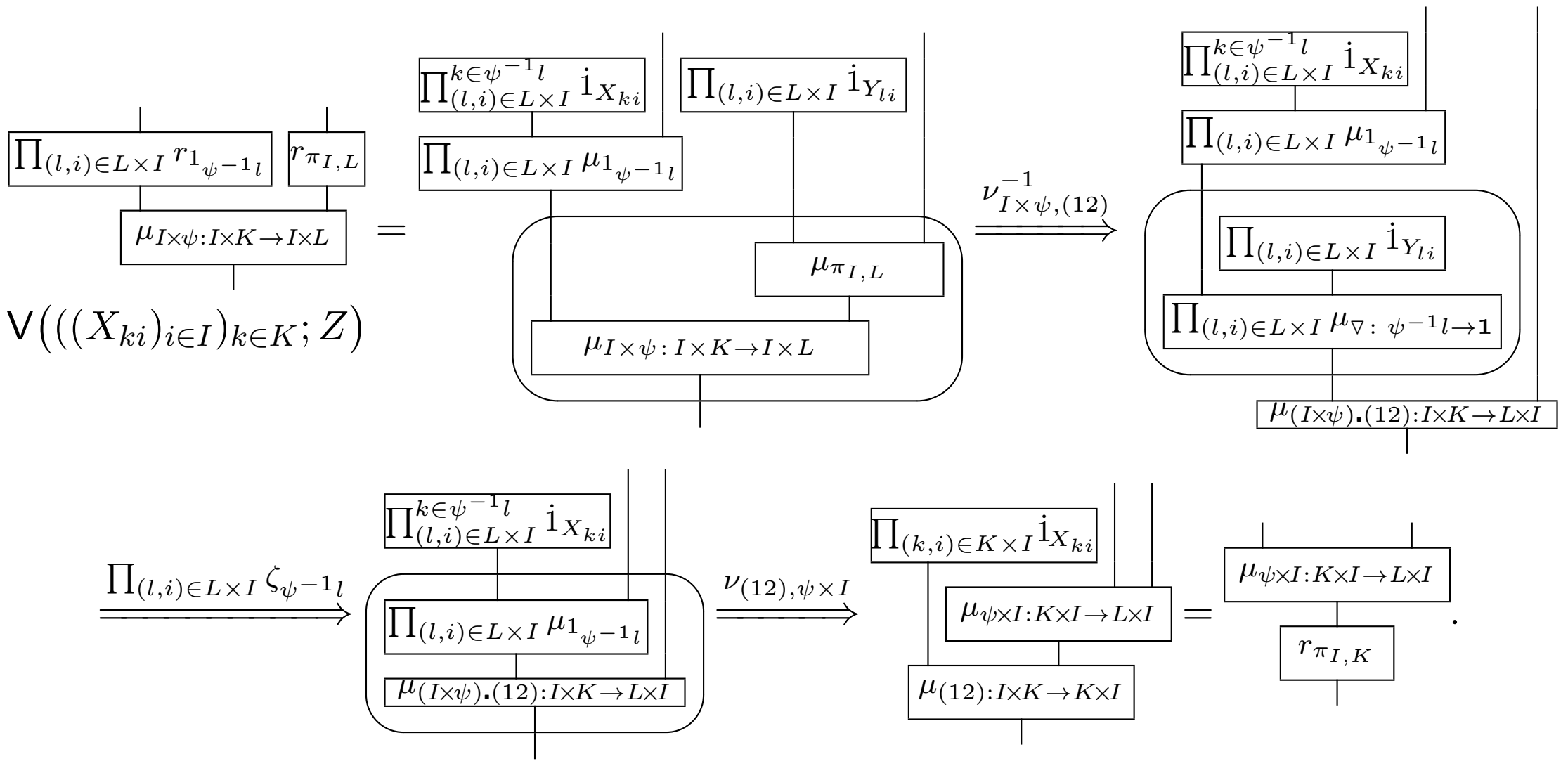

Here the first and the second $\nu$ refer to two decompositions of the same map:

$$\big(I\times K \xrightarrow{I\times\psi} I\times L \xrightarrow{(12)} L\times I\big) = \big(I\times K \xrightarrow{(12)} K\times I \xrightarrow{\psi\times I} L\times I\big).$$

□

# 3 Multinatural transformations and modifications

## 3.1 Multinatural transformations

**3.1.1 Definition.** Let $(F^l)_{l\in L}, G\colon \mathsf{C} \to \mathsf{D}$ be (a family of) symmetric weak multifunctors. A *multinatural transformation* $t\colon (F^l)_{l\in L} \to G\colon \mathsf{C} \to \mathsf{D}$ consists of

— objects $t_X \in \operatorname{Ob} \mathsf{D}((F^l X)_{l\in L}; GX)$;

— natural isomorphisms $t_{(X_i);Y}$ from

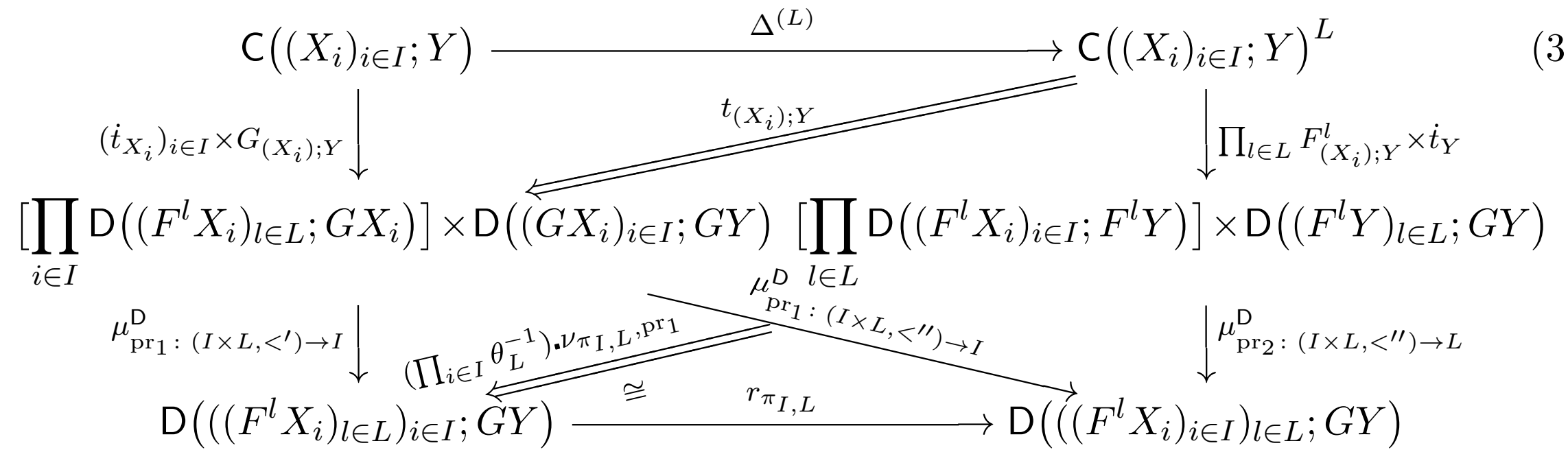

(3.1.1)

Here $\Delta^{(L)}\colon x \mapsto (x, x, \dots, x)$ is the diagonal functor. In Bartlett's notation we write

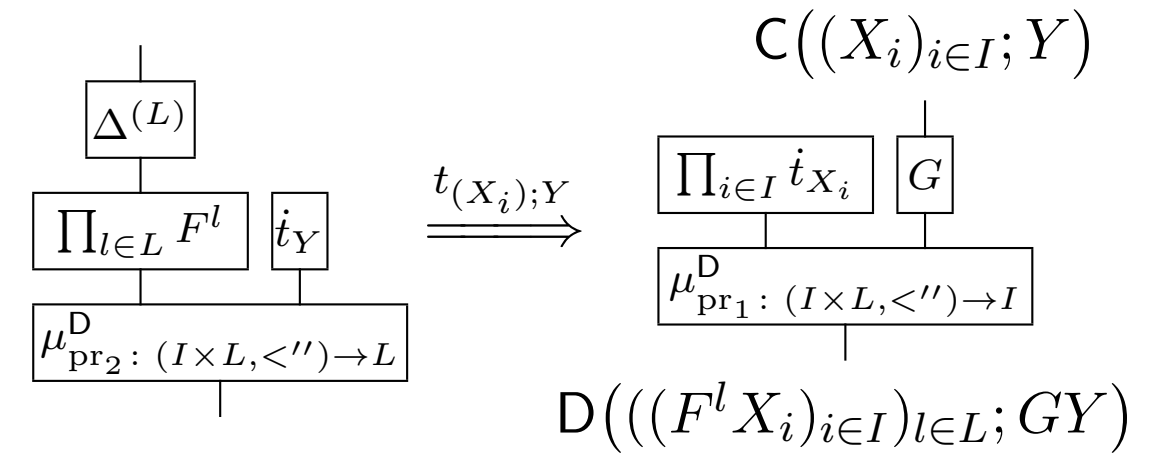

(3.1.2)

The following equations must hold for $\phi\colon I \to J \in \mathcal{S}_{\mathsf{sk}}$ and $X_i, Y_j, Z \in \operatorname{Ob} \mathsf{C}$

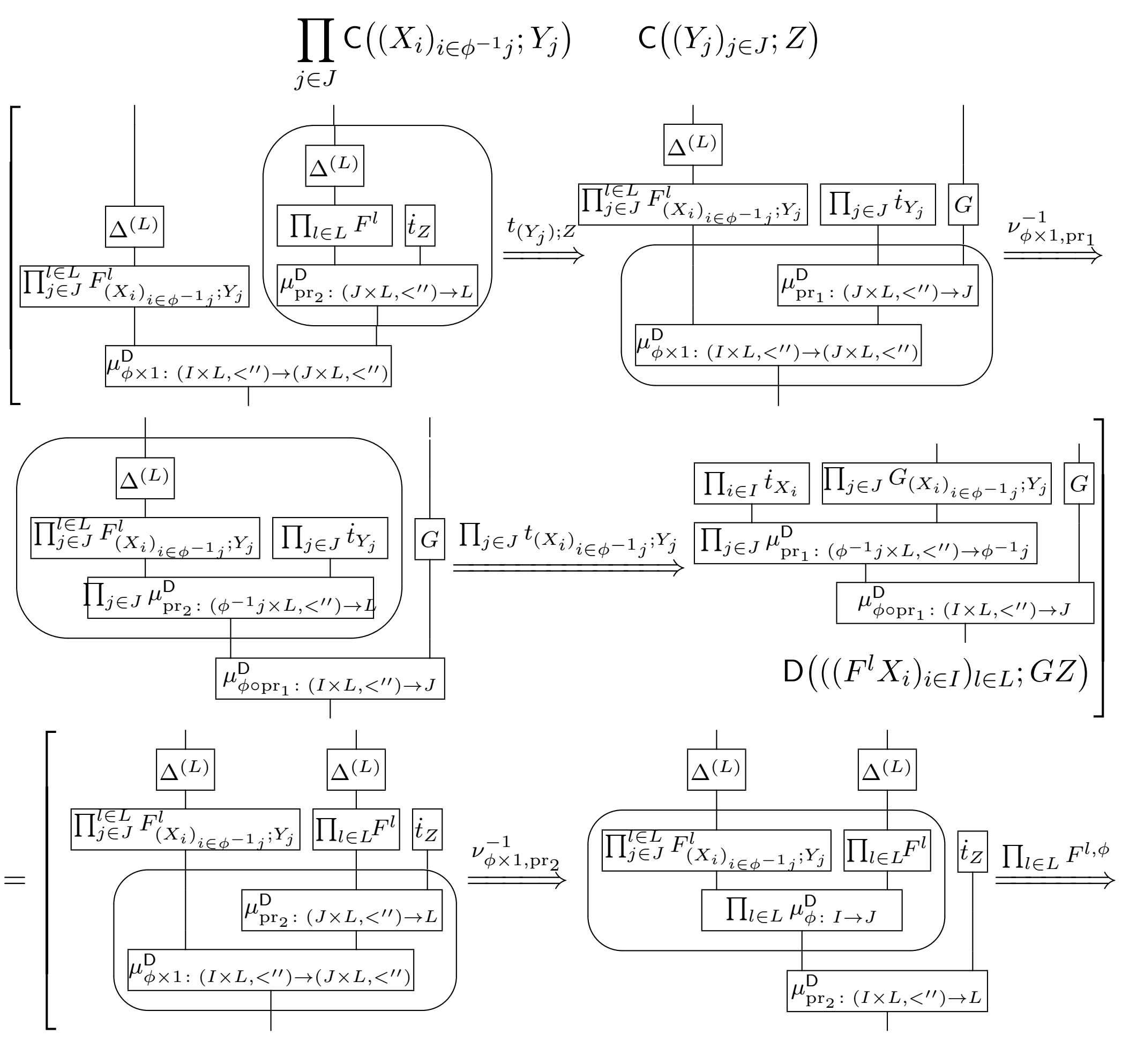

(3.1.3)

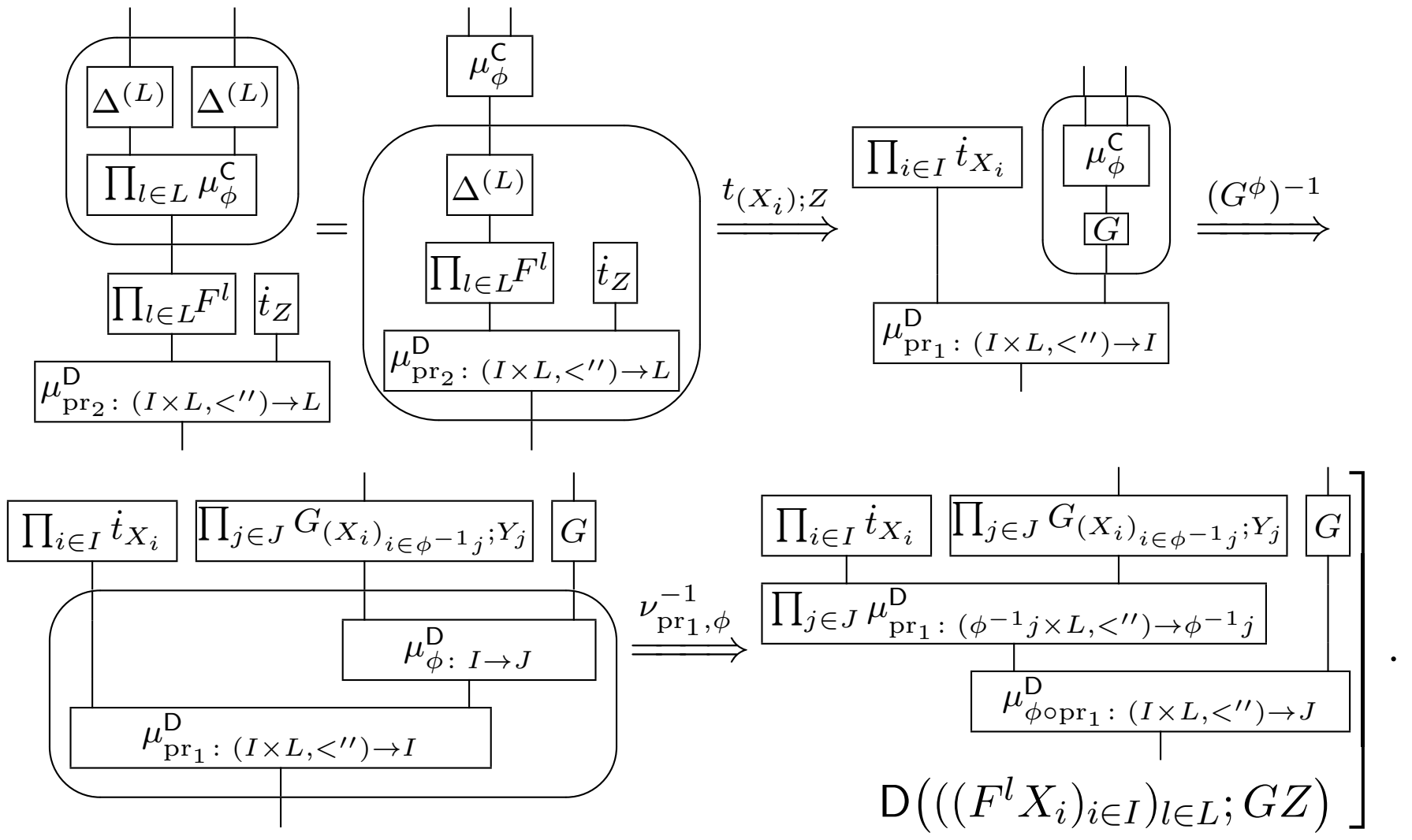

Another equation specifies the behaviour with respect to units:

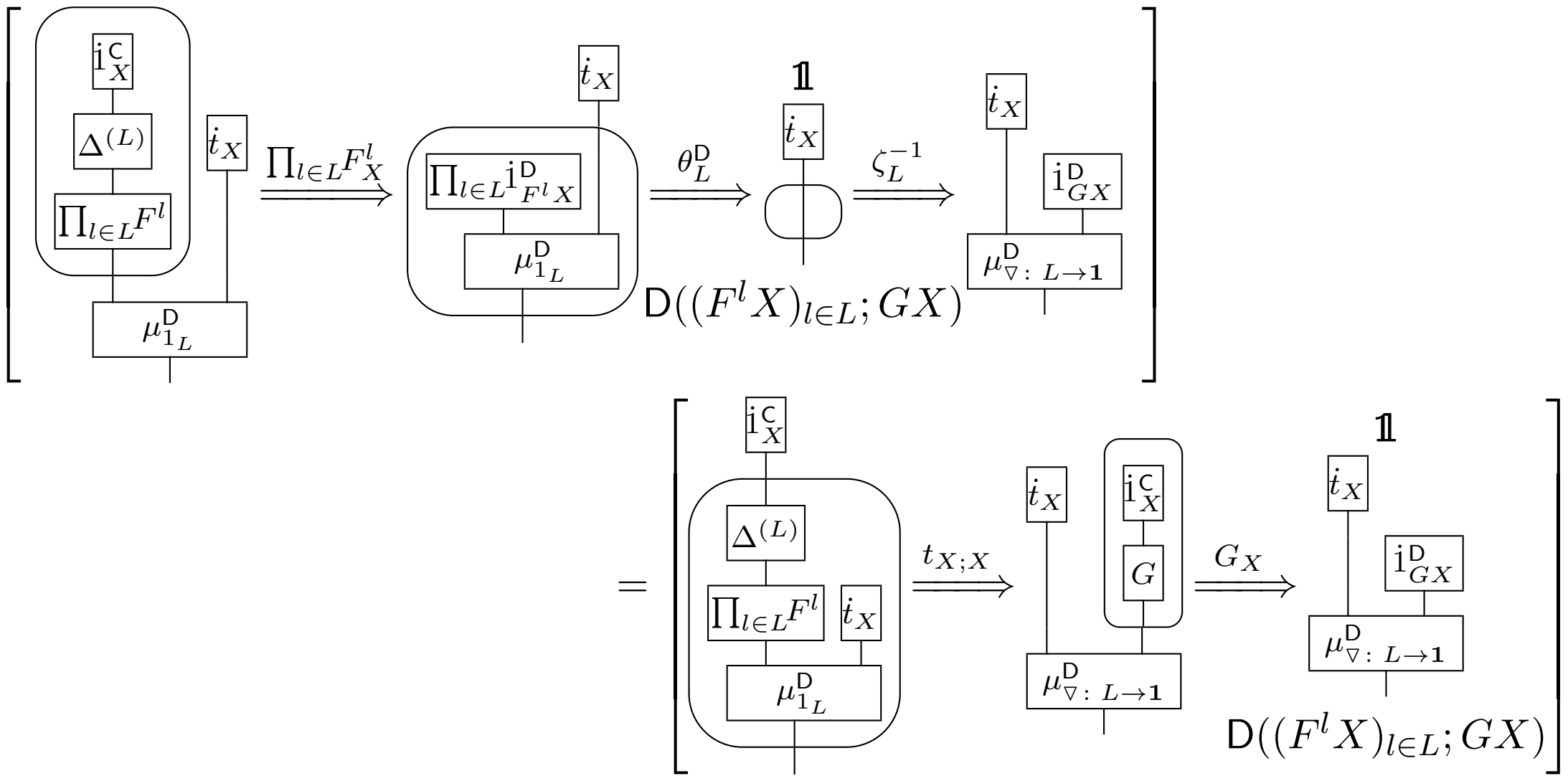

(3.1.4)

**3.1.2 Remark.** Objects $A$ of a category $\mathcal{C}$ are identified with functors $\dot{A} : \mathbb{1} \to \mathcal{C}$, $1 \mapsto A$. Morphisms $m\colon A \to B$ of a category $\mathcal{C}$ are identified with natural transformations $\dot{m}\colon \dot{A} \to \dot{B} : \mathbb{1} \to \mathcal{C}$, $\dot{m}_1 = m\colon A \to B$.

**3.1.3 Definition.** Let $t, p\colon (F^l)_{l\in L} \to G\colon \mathsf{C} \to \mathsf{D}$ be multinatural transformations of symmetric weak multifunctors. A *modification* $c\colon t \to p\colon (F^l)_{l\in L} \to G\colon \mathsf{C} \to \mathsf{D}$ is a family of

— morphisms $c_X \in \mathsf{D}\big((F^l X)_{l\in L}; GX\big)(t_X, p_X)$

such that the following square of natural transformations commutes:

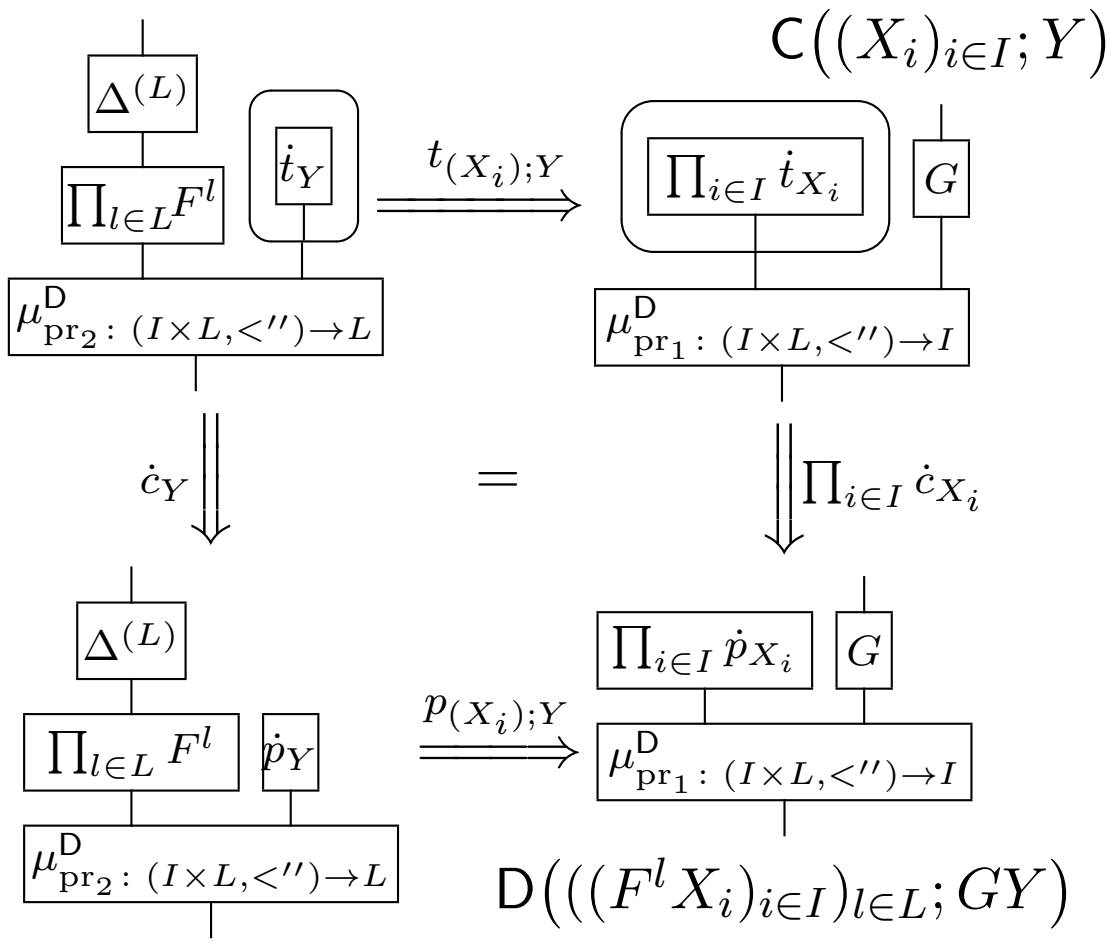

(3.1.5)

Comparing these notions with the definition of a monoidal bicategory we see that the latter requests more ingredients.

Institute of Mathematics, NAS Ukraine, 3 Tereshchenkivska st., Kyiv, 01024, Ukraine